\documentclass[12pt]{amsart}
\usepackage{amssymb}
\usepackage{amsthm}
\usepackage{amscd}

\usepackage{amsmath}
\usepackage{epic,eepic}
\usepackage{graphicx}
\DeclareMathAlphabet{\mathpzc}{OT1}{pzc}{m}{it}

\usepackage{mathrsfs}
\usepackage{latexsym}

\usepackage{pifont}

\theoremstyle{plain}
\newtheorem{lemma}{Lemma}[section]
\newtheorem{prop}[lemma]{Proposition}
\newtheorem{thm}[lemma]{Theorem}
\newtheorem{cor}[lemma]{Corollary}
\newtheorem{aplemma}{Lemma~A.\hspace{-1.5mm}}
\newtheorem{approp}{Proposition~A.\hspace{-1.5mm}}
\newtheorem{apthm}{Theorem~A.\hspace{-1.5mm}}
\newtheorem{apcor}{Corollary~A.\hspace{-1.5mm}}
\newtheorem{intthm}{Theorem}

\usepackage[all]{xy}
\usepackage {cancel}

\newcommand{\SSP}{\vspace{3mm}}
\newcommand{\LSP}{\vspace{5mm}}

\theoremstyle{definition}

\newtheorem{rema}[lemma]{Remark}

\newtheorem{remb}{Remark}

\newtheorem{defi}[lemma]{Definition}
\newtheorem{exa}[lemma]{Example}
\newtheorem{aprem}{Remark~A.\hspace{-1.5mm}}
\newtheorem{apdefi}{Definition~A.\hspace{-1.5mm}}
\newcommand{\bde}{\begin{defi}}
\newcommand{\ede}{\end{defi}\vspace{1mm}}
\newcommand{\ble}{\begin{lemma}}
\newcommand{\ele}{\end{lemma}}
\newcommand{\bpr}{\begin{prop}}
\newcommand{\epr}{\end{prop}}
\newcommand{\bt}{\begin{thm}}
\newcommand{\et}{\end{thm}}
\newcommand{\bco}{\begin{cor}}
\newcommand{\eco}{\end{cor}}
\newcommand{\bre}{\begin{rema}}
\newcommand{\ere}{\end{rema}}
\newcommand{\brea}{\begin{rema}}
\newcommand{\erea}{\end{rema}\vspace{1mm}}
\newcommand{\breb}{\begin{remb}}
\newcommand{\ereb}{\end{remb}\vspace{1mm}}
\newcommand{\bex}{\begin{exa}}
\newcommand{\eex}{\end{exa}}
\newcommand{\bpf}{\begin{proof}}
\newcommand{\epf}{\end{proof}\vspace{1mm}}

\newcommand{\bade}{\begin{apdefi}}
\newcommand{\eade}{\end{apdefi}}
\newcommand{\bale}{\begin{aplemma}}
\newcommand{\eale}{\end{aplemma}}
\newcommand{\bapr}{\begin{approp}}
\newcommand{\eapr}{\end{approp}}
\newcommand{\bat}{\begin{apthm}}
\newcommand{\eat}{\end{apthm}}
\newcommand{\baco}{\begin{apcor}}
\newcommand{\eaco}{\end{apcor}}
\newcommand{\bare}{\begin{aprem}}
\newcommand{\eare}{\end{aprem}}

\newcommand{\be}{\begin{enumerate}}
\newcommand{\ee}{\end{enumerate}}
\newcommand{\bcd}{\[\begin{CD}}
\newcommand{\ecd}{\end{CD}\]}
\newcommand{\bit}{\begin{itemize}}
\newcommand{\eit}{\end{itemize}}
\newcommand{\bq}{\begin{quote}}
\newcommand{\eq}{\end{quote}}
\newcommand{\ba}{\begin{array}}
\newcommand{\ea}{\end{array}}
\newcommand{\mcB}{\mathcal{B}}
\newcommand{\mcC}{\mathcal{C}}
\newcommand{\mcD}{\mathcal{D}}
\newcommand{\mcE}{\mathcal{E}}
\newcommand{\mcF}{\mathcal{F}}

\newcommand{\mcH}{\mathcal{H}}

\newcommand{\mcK}{\mathcal{K}}
\newcommand{\mcL}{\mathcal{L}}
\newcommand{\mcM}{\mathcal{M}}
\newcommand{\mcN}{\mathcal{N}}
\newcommand{\mcO}{\mathcal{O}}

\newcommand{\mcQ}{\mathcal{Q}}

\newcommand{\mcT}{\mathcal{T}}

\newcommand{\mcV}{\mathcal{V}}

\newcommand{\mbA}{\mathbb{A}}

\newcommand{\mbC}{\mathbb{C}}

\newcommand{\mbF}{\mathbb{F}}

\newcommand{\mbH}{\mathbb{H}}

\newcommand{\mbP}{\mathbb{P}}
\newcommand{\mbQ}{\mathbb{Q}}
\newcommand{\mbR}{\mathbb{R}}

\newcommand{\mbZ}{\mathbb{Z}}
\newcommand{\mfO}{\mathfrak{O}}

\newcommand{\mfb}{\mathfrak{b}}
\newcommand{\mfc}{\mathfrak{c}}

\newcommand{\mfg}{\mathfrak{g}}
\newcommand{\mfh}{\mathfrak{h}}

\newcommand{\mfl}{\mathfrak{l}}

\newcommand{\mfp}{\mathfrak{p}}
\newcommand{\mfq}{\mathfrak{q}}

\newcommand{\mfs}{\mathfrak{s}}
\newcommand{\mft}{\mathfrak{t}}

\newcommand{\msE}{\mathscr{E}}
\newcommand{\msF}{\mathscr{F}}

\newcommand{\msH}{\mathscr{H}}

\newcommand{\msU}{\mathscr{U}}

\newcommand{\msX}{\mathscr{X}}

\newcommand{\migi}{\rightarrow}
\newcommand{\longmigi}{\longrightarrow}

\newcommand{\isom}{\stackrel{\sim}{\migi}}

\newcommand{\migiincl}{\hookrightarrow}

\newcommand{\migisurj}{\twoheadrightarrow}

\newcommand{\mr}{\mathrm}
\newcommand{\hidden}[1]{\,}
\newcommand{\DE}{{^\dagger}\mcE}

\newcommand{\DV}{\mcV}

\newcommand{\qq}{\, q_1\hspace{-5.0mm}\rotatebox[origin=c]{-10}{$\swarrow$}\hspace{0.5mm}}

\newcommand{\C}{\mcC onn}

\begin{document}

\title[Infinitesimal deformations of opers  in positive characteristic]{Infinitesimal deformations of opers \\ in positive characteristic and \\ the de Rham cohomology of symmetric products}
\author{Yasuhiro Wakabayashi}
\date{}
\markboth{Yasuhiro Wakabayashi}{}
\maketitle
\footnotetext{Y. Wakabayashi: 
Graduate School of Information Science and Technology, Osaka University, Suita, Osaka 565-0871, Japan;}
\footnotetext{e-mail: {\tt wakabayashi@ist.osaka-u.ac.jp};}
\begin{abstract}
The Eichler-Shimura isomorphism describes a certain cohomology group with coefficients in a space of polynomials by using holomorphic modular/cusp forms. It determines a canonical decomposition of the corresponding de Rham cohomology group associated to a specific oper on a Riemann surface. One purpose of the present paper is to establish its analogue for opers in positive characteristic. We first discuss some basic properties on the (parabolic) de Rham cohomology groups and deformation spaces of $G$-opers (where $G$ is a semisimple algebraic group of adjoint type) in a general formulation. In particular, it is shown that the deformation space of a $G$-oper induced from an $\mathrm{SL}_2$-oper decomposes into a direct sum of the (parabolic) de Rham cohomology groups of its symmetric products. As a consequence, we obtain an Eichler-Shimura-type decomposition for dormant opers on general pointed stable curves by considering a transversal intersection of related spaces in the de Rham moduli space.
 
\end{abstract}
\tableofcontents 

%%%%%%%%%%%%%%%%%%%%%%%%%%%%%%%%%%%%%%%%%%%%%%%%%
%%%%%%%%%%%%%%%%%%%---[begin section]---%%%%%%%%%%%%%%
\section{Introduction}
\vspace{5mm}

%%%%%%%%%%%%%%%%%%%%%%%%%%%%%%%%%%%%%%%%%%%%%
%%%%%%%%%%%%%%---[begin section]---%%%%%%%%%%%%%%
\subsection{The Eichler-Shimura isomorphism for  a modular curve} \label{S01}

The Eichler-Shimura isomorphism describes a certain cohomology group with coefficients in a space of polynomials by using  holomorphic modular/cusp forms.
To be precise, let 
$N$ be an integer with $N \geq 5$, and 
$\Gamma$
the level $N$ congruence subgroup of $\mr{SL}_2 (\mbZ)$.
Given each $l \geq 2$,
we have  the $\mbC$-vector space
\begin{align}
V_{l-2} := \mr{Sym}^{l-2}(\mbC X + \mbC Y)
\end{align}
 of homogenous degree $l-2$ polynomials in $\mbC [X, Y]$.
It has a  $\Gamma$-action given by $\gamma (P) (X, Y) =P(a X + c Y, bX + d Y)$ for any $P \in V_{l-2}$ and $\gamma := \begin{pmatrix} a & b \\ c & d\end{pmatrix} \in \Gamma$.

The Eichler-Shimura isomorphism  (cf. ~\cite{Eic}, ~\cite{Shi}) asserts  a canonical isomorphism  for the cohomology group $H^1 (\Gamma, V_{l-2})$ of the resulting $\mbC [\Gamma]$-module $V_{l-2}$:
\begin{align} \label{Eq121}
M_l (\Gamma, \mbC) \oplus \overline{S_l (\Gamma, \mbC)} \isom H^1 (\Gamma, V_{l-2}),
\end{align} 
where $\overline{(-)}$ denotes the complex conjugation and 
$M_l (\Gamma, \mbC)$ (resp., $S_l (\Gamma, \mbC)$) denotes the $\mbC$-vector space of modular forms (resp., cusp forms) of weight $l$ and level $N$.

This isomorphism can be interpreted in geometric terms, as follows.
Let $Y_\Gamma$ denote  the Riemann surface defined as the quotient $\Gamma \backslash \mbH$ of 
the upper half plane $\mbH := \left\{z \in \mbC \, | \,  \mr{Im}(z) > 0\right\}$
  by the natural $\Gamma$-action.
Also, let $X_\Gamma := \Gamma \backslash \mbH^*$, where $\mbH^* := \mbH \cup \mbP^1 (\mbQ)$ (as a subset of $\mbP^1 (\mbC)$, equipped with a natural $\mr{SL}_2 (\mbR)$-action) has  a topology inducing 
a structure of  compact Riemann surface on the quotient $X_\Gamma$.
%The Riemann surface $Y_\Gamma^\mr{an}$ (resp., $X_\Gamma^{\mr{an}}$) comes from a smooth algebraic curve $Y_\Gamma$ (resp., $X_\Gamma$) over $\mbC$ via analytification.
The moduli interpretation of $Y_\Gamma$
brings  us 
a universal  
family of complex tori 
  $\pi : E \migi Y_\Gamma$ 
  of level $N$.

The relative de Rham cohomology  of this family forms
a holomorphic  vector bundle $\mcF_Y := \mbR^1 \pi_* \Omega_{E/Y_\Gamma}^\bullet$ of rank $2$; it carries a (flat) connection $\nabla_Y$ (called the {\it Gauss-Manin connection}), as well as  a $2$-step decreasing filtration $\{ \mcF_Y^j \}_{j=0}^2$ (called the {\it Hodge filtration}) given by $\mcF^0_Y := \mcF_Y$, $\mcF^{2}_Y := 0$, and $\mcF^1_Y  := \pi_* (\Omega_{E/Y_\Gamma})$, considered as a line subbundle of $\mcF_Y$ via the injection $\mcF^1_Y \migiincl \mcF_Y$ arising from the Hodge to de Rham spectral sequence  $'E_1^{a, b}:= \mbR^b \pi_* \Omega^a_{E/Y_\Gamma} \Rightarrow \mbR^{a + b}\pi_* \Omega_{E/Y_\Gamma}^\bullet $.

Each object in the data  $(\mcF_Y, \nabla_Y, \{ \mcF^j_Y \}_{j=0}^2)$ has  a natural extension  over  $X_\Gamma$
by admitting simple poles along the  cusps $D := X_\Gamma \setminus Y_\Gamma$;
%  defined as $X_\Gamma$ equipped with the log structure given by   the  cusps $D := X_\Gamma \setminus Y_\Gamma$; 
the resulting collection
 \begin{align} \label{Eq124}
 \msF_X^\heartsuit := (\mcF_X, \nabla_X, \{ \mcF_X^j \}_{j=0}^2)
 \end{align}
specifies a canonical  example of an $\mr{SL}_2$-oper 
 (cf. Remark \ref{Rem667}), and 
  the line bundle $\varTheta := \mcF_X^1$  satisfies  $\varTheta^{\otimes 2} \cong \Omega_{X} (D)$.

We shall write
$\mr{Sym}^{l-2}\msF_X$
for the $(l-2)$-nd symmetric power of  the flat bundle $\msF_X := (\mcF_X, \nabla_X)$.
The $1$-st   de Rham cohomology group 
$H_{\mr{dR}}^1 (X, \mr{Sym}^{l-2}\msF_X)$
 of 
 $\mr{Sym}^{l-2}\msF_X$
fits into the short exact sequence 
\begin{align} \label{Eq122}
0 \longmigi H^0 (X, \varTheta^{\otimes l})
\longmigi H_{\mr{dR}}^1 (X, \mr{Sym}^{l-2}\msF_X)
\longmigi H^0 (X, \varTheta^{\otimes l}(-D))^\vee
\longmigi 0
\end{align}
arising from the filtration on $\mr{Sym}^{l-2}\msF_X$ induced by   $\{ \mcF_X^j \}_j$.
 Since
  $\mr{Sym}^{l-2}\msF_X$ corresponds to the local system determined by the $\mbC [\Gamma]$-module  $V_{l-2}$ via the Riemann-Hilbert correspondence,
   there exists a canonical identification
  \begin{align}
  H^1_{\mr{dR}} (X, \mr{Sym}^{l-2} \msF_X) = 
  H^1 (\Gamma, V_{l-2})
  \end{align}
  (cf. ~\cite[Chap.\,II, Corollaire  6.10]{D}).
  On the other hand, 
   the elements of  $H^0 (X, \varTheta^{\otimes l})$ can be interpreted, via descent along the covering map $\mbH \migisurj Y_\Gamma$, as those of   $M_l (\Gamma, \mbC)$,
  and moreover,  the Petersson inner product yields 
   $H^0 (X, \varTheta^{\otimes l} (-D))^\vee =  \overline{S_l (\Gamma, \mbC)}$.
%   we have the 
% natural identifications
% \begin{align}
%H^1_{\mr{dR}} (X, \mr{Sym}^{l-2} \msF_X) = 
 % H^1 (\Gamma, V_{l-2}),
 % \hspace{3mm}
 % H^0 (X, \varTheta^{\otimes l}) = M_l (\Gamma, \mbC),
 % \hspace{3mm}
 % H^0 (X, \varTheta^{\otimes l} (-D))^\vee =  \overline{S_l (\Gamma, \mbC)},
 %\end{align}
%where the first ``$=$" follows from    
%and the third ``$=$" comes from the Petersson inner product.

Under these identifications,  the  Eichler-Shimura isomorphism (\ref{Eq121}) may be regarded as a canonical choice of  a   splitting of the short exact sequence (\ref{Eq122}) given by taking  periods.
See, e.g., ~\cite{BrHa}, ~\cite{KaSc}, ~\cite{Sch}  for algebraic treatments of the Eichler-Shimura isomorphism in terms of  de Rham cohomology.

%%%%%%%%%%%%%%%%%%%%%%%%%%%%%%%%%%%%%%%%%%%%%
%%%%%%%%%%%%%%---[begin section]---%%%%%%%%%%%%%%
\subsection{Analogous decomposition in positive characteristic} \label{S034}

The present paper aims  to 
investigate the short exact sequence (\ref{Eq122}) for opers in positive characteristic (including the case of  parabolic de Rham cohomology) in  relation to their  deformation spaces, as well as 
to establish an analogue of the Eichler-Shimura type decomposition for such opers. 
The present paper contains many arguments based on  those  in ~\cite{Wak8} and can be positioned  as a continuation of that work.

 Let  $p$ be an odd prime, $l$ an integer with $3 \leq l \leq p-1$,   $k$  an algebraically closed field of  characteristic $p$,  
 $(g, r)$ a pair of nonnegative integers with $2g-2 +r > 0$,
 and $\msX := (X, \{ \sigma_i \}_{i=1}^r)$
 an $r$-pointed stable curve over $k$ of genus $g$, where $X$ denotes a nodal curve over $k$ and $\{ \sigma_i \}_i$ is a set of marked points on $X$.
 (The present  paper deals not only with smooth curves but also with {\it pointed} and {\it  stable} curves, since the same arguments can be applied.
 However, it is not so essential to consider such a fairly general situation.)

Suppose that we are given  an  $\mr{SL}_2$-oper  $\msF^\heartsuit := (\mcF, \nabla, \{ \mcF^j \}_j)$ on $\msX$.
Then, the $2$-fold tensor product of the line bundle  $\varTheta := \mcF^1 \left(\subseteq \mcF\right)$ is isomorphic to
the sheaf  $\Omega_{X^\mr{log}}$ of logarithmic $1$-forms on the log curve $X^\mr{log}$ determined by  $\msX$.
%, which is a line bundle on $X$ with  $\varTheta^{\otimes 2} \cong \Omega_{X^\mr{log}}$.
Just as in the complex case, 
$\msF^\heartsuit$ associates 
the $1$-st  de Rham cohomology $H^1_{\mr{dR}} (X, \mr{Sym}^{l-2}\msF)$ of the $(l-2)$-nd symmetric power $\mr{Sym}^{l-2}\msF$ of  $\msF := (\mcF, \nabla)$, fitting into the  
 short exact sequence
\begin{align} \label{Eq120}
0 \longmigi H^0 (X, \varTheta^{\otimes l})
\longmigi  H^1_{\mr{dR}} (X, \mr{Sym}^{l-2}\msF)
\longmigi H^0 (X, \varTheta^{\otimes l}(-D))^\vee
\longmigi 0
\end{align}
(cf. Theorem \ref{Th67}, (i)), where $D$ denotes the divisor on $X$ given by the union of the $\sigma_i$'s.
(Similarly,  a version for the parabolic de Rham cohomology of $\mr{Sym}^{l-1}\msF$ holds, as asserted  in Theorem \ref{Th67}, (ii).)

Our question  regarding this sequence is described as follows:
\begin{quote}
{\it For what kind of $\mr{SL}_2$-opers $\msF^\heartsuit$ can we  construct a {\bf canonical} splitting of (\ref{Eq120}) arising from geometry  such as  the classical Eichler-Shimura isomorphism?}
\end{quote}
Even when
$\msF^\heartsuit$ is taken to be  the mod $p$ reduction of the ``$\msF_X^\heartsuit$" defined on a modular curve as above,
 the classical  Eichler-Shimura isomorphism would not be appropriate  in constructing a canonical splitting 
  because of its non-algebraic construction.

To get an answer to the above question, we discuss  
opers with vanishing $p$-curvature,
  called {\it dormant opers}.
Here,  recall that $p$-curvature  is an invariant for  (flat) connections measuring  the obstruction to the compatibility of $p$-power structures  appearing in certain associated spaces of infinitesimal symmetries;  this is a key ingredient of the $p$-adic Teichmuller theory developed by S. Mochizuki (cf. ~\cite{Mzk1}, ~\cite{Mzk2}), in which various moduli spaces of $\mr{SL}_2$-opers (or $\mr{PGL}_2$-opers) characterized by the degree of vanishing $p$-curvature  were investigated.
That work suggests that $\mr{SL}_2$-opers whose $p$-curvature are close to zero
may be thought of as candidates for suitable analogues of the canonical $\mr{SL}_2$-opers on  modular curves.
In fact, we prove the following assertion, which is 
one of the main results in the present paper (cf.  Theorem \ref{Co12} for a similar    assertion with slightly relaxed  assumptions).

\SSP
%---------------------------------------------------------------------[begin theorem]----------------------
\begin{intthm}[cf. Corollary \ref{ThAgg}] \label{ThA}
Let us keep the above notation.
Moreover, 
suppose that $l$ is even,  $\msF^\heartsuit$ is dormant (i.e. has vanishing $p$-curvature),  and $\msX$ is general in the moduli space $\overline{\mcM}_{g, r}$ of $r$-pointed stable curves of genus $g$. 
Then, there exists a {\bf canonical} splitting  of (\ref{Eq120}).
In particular,  we obtain  
a {\bf canonical} decomposition 
\begin{align} \label{Eqff131}
H^1_{\mr{dR}} (X, \mr{Sym}^{l-2}\msF) =   H^0  (X, \varTheta^{\otimes l}) \oplus H^0 (X, \varTheta^{\otimes l}(-D))^\vee
\end{align}
of  the $1$-st de Rham cohomology $H^1_{\mr{dR}} (X, \mr{Sym}^{l-2}\msF)$ of $\mr{Sym}^{l-2}\msF$.
 \end{intthm}
%------------------------------------------------------------------------------[end theorem]-------------
\SSP

Other relatively important results are Theorems \ref{Pr2}, \ref{Th9m}, \ref{Prc1}, \ref{Th67}, ~\ref{Th66}, \ref{Th77}, \ref{Co12}, and Corollaries \ref{Co34}, \ref{Co34f}.

\LSP

%%%%%%%%%%%%%%%%%%%%%%%%%%%%%%%%%%%%%%%%%%%%%
%%%%%%%%%%%%%%---[begin section]---%%%%%%%%%%%%%%
\subsection{Outline  of  our argument} \label{SS0009}
The present paper develops a general theory of opers on pointed stable curves and their infinitesimal deformations, that may be positioned as a continuation of the related work in ~\cite{Wak8}.
(The first half of our discussion proceeds under the situation where the characteristic of the base field  is not necessarily positive.)
Here, we briefly describe the arguments and results leading to
the above theorem.
(There are many assertions proved in this text that are not directly involved in proving Theorem \ref{ThA}. For example, the study of parabolic de Rham cohomology groups and their self-duality assertion  are included; see Corollary \ref{Coo78}.
Regarding the self-duality, we can find some previous results involved, e.g., ~\cite[Theorem 5.2, (a)]{Can}, ~\cite[Remark 3.2]{Og2}, and ~\cite[Theorem 2.7, (i)]{Sch}.)

Let $G$ be  a semisimple algebraic group of adjoint type with certain properties (cf.
the assumption imposed at the beginning of \S\,\ref{Ssde}).
The GIT quotient of its  Lie algebra $\mfg := \mr{Lie} (G)$ by the adjoint $G$-action defines a $k$-scheme $\mfc$.
We shall fix an element $\rho$ of $\mfc (k)^{\times r}$, which lifts to  an element $\widetilde{\rho} \in \mfg^{\times r}$ via the Kostant section (cf. (\ref{Eq202})).
 
 Denote by 
\begin{align}
\C_{G} 
\end{align}
(cf. (\ref{Eq3334}))
 the moduli stack 
classifying   flat $G$-bundles 
 on $\msX$.
As proved in Theorem \ref{Pr2},  
the moduli stack
\begin{align}
\mcO p_{G}
\ \left(\text{resp.,} \ \mcO p_{G, \rho} \right)
 \end{align}
classifying  $G$-opers (resp., $G$-opers of radii $\rho$) forms a locally closed substack of $\C_G$.

If  we are given a $G$-oper $\msE^\spadesuit$ classified by a closed point $q \in \mcO p_{G, \rho}$,
then the inclusion $\mcO p_{G} \migiincl \C_G$ 
induces, via differentiation, a $k$-linear inclusion $T_q\mcO p_G \migiincl T_q \C_G$
between the tangent spaces at $q$.
We investigate these spaces, as well as the relationship between them,  by describing relevant  infinitesimal deformations in terms of various  (de Rham) cohomology  groups  associated to  $\msE^\spadesuit$.
The resulting description allows us to obtain 
a  short exact sequence
\begin{align} \label{Et4}
0 \longmigi  T_q \mcO p_G 
&\xrightarrow{\mr{inclusion}}
 T_q \C_{G} \longmigi T_q^\vee \mcO p_{G, \rho} \longmigi 0  \end{align}
(cf. Theorem \ref{Th9m}),  where $T^\vee_q (-)$ denotes the cotangent space at $q$.

Next, since the base field  $k$ is of  characteristic $p$, 
one can define
the closed substack
\begin{align}
\C_G^{\psi = 0}
\end{align}
(cf. (\ref{Eq220})) of $\C_G$
   classifying flat $G$-bundles with vanishing $p$-curvature.
The intersection
\begin{align} \label{Eq35ggh}
\mcO p_G^{^\mr{Zzz...}} := \mcO p_G \cap \C_G^{\psi = 0} \ 
\end{align}
(cf. (\ref{Eq212}), (\ref{Eq211})) is nothing but  the moduli space of  dormant $G$-opers 
 and was substantially investigated  in ~\cite{Wak8} from  a viewpoint of enumerative geometry.
One important fact proved in  {\it loc.\,cit.} 
is the  \'{e}taleness of $\mcO p_G^{^\mr{Zzz...}}\!$ for a sufficiently general  $\msX$ (cf. ~\cite[Theorem 8.19]{Wak8}). 
(This fact also played  an essential  role in the proof of Joshi's conjecture; see ~\cite[Theorem  9.10]{Wak8}.)
Since the equality of dimensions  
\begin{align} \label{Ghr33}
\mr{dim}(T_q\mcO p_G) + \mr{dim}(T_q\C_G^{\psi = 0}) = \mr{dim}(T_q\C_G)
\end{align}
 holds (cf. Proposition \ref{Eqw3}), 
the  \`{e}taleness means that   $ \mcO p_G$ and $\C_G^{\psi = 0}$ 
 intersect transversally.
 Hence, 
 the composite
 \begin{align}
 T_q \C_{G}^{\psi = 0} \migiincl  T_q \C_{G} \migisurj T^\vee_q \mcO p_{G, \rho}
 \end{align}
 is an isomorphism, and the composite of its inverse and the inclusion 
 $T_q \C_{G}^{\psi = 0} \migiincl  T_q \C_{G}$ determines a split injection of (\ref{Et4}).
 It follows that 
 there exists a canonical decomposition
\begin{align} \label{Erkij}
T_q \mcO p_G \oplus 
T_q^\vee \mcO p_G  = T_q \C_{G} \
\end{align}
of 
$T_q \C_{G}$.

Finally, let us consider the case of  $G = \mr{PGL}_n$ with $n = \frac{p-1}{2}$.
(This is only the setting considered here for simplicity, and the main text will deal with a general $G$.)
Suppose further that 
 $\msE^\spadesuit$ comes from 
 the ``$\msF^\heartsuit$" as in the statement of Theorem \ref{ThA}
 via change of structure group  $\mr{SL}_2 \migi \mr{PGL}_n$ obtained  by integrating
 a suitable principal  $\mfs \mfl_2$-subalgebra of  $\mfp \mfg \mfl_n$.
Then,  a natural affine structure on $\mcO p_G$ yields
direct sum decompositions
\begin{align}
T_q \mcO p_G = \bigoplus_{j=1}^{n-1} H^0 (X, \varTheta^{\otimes 2j}),
\hspace{10mm}
T_q^\vee \mcO p_{G, \rho} = \bigoplus_{j=1}^{n-1} H^0 (X, \varTheta^{\otimes 2j}(-D))^\vee
\end{align}
(cf. ~\cite[Theorem 2.24]{Wak8}).
On the other hand, we show  (cf. (\ref{Eq3033}), (\ref{Eq834})) that
the $k$-vector space $T_q \C_{G}$ decomposes into a direct sum
\begin{align}
T_q \C_{G} = \bigoplus_{j=1}^{n-1} H_{\mr{dR}}^1 (X, \mr{Sym}^{2j} \msF)
\end{align}
according to the irreducible decomposition of the $\mfs \mfl_2$-algebra $\mfp \mfg \mfl_n$.
By various definitions involved,  (\ref{Erkij}) turns out to be  compatible with the respective structures of  direct sum decomposition. 
 Thus, the required decomposition (\ref{Eqff131}) can be obtained by considering its $(l/2)$-th factor.  
 
\SSP
%----------------------------------------------------------------------------------------------
\begin{rema} \label{Erdki}
As observed in ~\cite{Wak13}, 
it can be seen that the construction of the canonical decomposition (\ref{Eqff131}) discussed above is entirely similar to that of the classical Eichler-Shimura isomorphism, by regarding $\C_G^{\psi = 0}$ as a positive characteristic analogue of the moduli space of $G$-local systems  with real monodromy. 
\end{rema}
%----------------------------------------------------------------------------------------------

\LSP
%%%%%%%%%%%%%%%%%%%%%%%%%%%%%%%---[begin section]---%%%%%%%%%%%%%%
\subsection{Notation and Conventions}  \label{Ssde}
The basic terminology and notation in our discussion follow  ~\cite{Wak8}.
Throughout the present paper, we fix
   an algebraically closed field $k$,  a semisimple algebraic group $G$ over $k$,    a maximal torus $T$ of $G$, and 
 a Borel subgroup $B$ of $G$ containing $T$.
Write $\mfg$,  $\mfb$, and $\mft$ for the Lie algebras of $G$, $B$, and $T$, respectively.
Whenever  we discuss  the case of $G=\mr{SL}_n$ or $\mr{PGL}_n$,   
 the group $B$ (resp., $T$)  is assumed to be the subgroup of  $G$ consisting of upper triangular matrices (resp., diagonal matrices).

Regarding the algebraic group $G$, we assume that either ``$\mr{char}(k)=0$" or ``$\mr{char}(k) = p >2h_G$" or ``$G = \mr{PGL}_n$ with $1 < n < \mr{char}(k)$"
 is fulfilled, where  $h_G$ denotes   $1$ plus the height of a root with maximum height, as defined in  ~\cite[\S\,1.4.1]{Wak8}.
(If $G$ is simple, then $h_G$ coincides with its Coxeter  number.)
In particular,  when $\mr{char}(k)= p>0$, the prime $p$ is very good for $G$, in the sense of ~\cite[Chap.\,VI, Definition 1.6]{KW} (cf. the proof of Lemma \ref{Lem1}).

Next,
given an fs log scheme $T^\mr{log}$ and an fs  log scheme  $Y^\mr{log}$ over $T^\mr{log}$,
we shall write $\Omega_{Y^\mr{log}/T^\mr{log}}$ for the sheaf of logarithmic $1$-forms on $Y^\mr{log}$ over $T^\mr{log}$,  and write $\mcT_{Y^\mr{log}/T^\mr{log}}$ for its dual, i.e., the sheaf of locally defined logarithmic derivations on $\mcO_Y$ relative to $T^\mr{log}$.

We fix a pair of nonnegative integers $(g, r)$  with $2g-2 +r >0$ and an $r$-pointed stable curve  $\msX := (X, \{ \sigma_i \}_{i=1}^r)$  of genus $g$ over $k$, where $X$ denotes the underlying curve and $\{ \sigma_i \}_{i=1}^r$ denotes  the set  of its marked points.
Recall from ~\cite[Theorem 2.6]{KaFu}  that $\mr{Spec}(k)$ and $X$ admit canonical  log structures pulled-back, via the classifying map of $\msX$,  from the moduli stack of $r$-pointed stable curves of genus $g$ and the universal family of curves over that stack, respectively;
 the resulting log schemes are denoted by $\mr{Spec}(k)^\mr{log}$ and $X^\mr{log}$, respectively.
For simplicity, we write $\Omega := \Omega_{X^\mr{log}/\mr{Spec}(k)^\mr{log}}$ and 
$\mcT := \mcT_{X^\mr{log}/\mr{Spec}(k)^\mr{log}}$.
Also, write $d$ for 
the logarithmic universal derivation $\mcO_X \migi \Omega$. 

Let $D$ denote the divisor $D$ defined as the union of the $\sigma_i$'s.
For each $\mcO_X$-module $\mcF$, we set  ${^c}\mcF := \mcF \otimes_{\mcO_X} \mcO_X (-D)$.
In particular, ${^c}\Omega$ is isomorphic to the dualizing sheaf $\omega_X$ of $X/k$.

%%%%%%%%%%%%%%%%%%%%%%%%%%%%%%%%%%%%%%%%%%%%%%%%
%%%%%%%%%%%%%---[begin section]---%%%%%%%%%%%%%%
\vspace{10mm}
\section{Flat $G$-bundles and opers on a pointed stable curve}\label{S1}\SSP

In this section, we recall the definition of a flat $G$-bundle, as well as   a $G$-oper,  on a pointed stable curve.
We refer the reader  to, e.g.,  ~\cite{Wak8} for  detailed discussions on $G$-opers in such a general  formulation.
See also ~\cite{BD1}, ~\cite{JP} and  ~\cite{Mzk2}.
 
\LSP
%---------------------------[begin subsection]-------------
\subsection{Flat $G$-bundles on a pointed stable curve} \label{SS0111}

We shall fix a (right) $G$-bundle
 $\pi : \mcE \migi X$  on $X$.
Given a $k$-vector space $\mfh$ equipped with a $G$-action, we write
$\mfh_\mcE$ for the vector bundle on $X$ associated to the affine space $\mcE \times^G \mfh \left(:= (\mcE \times_k \mfh)/G \right)$.
In the case where $\mfh$ is taken to be the Lie algebra $\mfg$ equipped with the $G$-action via  the adjoint representation $\mr{Ad} : G \migi \mr{GL}(\mfg)$,
the resulting vector bundle $\mfg_\mcE$ is called the {\bf adjoint bundle} of $\mcE$.

By pulling-back the log structure of $X^\mr{log}$ via $\pi$,
one obtains a log structure on $\mcE$;  the resulting log scheme is denoted by $\mcE^\mr{log}$.
The $G$-action on $\mcE$ induces  a $G$-action on the direct image $\pi_* (\mcT_{\mcE^\mr{log}/\mr{Spec}(k)^\mr{log}})$ of $\mcT_{\mcE^\mr{log}/\mr{Spec}(k)^\mr{log}}$.
We shall set $\widetilde{\mcT}_{\mcE^\mr{log}} := \pi_* (\mcT_{\mcE^\mr{log}/\mr{Spec}(k)^\mr{log}})^G$, i.e., the subsheaf of $G$-invariant sections of $\pi_* (\mcT_{\mcE^\mr{log}/\mr{Spec}(k)^\mr{log}})$.
The differential of $\pi$ gives rise to   a canonical  short exact sequence of $\mcO_X$-modules
\begin{align} \label{Eqq10}
0 \longmigi \mfg_\mcE \longmigi \widetilde{\mcT}_{\mcE^\mr{log}}
 \xrightarrow{d_\mcE} \mcT \longmigi 0
 \end{align}
 (cf. ~\cite[\S\,1.2.5]{Wak8}).

 Recall that a {\it log $k$-connection}  (or, a {\it log connection}, for short) on $\mcE$ is, by definition,  an $\mcO_X$-linear morphism $\nabla : \mcT \migi \widetilde{\mcT}_{\mcE^\mr{log}}$ with $d_\mcE \circ \nabla = \mr{id}_{\mcT_{X^\mr{log}}}$.
 Since $\Omega$ is a line bundle,
 any  log connection is flat, i.e.,  has vanishing curvature (cf. ~\cite[Definition 1.23]{Wak8} for the definition of curvature).
 By a {\bf flat $G$-bundle} on $\msX$, we shall mean a pair $(\mcE, \nabla)$ consisting of a  (right) $G$-bundle $\mcE$ on $X$ and a  (flat) log  connection $\nabla$ on $\mcE$.
For example,  any flat $\mr{SL}_n$-bundle (where $n \geq 1$) on $\msX$ can be identified with a collection
\begin{align} \label{WWW19}
(\mcF, \nabla, \delta)
\end{align}
consisting of a rank $n$ vector bundle $\mcF$ on $X$, a log connection $\nabla$ on $\mcF$ in the usual  sense (i.e., a $k$-linear morphism $\nabla : \mcF \migi \Omega \otimes_{\mcO_X} \mcF$ satisfying the Leibnitz rule), and an isomorphism of flat bundles $\delta : (\mr{det}(\mcF), \mr{det}(\nabla)) \isom (\mcO_X, d)$, where $\mr{det}(\nabla)$ denotes the log connection on the determinant bundle $\mr{det}(\mcF)$ of $\mcF$ induced naturally from $\nabla$.

 We denote by
 \begin{align} \label{Eq3334}
\C_G
 \end{align}
 the moduli stack classifying flat $G$-bundles on $\msX$.
 It is verified (from, e.g.,  ~\cite[Theorem 1.0.1]{Wan}) that $\C_{G}$  may be represented by 
 an algebraic stack over $k$  in the sense of Artin, and locally of finite type over $k$.

\LSP
%---------------------------[begin subsection]-------------
\subsection{The Cartan decomposition of a Lie algebra} \label{SS051}

Denote
 %by $\mbG_m$ the multiplicative group over $k$ and 
 by  $\Gamma$
 %$\Gamma \subseteq \mr{Hom}(T, \mbG_m)$ ($:=$ the additive group of all characters of $T$)
 the set of simple roots in $B$ with respect to $T$.
For each character $\beta$ of $T$,
 %$\beta \in \mr{Hom}(T, \mbG_m)$,
  we set
\begin{align}
\mfg^\beta := \left\{ x \in \mfg \, | \, \mr{Ad}(t)(x) = \beta (t) \cdot x \ \text{for all} \ t \in T\right\}.
\end{align}
Recall that the Cartan decomposition of $\mfg$ is a Lie algebra grading  $\mfg = \bigoplus_{j \in \mbZ} \mfg_j$ such that $\mfg_j = \mfg_{-j} = 0$ for all $j > \mr{rk} (\mfg)$, where $\mr{rk}(\mfg)$ denotes the rank of $\mfg$, and $\mfg_0 = \mft$, $\mfg_1 = \bigoplus_{\alpha \in \Gamma} \mfg^\alpha$, $\mfg_{-1} = \bigoplus_{\alpha \in \Gamma} \mfg^{-\alpha}$.
It associates the  decreasing filtration  $\{ \mfg^j \}_{j \in \mbZ}$ on $\mfg$ given by  $\mfg^j := \bigoplus_{l \geq j} \mfg_l$.

Let us  fix a collection of data 
 $\,\qq := \{ x_\alpha \}_{\alpha\in \Gamma}$, 
 where each $x_\alpha$ denotes a generator of $\mfg^\alpha$, and write  
$q_{1} := \sum_{\alpha \in \Gamma} x_\alpha$.
(If we take two such collections, then they are conjugate by an element of $T$.
For that reason, the results obtained in the subsequent discussions are essentially independent of the choice of $\,\qq$.)
Also, by regarding  
 the fundamental coweight  $\check{\omega}_\alpha$ of $\alpha$ as an element of $\mft$ via differentiation, we obtain  $\check{\rho} := \sum_{\alpha \in \Gamma} \check{\omega}_\alpha \in \mft$.
Then,  there is  a unique  collection $\{ y_\alpha \}_{\alpha \in \Gamma}$ of generators  $y_\alpha \in \mfg^{-\alpha}$, such that the set 
 $\{  q_{-1}, 2 \check{\rho}, q_1  \}$, where $q_{-1}  := \sum_{\alpha \in \Gamma} y_\alpha \in \mfg_{-1}$,  forms an $\mfs \mfl_2$-triple.

If $G =\mr{PGL}_n$, then  
we identify its Lie algebra $\mfp \mfg \mfl_n$ with $\mfs \mfl_n$ via the  composite isomorphism $\mfs \mfl_n \migiincl \mfg \mfl_n \migisurj \mfp \mfg \mfl_n$,  and 
 assume that the triple   $\{  q_{-1}, 2 \check{\rho}, q_1  \}$ in $\mfs \mfl_n$ is always taken as 
\begin{align} \label{QQ0230}
& q_{-1} = \begin{pmatrix} 0 & 0 & \cdots  & 0 & 0 \\  1 & 0 & \cdots  & 0 & 0 \\  0 & 1 & \cdots  & 0 & 0 \\ \vdots  & \vdots  & \ddots  & \vdots  & \vdots  \\  0 & 0 & \cdots  & 1 & 0    \end{pmatrix}, \  
2 \check{\rho} = \begin{pmatrix} n-1 & 0 & 0 & \cdots   & 0 \\  0 & n-3 & 0 &  \cdots   & 0 \\  0 & 0 & n-5 & \cdots   & 0 \\ \vdots  & \vdots  & \vdots  & \ddots  & \vdots  \\  0 & 0 & 0 &  \cdots   & - (n-1)
\end{pmatrix}, 
\end{align} 
\begin{align}
 & 
  q_1 = \begin{pmatrix} 0 & n-1 & 0 & 0 & \cdots   & 0 \\  0 & 0 & 2(n-2) & 0 &  \cdots   & 0 \\  0 & 0 & 0  & 3(n-2) & \cdots   & 0 \\ \vdots  & \vdots  & \vdots  & \vdots &  \ddots  & \vdots  \\  0 & 0 & 0 & 0&  \cdots   & n-1   \\ 0 & 0 & 0 & 0&  \cdots   &0    \end{pmatrix}. 
  \notag
\end{align}

Denote by $\mfc_G$ (or  $\mfc$,  for simplicity) the GIT quotient
  of $\mfg$ by the adjoint $G$-action.
In particular, there exists a natural projection $\chi : \mfg \migisurj \mfc$.
This morphism
factors through the projection $\mfg \migi [\mfg/G]$ to the quotient stack $[\mfg/G]$ of $\mfg$ by the adjoint $G$-action.
 The resulting morphism  of $k$-stacks $[\mfg/G] \migi \mfc$ will be denoted by $[\chi]$.

Next, let us consider the space  $\mfg^{\mr{ad}(q_1)}$
 of $\mr{ad}(q_1)$-invariants, i.e.,
$\mfg^{\mr{ad}(q_1)} := \left\{ x \in \mfg \, | \, \mr{ad}(q_1) (x) = 0\right\}$.
It has dimension equal to $\mr{rk}(\mfg)$, and 
the Cartan decomposition $\mfg = \bigoplus_{j \in \mbZ} \mfg_j$ restricts to  a decomposition
$\mfg^{\mr{ad}(q_1)} = \bigoplus_{j \in \mbZ} \mfg_j^{\mr{ad}(q_1)}$
on $\mfg^{\mr{ad}(q_1)}$.
By the assumption on $G$ imposed in \S\,\ref{Ssde},
the   composite
\begin{align} \label{Eq141}
\mr{Kos} : \mfg^{\mr{ad}(q_1)} \xrightarrow{v \mapsto v + q_{-1}} \mfg \xrightarrow{\chi} \mfc.
\end{align}
becomes an  isomorphism of $k$-schemes.
Denote by $[0]_G$ (or simply, $[0]$, if there is no fear of confusion) the element of $\mfc (k)$ determined by the image of the zero section $0$ of $\mfg^{\mr{ad}(q_1)}$ via $\mr{Kos}$.

\LSP
%---------------------------[begin subsection]-------------
\subsection{Residues and radii} \label{SS0542}

 Let $(\mcE, \nabla)$ be a flat $G$-bundle on $\msX$.
 By a {\bf marking}  on  $(\mcE, \nabla)$, we mean   a collection of Lie algebra isomorphisms $\sigma_i^*(\mfg_\mcE)\isom \mfg$ ($i = 1, \cdots, r$).
 The notion of an isomorphism between flat $G$-bundles with marking can be defined in an evident manner.

 Under the assumption that $r > 0$,  let us choose  $i \in \{1, \cdots, r \}$, and moreover choose a local function $t \in \mcO_X$ defining the closed subscheme   $\mr{Im}(\sigma_i)$ of $X$. 
 Then, 
 the element 
 \begin{align} \label{Eq148}
 \mu_i^\nabla := \overline{\nabla \left(  t \frac{d}{d t}\right)}  \in  \sigma_i^*(\widetilde{\mcT}_{\mcE^\mr{log}})
 \end{align}
 lies in $\sigma_i^* (\mfg_\mcE)$  and does not depend on the choice of $t$ (cf. the discussion in ~\cite[\S\,1.6.3]{Wak8}).
 We refer to $\mu_i^\nabla$ as the {\bf residue} (or, the {\bf monodromy operator}, in the terminology of  ~\cite{Mzk2} and ~\cite{Wak8}) of $\nabla$ at $\sigma_i$.
If $(\mcE, \nabla)$ is equipped with a marking, then 
  $\mu_i^\nabla$  can be regarded as   an element of $\mfg$
 by passing to the isomorphism  $\sigma^*_i (\mfg_\mcE) \isom \mfg$ given by this marking.
 In the case of $r = 0$,  any log connection  is regarded as  being  equipped with the trivial marking and having  {\it residue $\emptyset$}.

Thus, for each $r$-tuple $\mu := (\mu_i)_{i=1}^r$ of elements of $\mfg$ (where $\mu := \emptyset$ if $r = 0$), we obtain 
 the moduli stack 
  \begin{align} \label{Eq4493}
\C_{G, \mu} 
 \end{align}
classifying  flat $G$-bundles on $\msX$ with marking whose residue  at $\sigma_i$ coincides with $\mu_i$ for every $i=1, \cdots, r$.
  
  Next,
  note that, for each $i=1, \cdots, r$,
  the pair $(\sigma_i^* (\mcE), \mu_i^\nabla)$ specifies a $k$-rational point of $[\mfg /G]$; it defines a $k$-rational point
\begin{align} \label{Eq4492}
\rho_i^\nabla := [\chi] \left( (\sigma_i^* (\mcE), \mu_i^\nabla)\right) \in \mfc (k)
\end{align}
  via the morphism 
 $[\chi] : [\mfg /G] \migi \mfc$.
This element is well-defined without choosing any marking, and 
we  refer to $\rho_i^\nabla$ as the {\bf radius} of $\nabla$ at $\sigma_i$.
Given an $r$-tuple  $\rho := (\rho_i)_{i=1}^r \in \mfc^{\times r}$, 
  we say that $\nabla$ is {\bf of radii $\rho$} if its radius at $\sigma_i$ coincides with $\rho_i$ for every $i =1, \cdots, r$.
  In the case of $r = 0$, any log connection  is regarded as having {\it radius $\emptyset$}.

  The stack $\C_{G}$ admits a closed substack 
  \begin{align} \label{Eq4491}
  \C_{G, \rho}
  \end{align}
   classifying flat $G$-bundles on $\msX$ of radii $\rho$.
 By  assigning the residues of each flat bundle to its radii  via $\chi : \mfg \migi \mfc$, we have
   a morphism of $k$-stacks
  \begin{align} \label{Eq4490}
  \C_{G, \mu} \migi \C_{G, \chi (\mu)},
  \end{align}
  where $\chi (\mu):= (\chi (\mu_i))_{i=1}^r$.

\LSP
%---------------------------[begin subsection]-------------
\subsection{$G$-opers on a pointed stable curve} \label{SS052}

Let $\mcE_B$ be a $B$-bundle on $X$.
Denote by $\mcE_G$ the $G$-bundle obtained from $\mcE_B$ via change of structure group by the inclusion $B \migiincl G$.
(Hence, $\mcE_B$ specifies a $B$-reduction of $\mcE_G$.)
By using the injection $\widetilde{\mcT}_{\mcE_B^\mr{log}} \migiincl \widetilde{\mcT}_{\mcE_G^\mr{log}}$ induced by the natural inclusion $\mcE_B \migiincl \mcE_G$,
we regard $\widetilde{\mcT}_{\mcE_B^\mr{log}}$ as an $\mcO_X$-submodule of $\widetilde{\mcT}_{\mcE_G^\mr{log}}$.
Since the $k$-vector subspace
$\mfg^j$ of $\mfg$ (for each $j \in \mbZ$)
 is closed under the adjoint $B$-action, it induces a subbundle $\mfg_{\mcE_B}^j$ of $\mfg_{\mcE_B} \left(= \mfg_{\mcE_G} \right)$.
For each $j \in \mbZ$, we shall set 
\begin{align}
\widetilde{\mcT}_{\mcE_G^\mr{log}}^j := \widetilde{\mcT}_{\mcE_B^\mr{log}} + \mfg_{\mcE_B}^j \left(\subseteq \widetilde{\mcT}_{\mcE_G^\mr{log}} \right).
\end{align}

The inclusion $\mfg_{\mcE_B} \left(= \mfg_{\mcE_G} \right) \migiincl \widetilde{\mcT}_{\mcE_G^\mr{log}}$ yields an isomorphism 
\begin{align} \label{Eqq49}
\mfg_{\mcE_B}^{j}/\mfg_{\mcE_B}^{j+1} \isom \widetilde{\mcT}_{\mcE_G^\mr{log}}^{j}/\widetilde{\mcT}_{\mcE_G^\mr{log}}^{j+1}.
\end{align}
On the other hand, since each $\mfg^{-\alpha}$ ($\alpha \in \Gamma$) is closed under the  adjoint $B$-action,  the decomposition $\mfg^{-1}/\mfg^0 = \bigoplus_{\alpha \in \Gamma} \mfg^{-\alpha}$ yields a canonical isomorphism $\mfg^{-1}_{\mcE_B}/\mfg^{0}_{\mcE_B} \isom \bigoplus_{\alpha \in \Gamma} \mfg_{\mcE_B}^{- \alpha}$.
By composing it with the inverse of  (\ref{Eqq49}) for $j=-1$,
we obtain  a decomposition
\begin{align} \label{Eq10}
\widetilde{\mcT}_{\mcE_G^\mr{log}}^{-1}/\widetilde{\mcT}_{\mcE_G^\mr{log}}^0 \isom \bigoplus_{\alpha \in \Gamma} \mfg_{\mcE_B}^{- \alpha}.
\end{align}

Now, let us consider a pair $\msE^\spadesuit := (\mcE_B, \nabla)$ consisting of a $B$-bundle $\mcE_B$ on $X$ and a log connection $\nabla$ on the $G$-bundle $\mcE_G := \mcE_B \times_B G$ induced by $\mcE_B$.
Recall that $\msE^\spadesuit$ is called a {\bf $G$-oper} on $\msX$(cf. ~\cite[Definition 2.1]{Wak8})  if it satisfies the following two conditions:
\begin{itemize}
\item
$\nabla (\mcT) \subseteq \widetilde{\mcT}_{\mcE_G^\mr{log}}^{-1}$;
\item
For any $\alpha \in \Gamma$, the composite
\begin{align}
\mcT\xrightarrow{\nabla} \widetilde{\mcT}_{\mcE_G^\mr{log}}^{-1} \migisurj \widetilde{\mcT}^{-1}_{\mcE_G^\mr{log}}/\widetilde{\mcT}^0_{\mcE_G^\mr{log}} \migisurj \mfg^{-\alpha}_{\mcE_B}
\end{align}
is an isomorphism, where the third arrow denotes the natural projection with respect to the decomposition (\ref{Eq10}).
\end{itemize}

If $r >0$ and $\rho$ is an element of $\mfc (k)^{\times r}$, then we say that a $G$-oper $\msE^\spadesuit := (\mcE_B, \nabla)$ is  {\bf of radii $\rho$} if the flat $G$-bundle $(\mcE_G, \nabla)$ is of radii $\rho$ (cf. ~\cite[Definition 2.32]{Wak8}).
In the case of  $r = 0$, any $G$-oper is regarded as  being {\bf of radius $\emptyset$}.

\SSP
%------------------------------------------------------------------------------------------
\begin{rema}[Case of $G = \mr{SL}_n$] \label{Rem667}
Each $\mr{SL}_n$-oper can be described in terms of vector bundles 
because
a $B$-reduction of an $\mr{SL}_n$-bundle corresponds to  a complete flag on the associated vector bundle.
To be precise,
a pair $(\mcE_B, \nabla)$ of a $B$-bundle $\mcE_B$ and a connection $\nabla$ on $\mcE_{\mr{SL}_n} := \mcE_B \times_B \mr{SL}_n$ can be translated into 
 a collection of data
\begin{align} \label{Eq58}
\msF^\heartsuit := (\mcF, \nabla, \{ \mcF^j \}_{j=0}^n, \delta)
\end{align}
such that $(\mcF, \nabla, \delta)$ is a triple as in (\ref{WWW19}) and 
$\{ \mcF^j \}_{j=0}^n$ denotes a decreasing filtration $0= \mcF^n \subseteq \mcF^{n-1} \subseteq \cdots \subseteq \mcF^1 \subseteq \mcF^0 = \mcF$ of $\mcF$
whose graded pieces 
$\mcF^j /\mcF^{j+1}$ ($i=0, \cdots, n-1$) are  line bundles.
Moreover, such a collection corresponds to an $\mr{SL}_n$-oper 
if it
satisfies  the following two conditions:
 \begin{itemize}
 \item
 For every $j = 1, \cdots, n-1$, the inclusion relation $\nabla (\mcF^j) \subseteq \Omega  \otimes_{\mcO_X}  \mcF^{j-1}$ holds;
 \item
For every $j= 1, \cdots, n-1$, the $\mcO_X$-linear morphism 
\begin{align} \label{Eqq212}
\mr{KS}_{\msF^\heartsuit}^j  : 
\mcF^j/\mcF^{j+1} \migi \Omega  \otimes_{\mcO_X} (\mcF^{j-1}/\mcF^{j})
\end{align}
 induced from  $\nabla$ (by taking account of the former  condition)
 is an isomorphism.
 (The isomorphism $\mr{KS}_{\msF^\heartsuit}^j$ is called the $j$-th {\it Kodaira-Spencer map} for $\msF^\heartsuit$.)
 \end{itemize}
 Also, such a collection without a choice of an isomorphism $\delta$ is called a {\bf $\mr{GL}_n$-oper} (cf. ~\cite[Definition 4.17]{Wak8}).
 Even when dealing with an $\mr{SL}_n$-oper, we will omit the datum   ``$\delta$" 
 %,  if there is no fear of confusion, 
 %omitted
  in its notation for simplicity.
 \end{rema}
 %------------------------------------------------------------------------------------------
%\SSP

%%%%%%%%%%%%%%%%%%%%%%%%%%%%%%%%%%%%%%%%%%%%%%%%
%%%%%%%%%%%%%---[begin section]---%%%%%%%%%%%%%%
\vspace{10mm}
\section{The moduli space of opers and change of structure group}\label{S2}\SSP

This section deals with the moduli stack of $G$-opers and recall
a construction of $G$-opers  using $\mr{PGL}_2$-opers (or $\mr{SL}_2$-opers).
%  together with a choice of an  $\mfs \mfl_2$-triple in  $\mfg$.

\LSP
%---------------------------[begin subsection]-------------
\subsection{The moduli space of $G$-opers} \label{SS055}

Denote by 
\begin{align} \label{Eq5009}
\mcO p_G 
\end{align}
the moduli stack classifying $G$-opers  on $\msX$.
 {\it  In the rest of the present paper, we assume that $G$ is of adjoint type unless otherwise stated.}
  Then, 
it follows from  ~\cite[Theorem A]{Wak8}
that 
$\mcO p_{G}$
  may be represented by  an affine scheme over $k$ (cf. \S\,\ref{SS036} for more details).
The assignment $(\mcE_B, \nabla) \mapsto (\mcE_G, \nabla)$ (for each $G$-oper $(\mcE_B, \nabla)$) determines  a morphism of $k$-stacks
\begin{align} \label{Eq201}
\mr{Imm}_G : \mcO p_{G} \migi \C_{G}.
\end{align}

Next, 
suppose 
 that we are given an element $\rho := (\rho_i)_{i=1}^r$ of $\mfc (k)^{\times r}$ (where $\rho := \emptyset$ if $r = 0$).
It determines a closed subscheme 
\begin{align} \label{Eq5007}
\mcO p_{G, \rho}
\end{align}
of $\mcO p_G$ classifying $G$-opers of radii $\rho$.
The morphism  $\mr{Imm}_G$ restricts to 
a  morphism of $k$-stacks
\begin{align} \label{Eq5001}
\mr{Imm}_{G, \rho} : \mcO p_{G, \rho} \migi \C_{G, \rho}.
\end{align}
 
 Let us take a $G$-oper $\msE^\spadesuit$ on $\msX$ of radii $\rho$.
 According to 
~\cite[Proposition 2.19]{Wak8},
 there exists a unique (up to isomorphism) pair $({^\dagger}\msE^\spadesuit, \mr{nor}_{\msE^\spadesuit})$ consisting of a   $\,\qq$-normal $G$-oper  ${^\dagger}\msE^\spadesuit := ({^\dagger}\mcE_B, \nabla)$ on $\msX$ and an isomorphism of $G$-opers $\mr{nor}_{\msE^\spadesuit} : {^\dagger}\msE^\spadesuit \isom \msE^\spadesuit$ (cf. ~\cite[Definition 2.14]{Wak8} for the definition of  $\,\qq$-normality).
Here, recall from ~\cite[\S\,2.4.3]{Wak8}  that, for each $i =1, \cdots, r$,  there exists a canonical $G$-equivariant  isomorphism  $\sigma^*_i ({^\dagger}\mcE_G) \isom G$,  where ${^\dagger}\mcE_G := {^\dagger}\mcE_B \times^B G$, that induces a Lie algebra  isomorphism  $\left(\sigma^*_i (\mfg_{{^\dagger}\mcE_G}) =  \right) \sigma^*_i (\mfg_{{^\dagger}\mcE_B}) \isom \mfg$.
In this way,  the underlying flat $G$-bundle $({^\dagger}\mcE_G, \nabla)$ of ${^\dagger}\msE^\spadesuit$
admits a canonical choice  of a marking.

Also,
by the definition of  $\,\qq$-normality,
 the residue of $\nabla$ at $\sigma_i$ coincides with $q_{-1} + \mr{Kos}^{-1}(\rho_i)$.
If we write 
\begin{align} \label{Eq202}
\widetilde{\rho}  \left(= (\widetilde{\rho}_i)_{i=1}^r \right) := (q_{-1}+ \mr{Kos}^{-1}(\rho_i))_{i=1}^r,
\end{align}
then
 the resulting assignment $\msE^\spadesuit  \mapsto ({^\dagger}\mcE_G, \nabla)$ determines  a morphism of $k$-stacks
\begin{align} \label{Eq200}
\mr{Imm}_{G, \widetilde{\rho}} : \mcO p_{G, \rho} \migi \C_{G, \widetilde{\rho}}.
\end{align}

The following assertion generalizes the fact in ~\cite[\S\,3.1.11, Remark, (iii)]{BD1}, dealing with the case where the underlying curve is an unpointed  smooth projective curve over the field of complex numbers.

\SSP
%------------------------------------------------------------------------------------------
\bt \label{Pr2}
(Recall that $G$ is assumed to be of adjoint type.)
The morphism of $k$-stacks  $\mr{Imm}_G$ (resp., $\mr{Imm}_{G, \rho}$; resp., $\mr{Imm}_{G, \widetilde{\rho}}$) is schematic and an immersion.
\et
%------------------------------------------------------------------------------------------
\begin{proof}
It suffices to consider the non-resp'd assertion because the proofs of the remaining ones are similar.
Denote by $\mcB un_G$ (resp., $\mcB un_B$) the moduli stack of $G$-bundles (resp., $B$-bundles) on $X$.
By forgetting the data of log connections, we obtain  a morphism of $k$-stacks
$\C_G \migi \mcB un_G$.
Note that the morphism  $\mcB un_B \migi \mcB un_G$ induced by 
 the change of  structure group using  $B \migiincl G$ is schematic 
  (cf. ~\cite[Corollary 3.2.4]{Wan}),
and $\mcO p_G$ can be represented by a locally closed substack of the fiber product $\C_G \times_{\mcB un_G} \mcB un_B$.
Hence,
the morphism $\mr{Imm}_G$ is verified to be  schematic.
The result of  Lemma \ref{Lem1} described below implies that $\mr{Imm}_G$ is   a monomorphism.
By Lemma \ref{Lem145},
this morphism also  satisfies the valuative criterion $(*)$ in the sense of ~\cite[Chap\,I, \S\,2.4]{Mzk2}.
In particular, it is a radimmersion (cf. ~\cite[Chap\,I, Theorem 2.12]{Mzk2}).
Consequently,   it follows from ~\cite[Chap\,I, Corollary  2.13]{Mzk2}
that  $\mr{Imm}_G$ defines an immersion.
\end{proof}
%------------------------------------------------------------------------------------------
\SSP

The following two lemmas were applied in the proof of the above proposition.

\SSP
%------------------------------------------------------------------------------------------
\ble \label{Lem1}
Let $\msE^\spadesuit_i := (\mcE_{B, i}, \nabla_i)$ ($i=1,2$) be $G$-opers on $\msX$ and
$\alpha_G : (\mcE_{G, 1}, \nabla_1) \isom (\mcE_{G, 2}, \nabla_2)$ (where $\mcE_{G, i} := \mcE_{B, i} \times^B G$) an isomorphism of flat $G$-bundles.
Then, there exists a unique isomorphism of $G$-opers $\alpha_B : \msE^\spadesuit_1 \isom \msE^\spadesuit_2$ inducing $\alpha_G$ via change of structure group by $B \migiincl G$.
In particular, if  $(\mcE_G, \nabla)$ is  a flat $G$-bundle  on $\msX$,  then a $B$-reduction $\mcE_B$ of $\mcE_G$
for which the pair $(\mcE_B, \nabla)$ defines a $G$-oper is (if it exists) uniquely determined up to isomorphism.
\ele
%------------------------------------------------------------------------------------------
\begin{proof}
We may assume, without loss of generality,  that $\mcE := \mcE_{G, 1} = \mcE_{G, 2}$, $\nabla := \nabla_1 = \nabla_2$, and that $\alpha_G$ 
coincides with the identity morphism of $\mcE$.
The two Borel reductions $\mcE_{B, 1}$, $\mcE_{B, 2}$ of $\mcE$ define the $\mcO_X$-submodules of $\mfg_\mcE$, i.e.,  $\mfg^0_{\mcE_{B, 1}} \left(= \mfb_{\mcE_{B, 1}}\right)$ and  $\mfg^0_{\mcE_{B, 2}} \left(= \mfb_{\mcE_{B, 2}}\right)$.
Notice here that, for each $i \in \{1, 2\}$ and $j \in \mbZ$,
the subquotient  $\mfg^j_{\mcE_{B, i}}/\mfg^{j+1}_{\mcE_{B, i}}$
is isomorphic to $\Omega^{\otimes j} \otimes_k\mfg_{j}$.
Hence, by using $\mr{deg}(\Omega)> 0$,
we can verify  the equality 
$\mfg^1_{\mcE_{B, 1}} = \mfg^1_{\mcE_{B, 2}}$.

Next, let us take a closed point $x$ of $X$, and fix an identification 
$\mcE_{B, 1}|_x = B$, which extends to 
$\mcE |_x = G$  and moreover induces  $\mfg^1_{\mcE_{B, 1}} |_x = \mfg^1$ via differentiation.
Then, we can write $\mcE_{B, 2}|_x = h B \left(\subseteq G \right)$ for some $h \in G$.
It follows that $\mfg^1_{\mcE_{B, 2}} |_x$ coincides with the image  $\mr{ad}(h)(\mfg^1)$
of $\mfg^1$ via the adjoint operator $\mr{ad}(h)$ given by $h$.
Since the characteristic of the base field $k$ is very good for $G$ if it is of positive characteristic,   
the equality $\mfg^1_{\mcE_{B, 1}} = \mfg^1_{\mcE_{B, 2}}$ proved above together with 
  the Springer isomorphism (cf. e.g., ~\cite[Theorem 3.3]{KW})   implies  
$N = h N h^{-1}$, where $N$ denotes the unipotent radical of $B$.
It follows that   $B = h B h^{-1}$ (cf. ~\cite[\S\,23.1, Corollary D]{Hum}), and we have  
$h \in B$ (cf. ~\cite[\S\,23.1, Theorem]{Hum}).
Therefore,  the equality  $\mcE_{B, 1}|_x = \mcE_{B, 2}|_x$ holds for all $x$.
This means  $\mcE_{B, 1} = \mcE_{B, 2}$, completing  the proof of this assertion.
\end{proof}
%------------------------------------------------------------------------------------------
\SSP

%------------------------------------------------------------------------------------------
\ble \label{Lem145}
Let $T$ be a trait, i.e., the spectrum of a discrete valuation ring.
Denote by $\eta$ its generic point and by $t$ its closed point.
Now, let us consider a collection
\begin{align}
(\msE_T, \mcE_{B, \eta}, \mcE_{B, t}),
\end{align}
where
\begin{itemize}
\item
$\msE_T := (\mcE_T, \nabla_T)$ denotes a flat $G$-bundle $\msE_T := (\mcE_T, \nabla_T)$ on the base-change $\msX_T := (X_T, \{ \sigma_{i, T} \}_{i=1}^r)$ of $\msX$ over $T$ (cf. ~\cite[\S\,1]{Wak8} for the formulation of a flat $G$-bundle on a family of pointed stable curves);
\item
$\mcE_{B, \eta}$ and  $\mcE_{B, t}$ are $B$-reductions of $\mcE_\eta := \mcE_T \times_T \eta$ and  $\mcE_{t} := \mcE_T \times_T t$, respectively, such that both $(\mcE_{B, \eta}, \nabla_T |_{\mcE_\eta})$ and $(\mcE_{B, t}, \nabla_T |_{\mcE_t})$ form $G$-opers.
\end{itemize}
Then, there exists a $B$-reduction $\mcE_{B, T}$ on $\mcE_T$ such that $\mcE_{B, T} \times_T \eta = \mcE_{B, \eta}$, $\mcE_{B, T} \times_T t = \mcE_{T, t}$, and   the pair $(\mcE_{B, T}, \nabla_T)$ forms a $G$-oper.
\ele
%------------------------------------------------------------------------------------------
\begin{proof}
Denote by $X_\eta$ the generic fiber of $X_T$.
Note that $\{ \mfg_{\mcE_\eta}^j \}_{j}$ and $\{ \mfg_{\mcE_t}^j \}_j$ coincide with the Harder-Narasimhan filtrations on $\mfg_{\mcE_\eta}$ and $\mfg_{\mcE_t}$, respectively.
According to ~\cite[Theorem 2.3.2]{HL},  $\{ \mfg_{\mcE_\eta}^j \}_{j}$ extends to a  decreasing filtration $\{ \mfg_{\mcE_T}^{j} \}_j$ (i.e., $\mfg_{\mcE_T}^{j} |_{X_\eta} = \mfg_{\mcE_\eta}^j$ for every $j$) on $\mfg_{\mcE_T}$ such that the subquotients $\mfg_{\mcE_T}^{j}/\mfg_{\mcE_T}^{j+1}$ are all $T$-flat.
By the upper semicontinuity of the Harder-Narasimhan type, 
the special fiber of $\{\mfg_{\mcE_T}^j \}_j$ must coincide with $\{ \mfg_{\mcE_t}^j \}_j$.
Here, recall that the isomorphism $G \isom \mr{Aut} (\mfg)^0$ (where $\mr{Aut} (\mfg)^0$ denotes the identity component of the group of Lie algebra automorphisms of $\mfg$) induced from the adjoint representation of $G$ restricts to an isomorphism of algebraic groups  $B \isom \left\{ h \in \mr{Aut}(\mfg)^0 \, | \, h (\mfb) \subseteq \mfb \right\}$ (cf. ~\cite[Remark 1.29]{Wak8}).
It follows that  the Lie subalgebra $\mfg_{\mcE_T}^0$ of $\mfg_{\mcE_T}$ determines a $B$-reduction $\mcE_{B, T}$ of $\mcE_{T}$, which  satisfies $\mcE_{B, T} \times_T \eta = \mcE_{B, \eta}$ and $\mcE_{B, t} \times_T t = \mcE_{B, t}$.
Since the generic and special fibers of $(\mcE_{B, T}, \nabla_T)$ coincide with the $G$-opers $(\mcE_{B, \eta}, \nabla_T |_{\mcE_\eta})$ and $(\mcE_{B, t}, \nabla_T |_{\mcE_t})$, respectively, the pair $(\mcE_{B, T}, \nabla_T)$ turns out to form a $G$-oper.
This completes the proof of this lemma.
\end{proof}
%------------------------------------------------------------------------------------------

\LSP
%----------------------------------------------------------------------[begin subsection]-------------
\subsection{Change of structure group} \label{SS0r2}

In the case of $G = \mr{PGL}_2$,
  we use the notations $G^\odot$, $B^\odot$, $\mfg^\odot $, and  $\mfb^\odot$ instead of $G$, $B$,  $\mfg$, and $\mfb$, respectively, for simplicity.
The $\mfs \mfl_2$-triple $\{  q_{-1}, 2  \check{\rho}, q_1  \}$ in $\mfg$
 associated to  the fixed collection $\ \qq$ 
   determines an injection
\begin{align} \label{injg} \iota_\mfg : \mfg^\odot \migiincl \mfg. 
\end{align}
To be precise, 
 $\iota_\mfg$ is given by  
$\iota_\mfg \left(\begin{bmatrix} 0 & 1 \\ 0 & 0 \end{bmatrix} \right) = q_1$,
$\iota_\mfg \left(\begin{bmatrix} 1 & 0 \\ 0 & -1 \end{bmatrix} \right) = 2 \check{\rho}$,
and $\iota_\mfg \left(\begin{bmatrix} 0 & 0 \\ 1 & 0 \end{bmatrix} \right) = q_{-1}$.
In particular, we have
 $\iota_\mfg (\mfb^\odot) \subseteq \mfb$.
 According to ~\cite[Proposition 2.10]{Wak8},  there exists 
 a unique  injective morphism of algebraic $k$-groups
 $\iota_B : B^\odot \migiincl B$ 
 whose differential induces 
   the restriction $\mfb^\odot \migiincl \mfb$ of $\iota_\mfg$.

Let  $\msE^\spadesuit_\odot := (\mcE_{B^\odot}, \nabla_\odot)$ be a $G^\odot$-oper   on $\msX$.
 We shall write 
 $\mcE_{G^\odot} := \mcE_{B^\odot} \times^{B^\odot} G^\odot$,
$\mcE_{B} := \mcE_{B^\odot} \times^{B^\odot, \iota_B} B$, 
and $\mcE_G := \mcE_B \times^{B} G$.
The injection $\widetilde{\mcT}_{\mcE^\mr{log}_{B^\odot}} \migiincl \widetilde{\mcT}_{\mcE^\mr{log}_{B}}$ induced by $\iota_B$ extends to an $\mcO_X$-linear  injection
$d \iota_{G} : \widetilde{\mcT}_{\mcE^\mr{log}_{G^\odot}} \migiincl 
\widetilde{\mcT}_{\mcE^\mr{log}_{G}}$.
The composite
$\iota_{G*}(\nabla_\odot) : \mcT_{X^\mr{log}} \xrightarrow{\nabla_\odot} \widetilde{\mcT}_{\mcE^\mr{log}_{G^\odot}} \xrightarrow{d \iota_{G}}
\widetilde{\mcT}_{\mcE^\mr{log}_{G}}$
specifies a log connection on $\mcE_G$, and 
 the resulting pair 
\begin{align} \label{associated} 
\iota_{G *}(\msE^\spadesuit_\odot) :=
 (\mcE_{B},  \iota_{G*}(\nabla_\odot))
 \end{align}
  forms a $G$-oper on $\msX$.
The assignment $\msE_\odot^\spadesuit \mapsto \iota_{G *}(\msE_\odot^\spadesuit)$
 determines a closed immersion between $k$-schemes 
\begin{align} \label{Eq405}
%\iota_G^{\mcO p} : 
\mcO p_{G^\odot} \migiincl  \mcO p_G
\end{align}
(cf. ~\cite[Theorem 2.24, (ii)]{Wak8}).

If $\iota_\mfc$ denotes  the morphism $\mfc_{G^\odot} \migi \mfc_{G}$ induced by $\iota_\mfg$ via the adjoint quotients,
then, for each $\rho_\odot := (\rho_{\odot, i})_{i=1}^r \in \mfc_{G^\odot} (k)^{\times r}$,
the morphism (\ref{Eq405}) restricts to a closed immersion
\begin{align} \label{Eq281}
%\iota_{G, \rho_\odot}^{\mcO p} :
 \mcO p_{G^\odot, \rho_\odot} \migiincl \mcO p_{G, \iota_\mfc (\rho_\odot)}
\end{align}
by putting $\iota_{\mfc} (\rho_\odot) := (\iota_{\mfc}(\rho_{\odot, i}))_{i=1}^r$, where $\rho_\odot := \emptyset$ and $\iota_{\mfc} (\rho_\odot) := \emptyset$ if $r = 0$ (cf. ~\cite[Remark 2.37]{Wak8}).
In particular, we obtain  a closed immersion 
$\mcO p_{G^\odot, [0]^{\times r}} \migiincl \mcO p_{G, [0]^{\times r}}$,
where
$[0]^{\times r} := ([0], \cdots, [0])$ in $\mfc_{G^\odot} (k)^{\times r}$ or 
$\mfc_{G} (k)^{\times r}$.

\LSP
%---------------------------[begin subsection]-------------
\subsection{The affine structure on the moduli space} \label{SS036}

Note that $\mfg^{\mr{ad}(q_1)}$ is closed under the $B^\odot$-action via $\iota_B$.
Therefore, 
if ${^\dagger}\mcE_{B^\odot}$ denotes the  $B^\odot$-bundle on $X$ constructed in ~\cite[(190)]{Wak8}, then
it induces a rank $\mr{rk}(\mfg)$ vector bundle
\begin{align}   \label{QQ801}
 \DV_{G}
:=  \Omega\otimes_{\mcO_X}  (\mfg^{\mr{ad}(q_1)})_{\DE_{B^\odot}}
\end{align}
on $X$ (cf. ~\cite[\S\,2.4.5]{Wak8}).
For example, $\DV_{G^\odot}$ is canonically isomorphic to  $\Omega^{\otimes 2}$.
The decomposition    $\mfg^{\mr{ad}(q_1)} = \bigoplus_{j \in \mbZ} \mfg_j^{\mr{ad}(q_1)}$
determines   
a decomposition 
$\DV_{G} = \bigoplus_{j \in \mbZ}  \DV_{G, j}$; it gives 
a  decreasing filtration $\{ \DV^j_{G} \}_{j \in \mbZ}$ on $\DV_G$
by putting  $\DV^j_{G} := \bigoplus_{j' \geq j}  \DV_{G, j'}$.
Also, for each $j \in \mbZ$, there exists an isomorphism 
\begin{align} \label{QR020}
\DV_{G, j} \isom  \Omega^{\otimes (j+1)} \otimes_k \mfg^{\mr{ad}(q_1)}_j.
\end{align}
 
  Recall that $\mcO p_{G^\odot}$ admits a structure of affine space  modeled on the space of  quadratic log differentials $H^0 (X, \Omega^{\otimes 2}) \left(= H^0 (X, \DV_{G^\odot}) \right)$ (cf. ~\cite[\S\,2.3.4]{Wak8}).
In what follows, let us review a canonical affine structure   on $\mcO p_G$ generalizing this fact (cf. ~\cite[Theorem 2.24, (i)]{Wak8}).
 The $B^\odot$-equivariant inclusion $\mfg^{\mr{ad}(q_1)} \migiincl \mfg$ yields
an $\mcO_X$-linear injection
\begin{align}    \label{inclV}
\varsigma   : \DV_{G} \migiincl \left( \Omega\otimes_{\mcO_X}   \mfg_{\DE_{B^\odot}}= \right)\Omega\otimes_{\mcO_X}   \mfg_{\DE_{B}},
\end{align}
where ${^\dagger}\mcE_B := {^\dagger}\mcE_{B^\odot} \times^{B^\odot, \iota_B} B$.
Also, $\iota_\mfg$ 
  induces 
 an $\mcO_X$-linear injection
%\begin{align}  \label{inclOmega}
$\Omega^{\otimes 2} \left(=  \DV_{G^\odot}  \right)  \migiincl \DV_{G}$.  
%\end{align}
By using
these injections,
 %$\varsigma$ (resp., $\iota_{\mfg}^\mcV$), 
 we shall regard $\DV_{G}$ and $\Omega^{\otimes 2}$ as  $\mcO_X$-submodules of $\Omega\otimes_{\mcO_X}   \mfg_{\DE_{B}}$ and  $\DV_{G}$, respectively.

Next, let us take a $\,\qq$-normal  $G^\odot$-oper ${^\dagger}\msE_{\odot}^\spadesuit := ({^\dagger}\mcE_{B^\odot}, \nabla_\odot)$ on $\msX$.
Given an element $R$ of $H^0  (X, \DV_G)$, regarded as 
 an
$\mcO_X$-linear morphism 
$\mcT \migi (\mfg^{\mr{ad}(q_1)})_{{^\dagger}\mcE_{B^\odot}} \left(\subseteq \widetilde{\mcT}_{\DE^\mr{log}_{G}} \right)$,
one may verify that 
the collection
\begin{align}
\iota_{G*}(\msE^\spadesuit_{\odot})_{+R} := ({^\dagger}\mcE_{B}, \iota_{G*}(\nabla_\odot) + R)
\end{align}
forms a $\,\qq$-normal $G$-oper.
The resulting assignment $(\msE^\spadesuit_\odot, R) \mapsto \iota_{G*}(\msE^\spadesuit_{\odot})_{+R}$ defines an isomorphism of $k$-schemes
\begin{align} \label{Eq42}
  \mcO p_{G^\odot} \times^{H^0 (X, \Omega^{\otimes 2})} H^0 (X,  \DV_{G}) \isom  \mcO p_{G}.
 \end{align}
 In particular, the natural  $H^0 (X, \DV_G)$-action on 
 $\mcO p_{G^\odot} \times^{H^0 (X, \Omega^{\otimes 2})} H^0 (X,  \DV_{G})$ yields, via (\ref{Eq42}), an affine structure on $\mcO p_G$ modeled on that space (cf. ~\cite[Theorem 2.24]{Wak8}).
 
 Moreover, if $\rho$ is an element of $\mfc (k)^{\times r}$ (where $\rho := \emptyset$ if $r = 0$),
 then this structure restrict to an affine structure on $\mcO p_{G, \rho} \left(\subseteq \mcO p_G \right)$ modeled on $H^0(X, {^c}\DV_G) \left(\subseteq  H^0(X, \DV_G) \right)$ (cf. ~\cite[Theorem 2.36]{Wak8}).

\LSP
%---------------------------[begin subsection]-------------
\subsection{The $\mr{SL}_{n}$-oper associated to an $\mr{SL}_2$-oper} \label{SS03g61}

Let $n$ be a positive integer with $1 < n < p$ and $\msF^\heartsuit_\odot := (\mcF_\odot, \nabla_\odot, \{ \mcF_\odot^j \}_{j=0}^2)$  an $\mr{SL}_2$-oper on $\msX$ (cf. Remark \ref{Rem667}).
Write $\varTheta := \mcF^1_\odot$, which is a  line bundle on $X$.
Then, we obtain   a composite isomorphism
\begin{align} \label{Eq803}
\Omega \isom \varTheta \otimes_{\mcO_X} (\mcF_\odot/\varTheta)^\vee \isom \varTheta^{\otimes 2},
\end{align}
where the first arrow denotes the isomorphism given by   the $1$-st  Kodaira-Spencer map 
$\mr{KS}_{\msF_\odot^\heartsuit}^1 : \varTheta \isom \Omega \otimes_{\mcO_X} (\mcF_\odot/ \varTheta)$ for  $\msF^\heartsuit$ (cf. (\ref{Eqq212})) and the second arrow 
arises from the fixed isomorphism $\left(\varTheta \otimes_{\mcO_X} (\mcF_\odot/\varTheta) =  \right) \mr{det}(\mcF_\odot) \isom \mcO_X$.
Thus, the line bundle $\varTheta$ together with (\ref{Eq803}) specifies a theta characteristic  of $X^\mr{log}$ in the sense of ~\cite[Example 4.34]{Wak8};
we shall refer to it as the {\bf theta characteristic associated to $\msF^\heartsuit_\odot$}.

The $(n-1)$-st symmetric product $\mr{Sym}^{n-1} \mcF_\odot$ of $\mcF_\odot$ over $\mcO_X$ forms a rank $n$ vector bundle on $X$;
it admits a log connection $\mr{Sym}^{n-1}\nabla_\odot$ induced naturally by $\nabla_\odot$.
Moreover, $\mr{Sym}^{n-1}\mcF_\odot$ has  an $n$-step decreasing filtration $\{ (\mr{Sym}^{n-1}\mcF_\odot)^j \}_{j=0}^n$ defined in such a way that
 $(\mr{Sym}^{n-1}\mcF_\odot)^0 := \mr{Sym}^{n-1}\mcF_\odot$, $(\mr{Sym}^{n-1} \mcF_\odot)^n = 0$, and $(\mr{Sym}^{n-1}\mcF_\odot)^j$ (for each $j =1, \cdots, n-1$) is defined as the image of 
$\varTheta^{\otimes j} \otimes \mcF_\odot^{\otimes (n-1-j)}$ via the natural quotient $\mcF_\odot^{\otimes (n-1)} \migisurj \mr{Sym}^{n-1}\mcF_\odot$.
For each $j = 0, \cdots, n-1$, there exists a canonical composite of isomorphisms
\begin{align} \label{Eq901}
(\mr{Sym}^{n-1}\mcF_\odot)^j /(\mr{Sym}^{n-1}\mcF_\odot)^{j+1} 
\isom \varTheta^{\otimes j} \otimes_{\mcO_X} (\mcF_\odot/\varTheta)^{\otimes (n-1 - j)}
\isom \varTheta^{\otimes (2j -n +1)}.
\end{align}

Moreover, by the assumption $n < p$, the collection
\begin{align} \label{Eq288}
\mr{Sym}^{n-1}\msF_\odot^\heartsuit := (\mr{Sym}^{n-1}\mcF_\odot, \mr{Sym}^{n-1}\nabla_\odot, \{ (\mr{Sym}^{n-1}\mcF_\odot)^j \}_{j=0}^n)
\end{align}
 forms an $\mr{SL}_n$-oper on $\msX$, called the {\bf $(n-1)$-st symmetric product of $\msF^\heartsuit_\odot$}.
The resulting assignment $\msF^\heartsuit_\odot \mapsto \mr{Sym}^{n-1} \msF^\heartsuit_\odot$
determines a morphism
\begin{align} \label{Eq289}
\mcO p_{\mr{SL}_2} \migi \mcO p_{\mr{SL}_n}.
\end{align}
This morphism makes  the following square diagram commute:
\begin{align} \label{Eq345}
\vcenter{\xymatrix@C=46pt@R=36pt{
\mcO p_{\mr{SL}_2} \ar[r]^-{(\ref{Eq289})}  \ar[d] & \mcO p_{\mr{SL}_n}\ar[d]
\\
\mcO p_{G^\odot} \ar[r]_-{(\ref{Eq405})} & \mcO p_{\mr{PGL}_n},
}}
\end{align}
where the right-hand and left-hand vertical arrows are obtained via projectivization.

\LSP
%---------------------------[begin subsection]-------------
\subsection{The $\mr{SL}_{2l+1}$-oper associated to a $\mr{PGL}_2$-oper} \label{SS03f61}

This subsection discusses a construction of an $\mr{SL}_{n}$-oper for an odd $n$   using a $G^\odot$-oper that 
 is compatible with (\ref{Eq405}) and (\ref{Eq289}).
% $\widetilde{\iota}^{\mcO p}_{\mfs \mfl_n}$ 
 %via the projection $\mcO p_{\mr{SL}_2} \migi \mcO p_{G^\odot}$.
  The advantage of considering this construction is that 
  $\msX$ always admits $G^\odot$-opers in contrast to the case of $\mr{SL}_{2}$-opers.
  (In fact, there is a nodal curve that does not have  theta characteristics; this  implies $\mcO p_{\mr{PGL}_2} = \emptyset$ for such a curve.)
%  there can be situations in which 
%$\msX$ does not  have $\mr{SL}_2$-opers, or equivalently, does not have  theta characteristics
%(in contrast to the nonemptiness of  $\mcO p_{\mr{PGL}_2}$).

Let $\msE^\spadesuit_\odot := (\mcE_{B^\odot}, \nabla_\odot)$ be a  $G^\odot$-oper  on $\msX$.
Write $\mcE_\odot := \mcE_{B^\odot} \times^{B^\odot} G^\odot$ and
$\msE_\odot := (\mcE_\odot, \nabla_\odot)$.
We can choose a  $\mr{GL}_2$-oper 
$\msF^\heartsuit_\odot := (\mcF_\odot, \nabla_{\mcF_\odot}, \{ \mcF^j_\odot \}_{j=0}^2)$ on $\msX$ whose projectivization is isomorphic to $\msE^\spadesuit_\odot$ (cf. ~\cite[Theorem 4.66]{Wak8}).
For an integer $l$ with $1 \leq l <  \frac{p-1}{2}$,
the previous discussion can be  applied to obtain a  $\mr{GL}_{2l+1}$-oper 
\begin{align}
\mr{Sym}^{2l}\msF^\heartsuit_\odot := (\mr{Sym}^{2l}\mcF_\odot, \mr{Sym}^{2l}\nabla_{\mcF_\odot}, \{ (\mr{Sym}^{2l}\mcF_\odot)^j \}_{i=0}^{2l+1}).
\end{align}
In particular,  there exists a sequence of natural isomorphisms 
\begin{align} \label{Eq877}
\mr{det}(\mr{Sym}^{2l}\mcF_\odot) &\isom \bigotimes_{j=0}^{2l} (\mr{Sym}^{2l}\mcF_\odot)^j/(\mr{Sym}^{2l}\mcF_\odot)^{j+1} \\
&\isom  \bigoplus_{j=0}^{2l} (\mcF^1_\odot)^{\otimes 2l} \otimes_{\mcO_X} \mcT^{\otimes 2l-j} \notag \\
&\isom
(\mcF^1_\odot)^{\otimes 2l(2l+1)} \otimes_{\mcO_X} \mcT^{\otimes l (2l+1)} \notag \\
& \isom ((\mcF^{1}_\odot)^{\otimes 2l} \otimes_{\mcO_X} \mcT^{\otimes l})^{\otimes (2l+1)}. \notag
\end{align}
Since $p\nmid 2l+1$, there exists a log connection $\nabla_\mcL$ on $\mcL := (\mcF^{1})^{\otimes 2l} \otimes_{\mcO_X} \mcT^{\otimes l}$ whose $(2l+1)$-fold tensor product
 corresponds to $\mr{Sym}^{2l}\nabla_\odot$ via  (\ref{Eq877}) (cf. ~\cite[Proposition 4.22, (i)]{Wak8}).
Let us  write 
\begin{align}
\mr{Sym}^{2l}\mcE_\odot := \mcL^\vee \otimes_{\mcO_X} \mr{Sym}^{2l} \mcF_\odot \ \ \ \text{and} \ \   \ (\mr{Sym}^{2l}\mcE_\odot)^j := \mcL^\vee \otimes_{\mcO_X} (\mr{Sym}^{2l} \mcF_\odot)^j
\end{align}
($j=0, \cdots, 2l+1$).
Then,
there exist canonical isomorphisms
\begin{align} \label{Eq913}
(\mr{Sym}^{2l}\mcE_\odot)^j /(\mr{Sym}^{2l}\mcE_\odot)^{j+1} \isom \Omega^{\otimes j-l} 
\end{align}
($j=0, \cdots, 2l$).
If $\nabla^\vee_\mcL$ denotes the log connection on the dual $\mcL^\vee$ induced naturally by $\nabla_\mcL$,
the log connection 
$\mr{det}(\nabla_\mcL^\vee \otimes \mr{Sym}^{2l}\nabla_{\mcF_\odot})$
on the determinant  $\mr{det}(\mr{Sym}^{2l}\mcE_\odot)$ is compatible with the trivial connection $d$ on $\mcO_X$ via
the composite isomorphism
\begin{align}
\mr{det}(\mr{Sym}^{2l}\mcE_\odot) \isom \bigotimes_{j=0}^{2l} 
(\mr{Sym}^{2l}\mcE_\odot)^j /(\mr{Sym}^{2l}\mcE_\odot)^{j+1}
\isom \bigotimes_{j=0}^{2l}  \Omega^{\otimes j-l} \isom \mcO_X.
\end{align}
Write  $\mr{Sym}^{2l}\nabla_\odot$ for the log connection 
$\nabla_\mcL^\vee \otimes \mr{Sym}^{2l}\nabla_{\mcF_\odot}$ on
$\mr{Sym}^{2l}\mcE_\odot$.

Then, the resulting pair $\mr{Sym}^{2l}\msE_\odot := (\mr{Sym}^{2l}\mcE_\odot, \mr{Sym}^{2l}\nabla_\odot)$ specifies a flat bundle, and 
moreover, 
the collection of data
\begin{align} \label{Eq9d00}
\mr{Sym}^{2l}\msE_\odot^\spadesuit 
:= (\mr{Sym}^{2l}\mcE_\odot, \mr{Sym}^{2l}\nabla_\odot, \{ (\mr{Sym}^{2l}\mcE_\odot)^j \}_{j=0}^{2l+1})
\end{align}
has  a structure of 
 $\mr{SL}_{2l+1}$-oper.
 The isomorphism classes of $\mr{Sym}^{2l}\msE_\odot$ and $\mr{Sym}^{2l}\msE_\odot^\spadesuit$ depend only on 
 $\msE^\spadesuit_\odot$ (i.e,, not on the choice of $\msF^\heartsuit_\odot$).
 We refer to $\mr{Sym}^{2l}\msE_\odot^\spadesuit$ as 
   the {\bf $2l$-th symmetric product of $\msE_\odot^\spadesuit$}.
 
 The resulting assignment $\msE^\spadesuit_\odot  \mapsto \mr{Sym}^{2l}\msE^\spadesuit_\odot$ determines a well-defined  morphism
\begin{align} \label{Eq903}
\mcO p_{G^\odot} \migi \mcO p_{\mr{SL}_{2l+1}},
\end{align}
which makes the following diagram commute:
\begin{align} \label{Eq34g5}
\vcenter{\xymatrix@C=46pt@R=36pt{
\mcO p_{\mr{SL}_2} \ar[r]^-{(\ref{Eq289})}  \ar[d]_-{\mr{projectivization}} & \mcO p_{\mr{SL}_{2l+l}}\ar[d]^-{\mr{projectivization}}
\\
\mcO p_{G^\odot} \ar[r]_-{(\ref{Eq405})} \ar[ur]^-{(\ref{Eq903})} & \mcO p_{\mr{PGL}_{2l+1}}.
}}
\end{align}

%%%%%%%%%%%%%%%%%%%%%%%%%%%%%%%%%%%%%%%%%%%%%
%%%%%%%%%%%%%---[begin section]---%%%%%%%%%%%%%%%%%%%%%%%
\vspace{10mm}
\section{Infinitesimal deformations of opers}\label{S230}\SSP

This section deals with  infinitesimal deformations  of a flat $G$-bundle, as well as a $G$-oper.
In particular,  we verify dualities between certain deformation spaces of a flat $G$-bundle with prescribed residues/radii.
These results 
 are probably already known at least in the case of characteristic $0$ or in some restricted situations. However, 
 the author could not find any suitable  references  consistent with our situation, so we here  discuss them.

\LSP
%---------------------------[begin subsection]-------------
\subsection{De Rham cohomology/Parabolic de Rham cohomology} \label{SSa22}

Let $\nabla' : \mcK^0 \migi \mcK^1$ be 
 a morphism of sheaves on $X$.
 It may be regarded as a complex concentrated at degrees $0$ and $1$;  this complex is denoted by 
$\mcK^\bullet [\nabla']$.
(In particular, we have  $\mcK^j [\nabla'] := \mcK^j$ for $j= 0,1$).
For each $l \geq 0$, we shall write $\mbH^l (X, \mcK^\bullet [\nabla'])$ for the $l$-th hypercohomology group of the complex $\mcK^\bullet [\nabla']$.

Given an integer $l$ and a sheaf $\mcF$ on $X$,
we define the complex $\mcF [l]$ to be $\mcF$ (considered as a complex concentrated at degree $0$) shifted down by $l$, so that $\mcF [l]^{-l} = \mcF$ and $\mcF [l]^i = 0$ for $i \neq l$. 

For a flat bundle $\msF := (\mcF, \nabla)$  on $\msX$, we shall set
$\Omega^{0}_{\mr{par}} (\mcF) := \mcF$ and  
\begin{align} \label{Eq651}
\Omega^{1}_{\mr{par}}(\mcF) := \Omega \otimes_{\mcO_X}  {^c}\mcF + \mr{Im}(\nabla) \left(\subseteq \Omega \otimes_{\mcO_X} \mcF\right).
\end{align}
Then, $\nabla$ restricts to a $k$-linear morphism
\begin{align} \label{Eq650}
\nabla_{\mr{par}} : \Omega^{0}_{\mr{per}}(\mcF) \migi \Omega^{1}_{\mr{per}} (\mcF).
\end{align}
The $k$-vector space
\begin{align} \label{Eq5010}
H^l_{\mr{dR}} (X, \msF) := \mbH^l (X, \mcK^\bullet [\nabla])  \ \left(\text{resp.,} \ 
 H^l_{\mr{dR}, \mr{par}} (X, \msF) := \mbH^l (X, \mcK^\bullet [\nabla_{\mr{par}}])
 \right)
\end{align}
is called the {\bf $l$-th de Rham cohomology group} (resp., the {\bf $l$-th  parabolic de Rham cohomology group}) of $\msF$.
The $k$-linear morphism  $\mbH^1 (X, \mcK^\bullet [\nabla_\mr{par}]) \migi \mbH^1 (X, \mcK^\bullet [\nabla])$ induced by 
the inclusion of complexes $\mcK^\bullet [\nabla_\mr{par}] \migi \mcK^\bullet [\nabla]$ is  injective.
By using this injection, we occasionally regard $\mbH^1 (X, \mcK^\bullet [\nabla_\mr{par}])$ as a subspace of $\mbH^1 (X, \mcK^\bullet [\nabla])$.

Here, 
we prove  a comparison between the parabolic de Rham cohomology and the hypercohomology of 
the restriction 
\begin{align} \label{Eqqs}
{^c}\nabla_{\mr{par}} : 
\Omega_\mr{par}^{0, c} (\mcF) \migi \Omega_\mr{par}^{1, c} (\mcF) 
\end{align}
of $\nabla$, where 
$\Omega_\mr{par}^{0, c} (\mcF) := \nabla^{-1}(\Omega\otimes_{\mcO_X} {^c}\mcF)$ and $\Omega_\mr{par}^{1, c} (\mcF) := \Omega \otimes_{\mcO_X} {^c}\mcF$.

\SSP
%-------------------------------------------------------------------------------------
\bpr \label{Prp7}
For each $l \geq 0$,
the morphism of $k$-vector spaces
\begin{align} \label{Eq8010}
\mbH^l (X, \mcK^\bullet [{^c}\nabla_{\mr{par}}]) \isom
\mbH^l (X, \mcK^\bullet [\nabla_{\mr{par}}])
\end{align}
induced by the inclusion $\mcK^\bullet [{^c}\nabla_{\mr{par}}] \migiincl  \mcK^\bullet [\nabla_{\mr{par}}]$ is an isomorphism.
\epr
%-------------------------------------------------------------------------------------
\begin{proof}
Let us choose   $i \in \{1, \cdots, r \}$.
Then, $\sigma^*_i (\Omega \otimes_{\mcO_X} \mcF)$ may be identified with $\sigma^*_i (\mcF)$ via 
$\sigma_i^*(\Omega \otimes_{\mcO_X} \mcF) \cong \sigma_i^* (\Omega) \otimes_k \sigma_i^* (\mcF) \cong \sigma_i^*(\mcF)$, where the second ``$\cong$" follows from the residue map $\sigma^*_i (\Omega) \isom k$.
Under this identification, 
the log connection $\nabla$ specifies   an endomorphism
$\nu_i : \sigma^*_i (\mcF) \migi \sigma_i^* (\mcF)  \left(= \sigma_i^*(\Omega \otimes_{\mcO_X} \mcF) \right)$ of $\sigma_i^* (\mcF)$.
By taking quotients of the both sides of  the equality  $\Omega_{\mr{par}}^0 (\mcF)/{^c}\mcF = \bigoplus_{i=1}^r \sigma_{i*}(\sigma_i^*(\mcF))$, 
we obtain 
\begin{align} \label{Eq8012i}
\Omega_{\mr{par}}^0 (\mcF) / \Omega_{\mr{par}}^{0, c} (\mcF) = \bigoplus_{i=1}^r \sigma_{i*}(\sigma_i^*(\mcF)/\mr{Ker}(\nu_i)).
\end{align}
Also, the identification $(\Omega \otimes_{\mcO_X} \mcF)/\Omega_\mr{par}^{1, c} (\mcF)  = \left( \bigoplus_{i=1}^r \sigma_{i*} (\sigma^*_i (\Omega) \otimes_k \sigma^*_i (\mfg_\mcE)) =\right) \bigoplus_{i=1}^r \sigma_{i*}(\mcF)$
restricts to 
\begin{align} \label{Eq222ii}
\Omega_{\mr{par}}^1 (\mcF)/ \Omega_{\mr{par}}^{1, c} (\mcF)
= \bigoplus_{i=1}^r \sigma_{i*}(\mr{Im}(\nu_i)).
\end{align}
By using (\ref{Eq8012i}) and (\ref{Eq222ii}), we obtain a morphism of short exact sequences
\begin{align} \label{Eqq300}
\vcenter{\xymatrix@C=46pt@R=36pt{
0 \ar[r] & \Omega_\mr{par}^{0, c} (\mcF) \ar[r]^-{\mr{inclusion}} \ar[d]^-{{^c}\nabla_{\mr{par}}} & \Omega_{\mr{par}}^0 (\mcF) \ar[r] \ar[d]^-{\nabla_\mr{par}} & \bigoplus_{i=1}^r \sigma_{i*}(\sigma^*(\mcF) / \mr{Ker}(\nu_i)) \ar[r] \ar[d]^-{\bigoplus_{i=1}^r \sigma_{i*}(\nu_i)}& 0
\\
0 \ar[r] & \Omega^{1, c}_\mr{par} (\mcF) \ar[r]_-{\mr{inclusion}}  & \Omega_{\mr{par}}^1 (\mcF) \ar[r]& \bigoplus_{i=1}^r \sigma_{i*}(\mr{Im}(\nu_i)) \ar[r] & 0. 
}}
\end{align}
Since the right-hand vertical arrow in this diagram  is an isomorphism,
the inclusion   $\mcK^\bullet [{^c}\nabla_{\mr{par}}] \migiincl \mcK^\bullet [\nabla_{\mr{par}}]$ turns out to be  a quasi-isomorphism.
This  means  that
the induced  morphism $\mbH^1 (X, \mcK^\bullet [{^c}\nabla_{\mr{par}}]) \migi 
\mbH^1 (X, \mcK^\bullet [\nabla_{\mr{par}}])$  is  an isomorphism.  
\end{proof}
%-------------------------------------------------------------------------------------
%\SSP

\LSP
%---------------------------[begin subsection]-------------
\subsection{Duality for  de Rham cohomology} \label{SS055}

Now, let $(\mcE, \nabla)$ be a flat $G$-bundle on $\msX$.
The log connection $\nabla$ induces a log connection 
$\nabla^\mr{ad} : \mfg_{\mcE} \migi \Omega\otimes_{\mcO_X} \mfg_{\mcE}$
 on the adjoint  bundle $ \mfg_{\mcE}$;
 it restricts to a log connection
\begin{align} \label{Eq454}
{^c} \nabla^\mr{ad} : {^c}\mfg_\mcE  \migi \Omega \otimes_{\mcO_X} {^c}\mfg_\mcE 
\end{align}
 on ${^c}\mfg_\mcE \left(= \mfg_\mcE (-D)\right)$.
 If ${^c}d$ denotes the restriction  ${^c}\mcO_X  \migi \Omega \otimes_{\mcO_X} {^c}\mcO_X  \left(= \omega_X \right)$ of the universal logarithmic derivation $d$,
 then
 the pairing  $({^c}\mfg_\mcE, {^c}\nabla^\mr{ad}) \otimes_{\mcO_X} (\mfg_\mcE, \nabla^\mr{ad})\migi ({^c}\mcO_X, {^c}d)$ arising from  the Killing form $\kappa$ on $\mfg$ yields 
a composite $k$-linear pairing
\begin{align} \label{Eq24012}
\mbH^1 (X, \mcK^\bullet [{^c}\nabla^\mr{ad}]) \times \mbH^1 (X, \mcK^\bullet [\nabla^\mr{ad}]) \xrightarrow{\cup} \mbH^2 (X, \mcK^\bullet [{^c}d]) \isom H^1 (X, \omega_X) \isom k.
\end{align}
The following assertion for the case of flat $G$-bundles with vanishing $p$-curvature can be found in ~\cite[Corollary 6.16]{Wak8}.

\SSP
%--------------------------------------------------------------------------------------
\bpr \label{Prop90}
The  pairing  (\ref{Eq24012})
is nondegenerate.
That is to say, 
the $k$-linear morphism
\begin{align} \label{Eww3}
\mbH^1 (X, \mcK^\bullet [{^c}\nabla^\mr{ad}]) \migi \mbH^1 (X, \mcK^\bullet [\nabla^{\mr{ad}}])^\vee
\end{align}
determined  by (\ref{Eq24012})
 is an isomorphism.
\epr
%--------------------------------------------------------------------------------------
\begin{proof}
Let us consider the diagram
\begin{align} \label{Eq3i90}
\vcenter{\xymatrix@C=7pt@R=36pt{
H^0 (X, {^c}\mfg_\mcE) \ar[r] \ar[d] & H^0 (X, \Omega\otimes_{\mcO_X} {^c}\mfg_\mcE) \ar[r] \ar[d] & \mbH^1 (X, \mcK^\bullet [{^c}\nabla^\mr{ad}]) \ar[r] \ar[d]^-{(\ref{Eww3})} &H^1 (X, {^c}\mfg_\mcE)  \ar[r] \ar[d] &H^1 (X, \Omega\otimes_{\mcO_X} {^c}\mfg_\mcE) \ar[d]  \\
H^1 (X, \Omega\otimes_{\mcO_X} \mfg_\mcE)^\vee \ar[r] & H^1 (X, \mfg_\mcE)^\vee \ar[r] &\mbH^1 (X, \mcK^\bullet [\nabla^\mr{ad}])^\vee \ar[r] &
H^0 (X, \Omega\otimes_{\mcO_X} \mfg_\mcE)^\vee \ar[r]  & H^0 (X, \mfg_\mcE)^\vee,
}}
\end{align}
where the upper and lower horizontal sequences are the exact sequences arising from 
the Hodge to de Rham spectral sequences 
$H^b (X, \mcK^a[{^c}\nabla^\mr{ad}]) \Rightarrow \mbH^{a+b}(X, \mcK^\bullet[{^c}\nabla^\mr{ad}])$ and $H^b (X, \mcK^a[\nabla^\mr{ad}]) \Rightarrow \mbH^{a+b}(X, \mcK^\bullet[\nabla^\mr{ad}])$,
respectively. 
The assumed condition on $G$  described  in \S\,\ref{Ssde} implies that the Killing form $\kappa$ is nondegenerate, so it induces an isomorphism $\mfg_\mcE \isom \mfg_\mcE^\vee$.
All the vertical arrows in (\ref{Eq3i90}) except for the middle one  are isomorphisms because they  come from the pairings given by  Serre duality.
It follows that 
  (\ref{Eww3})  is  an isomorphism by the five lemma, and this completes the proof of this assertion.
\end{proof}
%--------------------------------------------------------------------------------------
\SSP

Denote by $q$ the $k$-rational point of $\C_G$ classifying 
$(\mcE, \nabla)$.
Also, let
 $T_q \C_G$ (resp., $T_q^\vee \C_G$) denote the tangent space (resp., the cotangent space)  of $\C_G$  at $q$.
As discussed in   ~\cite[\S\,6.1.4]{Wak8} (or, e.g., ~\cite[Proposition 3.6]{O3} for the case where $\msX$ is unpointed and smooth), 
 the  underlying set of  $\mbH^1 (X, \mcK^\bullet [\nabla^\mr{ad}])$
   may be identified with the space of first-order deformations of  $(\mcE, \nabla)$.  That is to say,  there exists a canonical  isomorphism of $k$-vector spaces
\begin{align} \label{Eq3033}
\gamma : T_q \C_G \isom \mbH^1 (X, \mcK^\bullet [\nabla^\mr{ad}]).
\end{align}

If $(\mcE, \nabla)$ is equipped with a structure of marking and its residues are given by   $\mu := (\mu_i )_{i=1}^r \in \mfg^{\times r}$ (where $\mu := \emptyset$ if $r = 0$),
then
a similar discussion  yields a canonical isomorphism of $k$-vector spaces
\begin{align} \label{Eq350}
\gamma_\mu : T_q \C_{G, \mu} \isom \mbH^1 (X, \mcK^\bullet [{^c}\nabla^\mr{ad}]),
\end{align}
where
 we use the same notation ``$q$" to denote
 the $k$-rational point of  $\C_{G, \mu}$ classifying $(\mcE, \nabla)$ (with marking).
In particular,  by using the isomorphisms (\ref{Eww3}), (\ref{Eq3033}), and (\ref{Eq350}),
we obtain an isomorphism
\begin{align} \label{Er67}
T_q \C_{G, \mu}  \isom  T_q^\vee \C_G.
\end{align}

\LSP
%---------------------------[begin subsection]-------------
\subsection{Duality for  parabolic de Rham cohomology} \label{SS055g}

Next, 
let us consider the case of parabolic de Rham cohomology.

We shall denote by 
\begin{align} \label{Eqq338}
\kappa^\triangleright : \left((\Omega \otimes_{\mcO_X} {^c}\mfg_\mcE) \times \mfg_\mcE =  \right) \Omega^{1, c}_\mr{par} (\mfg_\mcE) \times \Omega_\mr{par}^0 (\mfg_\mcE) \migi {^c}\Omega
\end{align}
 the $\mcO_X$-bilinear pairing induced from the Killing form $\kappa$ of $\mfg$.
 Since $\kappa$ is nondegenerate, we see that $\kappa^\triangleright$ is nondegenerate.
By applying  Serre  duality to $\kappa^\triangleright$, we obtain a nondegenerate pairing of $k$-vector spaces
 \begin{align} \label{Eqw5}
 H^{1-l} (X, \Omega_\mr{par}^{1, c} (\mfg_\mcE)) \otimes_k H^l (X, \Omega^0_\mr{par}(\mfg_\mcE)) \migi \left(H^1 (X, {^c}\Omega) = \right) k
 \end{align}
 for  $l =0,1$.

Next, we shall write $U := X \setminus \bigcup \{ \sigma_i \}_{i=1}^r$ and  write $\iota$ for the inclusion $U \migiincl X$.
The $\mcO_X$-module $ \iota_*(\iota^* (\mfg_\mcE))$  (resp., $\iota_*(\iota^* (\Omega \otimes_{\mcO_X} \mfg_\mcE))$) contains $\Omega^{0, c}_\mr{par}(\mfg_\mcE)$ (resp., $\Omega^1_\mr{par}(\mfg_\mcE)$) as an $\mcO_X$-submodule. 
Here, let us consider the  natural 
skew-symmetric $\mcO_X$-bilinear pairing
\begin{align}  \label{Eqw4}
\widetilde{\kappa}^\blacktriangleright : \iota_*(\iota^* (\mfg_\mcE)) \times \iota_*(\iota^* (\Omega \otimes_{\mcO_X} \mfg_\mcE)) \migi \iota_*(\iota^*(\Omega))
\end{align}
arising from  $\kappa$.

\SSP
%-------------------------------------------------------------------------------------
\bpr \label{Prop99}
Assume (cf. Remark \ref{Rem55}) that 
there exists a $B$-reduction $\mcE_B$ of $\mcE$ for which the pair $\msE^\spadesuit := (\mcE_B, \nabla)$ forms a $G$-oper 
  of radii $[0]^{\times r} := ([0], [0], \cdots, [0]) \in \mfc (k)^{\times r}$ (where $[0]^{\times r} := \emptyset$ if $r = 0$).
Then,  $\widetilde{\kappa}^\blacktriangleright$
restricts to 
a nondegenerate $\mcO_X$-bilinear pairing
\begin{align} \label{Eqw3}
\kappa^\blacktriangleright : \Omega^{0, c}_\mr{par}(\mfg_\mcE) \times \Omega^1_\mr{par}(\mfg_\mcE) \migi {^c}\Omega.
\end{align}
In particular, for $l =0, 1$,  the morphism 
\begin{align} \label{Eqw1}
H^{1-l} (X, \Omega^{0, c}_\mr{par}(\mfg)) \otimes_k H^l (X, \Omega^1_\mr{par}(\mfg_\mcE)) \migi \left(H^1 (X, {^c}\Omega) =  \right) k
\end{align}
obtained by applying Serre duality to this pairing is nondegenerate.
\epr
%-------------------------------------------------------------------------------------
\begin{proof}
After replacing $\msE^\spadesuit$ with its $\,\qq$-normalization (cf. ~\cite[Proposition 2.18]{Wak8}), we may assume that
$\msE^\spadesuit$ is $\qq$-normal.
In particular, since $\nabla$ is of radii $[0]^{\times r}$,
its residue at $\sigma_i$ coincides with $q_{-1} \in \mfg$ for every $i=1, \cdots, r$.

We first consider the former assertion.
Just as in the case of $\kappa^\triangleright$, 
 the restriction of $\widetilde{\kappa}^\blacktriangleright$ to $U$ is nondegenerate.
Thus, the problem is  reduced to examining  the pairing $\widetilde{\kappa}^\blacktriangleright$  around   the marked points.
Let us choose $i \in \{1, \cdots, r \}$ and choose a local function  $t \in \mcO_X$ defining the closed subscheme  $\mr{Im}(\sigma_i) \subseteq X$.
Then, the formal neighborhood $Q$ of $\sigma_i$ is naturally isomorphic to $\mr{Spec}(k[\![t]\!])$, and we have  ${^c}\Omega |_Q = k[\![t]\!] dt$.
The choice of such a function $t$ specifies a canonical identification $\mfg_\mcE |_Q = k[\![t]\!] \otimes_k \mfg$ 
(cf. ~\cite[\S\,2.4.3]{Wak8}), which also gives
 $\iota_* (\iota^* (\mfg_\mcE)) |_{Q} = k (\!(t)\!) \otimes_k \mfg$.
Under these identifications,
the restriction of  $\widetilde{\kappa}^\blacktriangleright$ to $Q$ coincides with the Killing form $\kappa$ tensored with $k(\!(t)\!)$.
Here, recall  (cf. ~\cite[\S\,2.4.4 and  \S\,6.2.3]{Wak8}) that there exist 
direct sum 
decompositions of $\mfg$: 
\begin{align} \label{Eqq721}
\mfg = \mr{Ker}(\mr{ad}(q_{-1})) \oplus \mr{Im}(\mr{ad}(q_1)) = \mr{Ker}(\mr{ad}(q_{1})) \oplus \mr{Im}(\mr{ad}(q_{-1})),
\end{align}
and that the isomorphism $\mfg \isom \mfg^\vee$ arising from $\kappa$ 
 restricts to $\mr{Ker}(\mr{ad}(q_{-1})) \isom \mr{Ker}(\mr{ad}(q_{1}))^\vee$, as well as to  $\mr{Im}(\mr{ad}(q_1))\isom \mr{Im}(\mr{ad}(q_{-1}))^\vee$.
In particular, by
the equality  $\mu_i^\nabla = q_{-1}$, 
  the first (resp., second) decomposition of (\ref{Eqq721})  yields
\begin{align} \label{Eqq291}
\Omega_\mr{par}^{0, c} (\mfg_\mcE) |_Q = k[\![t]\!]  \otimes_k\mr{Ker}(\mr{ad}(q_{-1})) \oplus k[\![t]\!] t\otimes_k\mr{Im}(\mr{ad}(q_1))\hspace{11mm} \\
\left(\text{resp.,} \ \Omega_\mr{par}^{1} (\mfg_\mcE) |_Q = k[\![t]\!] dt \otimes_k\mr{Ker}(\mr{ad}(q_{1})) \oplus k[\![t]\!] \frac{dt}{t}\otimes_k\mr{Im}(\mr{ad}(q_{-1}))\right). \notag
\end{align}
 According to this description, 
 each  section $v$ (resp., $u$) of 
$\Omega_\mr{par}^{0, c} (\mfg_\mcE) |_Q$ (resp., $\Omega_\mr{par}^{1} (\mfg_\mcE) |_Q$) may be described as the sum 
$v = v_1 + t \cdot v_2$ (resp., $u = u_1 + \frac{1}{t} \cdot u_2$) for some $v_1 \in k[\![t]\!]  \otimes_k\mr{Ker}(\mr{ad}(q_{-1}))$,
$v_2 \in k[\![t]\!]\otimes_k\mr{Im}(\mr{ad}(q_1))$ (resp., $u_1 \in k[\![t]\!] dt \otimes_k\mr{Ker}(\mr{ad}(q_{1}))$, $u_2 \in k[\![t]\!] dt \otimes_k\mr{Im}(\mr{ad}(q_{-1}))$). 
Then, we have 
\begin{align}
\widetilde{\kappa}^\blacktriangleright (v, u) &=
\widetilde{\kappa}^\blacktriangleright (v_1, u_1) + \widetilde{\kappa}^\blacktriangleright (t \cdot v_2, \frac{1}{t} \cdot u_2) \\
 &= \widetilde{\kappa}^\blacktriangleright (v_1, u_1) + \widetilde{\kappa}^\blacktriangleright (v_2, u_2) \notag \\
 &= \widetilde{\kappa}^\blacktriangleright (v_1 + v_2, u_1 + u_2). \notag
\end{align} 
It follows that  $\widetilde{\kappa}^\blacktriangleright (v, u)$ lies in  $k[\![t]\!] dt \left(= {^c}\Omega |_Q\right)$, and 
 the restriction of  $\widetilde{\kappa}^\blacktriangleright$  over $Q$ is nondegenerate.
This completes the proof of the former assertion.

The latter assertion follows directly  from the former assertion.
\end{proof}
%-------------------------------------------------------------------------------------
\SSP

%-------------------------------------------------------------------------------------
\begin{rema}[Flat $G$-bundles with parabolic structure] \label{Rem55}
By  the similarity of the proof, 
the same assertion as Proposition \ref{Prop99} remains  true after replacing 
the assumption (i.e., the existence of ``$\mcE_B$") with the following condition (regarded as  a kind of parabolic structure):
the flat $G$-bundle $(\mcE, \nabla)$ is equipped with a marking, and its residue at  $\sigma_i$ coincides with $q_{-1} \in \mfg$ (via the fixed  marking) for every $i=1, \cdots, r$.
\end{rema}
%-------------------------------------------------------------------------------------
\SSP

Until the end of this subsection, we keep the assumption in the above proposition.
Let $\Box$ denote either the absence or presence of ``$c\,$".
 One may calculate $\mbH^1 (X, \mcK^\bullet [{^\Box}\nabla_\mr{par}^\mr{ad}])$
as the total cohomology of the \v{C}ech double complex $\mr{Tot}^\bullet (\check{C} (\msU, \mcK^\bullet [{^\Box}\nabla_\mr{par}^\mr{ad}]))$ associated to $\mcK^\bullet [{^\Box}\nabla_\mr{par}^\mr{ad}]$, where $\msU := \{ U_\alpha \}_{\alpha \in I}$ is a finite affine open covering  of $X$.
Denote by   $I_2$ the set of pairs $(\alpha, \beta) \in I \times I$ with $U_{\alpha \beta} := U_\alpha \cap U_\beta \neq \emptyset$.
Then, each element $v$ of $\mbH^1 (X, \mcK^\bullet [{^\Box}\nabla_\mr{par}^\mr{ad}])$ may be given by a collection of data
\begin{align} \label{Eqq321}
v = (\{ \partial_{\alpha \beta} \}_{(\alpha, \beta)}, \{ \delta_\alpha \}_{\alpha})
\end{align}
consisting of a \v{C}ech $1$-cocycle $\{ \partial_{\alpha \beta}\}_{(\alpha, \beta) \in I_2} \in \check{C}^1 (\msU, \Omega_\mr{par}^{0, \Box} (\mfg_\mcE))$ with $\partial_{\alpha \beta} \in H^0 (U_{\alpha \beta}, \Omega_\mr{par}^{0, \Box} (\mfg_\mcE))$ and a \v{C}ech $0$-cochain $\{ \delta_\alpha \}_{\alpha \in I} \in \check{C}^0 (\msU, \Omega^{1, \Box}_\mr{par}(\mfg_\mcE))$ with $\delta_\alpha \in H^0 (U_\alpha, \Omega_{\mr{par}}^
{1, \Box} (\mfg_\mcE))$ which agree under ${^\Box}\nabla_\mr{par}^\mr{ad}$ and the \v{C}ech coboundary map.

Using the description in terms of \v{C}ech double complexes,
one can obtain a skew-symmetric $k$-bilinear pairing 
\begin{align} \label{Eqq1}
 \mbH^1 (X, \mcK^\bullet [{^c}\nabla^\mr{ad}_\mr{par}])  \times \mbH^1 (X, \mcK^\bullet [\nabla^\mr{ad}_\mr{par}])  \migi 
\left(H^1 (X, {^c}\Omega)  =\right) k
\end{align}
given by assigning $(\{ \partial_{\alpha \beta} \}_{(\alpha, \beta)}, \{ \delta_\alpha \}_{\alpha}) \otimes (\{ \partial'_{\alpha \beta} \}_{(\alpha, \beta)}, \{ \delta'_\alpha \}_{\alpha}) \mapsto \kappa^\blacktriangleright  (\partial_{\alpha \beta}, \delta'_\beta) -  \kappa^{\triangleright} (\delta_\alpha, \partial'_{\alpha \beta})$.

\SSP
%--------------------------------------------------------------------------------
\bpr \label{Prop66}
Let us keep the assumption in Proposition \ref{Prop99} (or impose the condition described in Remark \ref{Rem55} instead).
Then, the pairing 
(\ref{Eqq1})
  is nondegenerate.
In particular, 
the morphism of  $k$-vector spaces
\begin{align} \label{Eqq233}
 \mbH^1 (X, \mcK^\bullet [{^c}\nabla^\mr{ad}_\mr{par}]) \migi \mbH^1 (X, \mcK^\bullet [\nabla^\mr{ad}_\mr{par}])^\vee
\end{align}
induced by (\ref{Eqq1}) is an isomorphism.
\epr
%--------------------------------------------------------------------------------
\begin{proof}
Note that 
the morphism  (\ref{Eqq233}) fits into 
 the  morphism of exact sequences
\begin{align} \label{Eqq301}
\vcenter{\xymatrix@C=7pt@R=36pt{
H^0 (X, \Omega_\mr{par}^{0, c}(\mfg_\mcE)) \ar[r] \ar[d]^-{\wr} & H^0 (X, \Omega_\mr{par}^{1, c} (\mfg_\mcE)) \ar[r] \ar[d]^-{\wr} & \mbH^1 (X, \mcK^\bullet [{^c}\nabla_\mr{par}^\mr{ad}]) \ar[r] \ar[d]^-{ (\ref{Eqq233})} & H^1 (X, \Omega_\mr{par}^{0, c}(\mfg_\mcE)) \ar[d]^-{\wr} \ar[r] & H^1 (X, \Omega_\mr{par}^{1, c} (\mfg_\mcE)) \ar[d]^-{\wr}
\\
H^1 (X, \Omega_\mr{par}^1 (\mfg_\mcE))^\vee \ar[r] &H^1 (X, \Omega_\mr{par}^0 (\mfg_\mcE))^\vee \ar[r]& \mbH^1 (X, \mcK^\bullet [\nabla_\mr{par}^\mr{ad}])^\vee \ar[r] & H^0 (X, \Omega_\mr{par}^1 (\mfg_\mcE))^\vee \ar[r] & H^0 (X, \Omega_\mr{par}^0 (\mfg_\mcE))^\vee, 
}}
\end{align}
where 
\begin{itemize}
\item
the upper and lower horizontal sequences  are the exact sequences induced  from 
the Hodge to de Rham spectral sequences 
$E_1^{a, b} = H^b (X, \mcK^a[{^c}\nabla^\mr{ad}_\mr{par}]) \Rightarrow \mbH^{a+b}(X, \mcK^\bullet[{^c}\nabla^\mr{ad}_\mr{par}])$ and $E_1^{a, b} = H^b (X, \mcK^a[\nabla_\mr{par}^\mr{ad}]) \Rightarrow \mbH^{a+b}(X, \mcK^\bullet[\nabla^\mr{ad}_\mr{par}])$,
respectively;
\item
the second and fifth  vertical arrows from the left are the isomorphisms arising from the nondegenerate pairing (\ref{Eqw5}) for $l=1$ and $0$, respectively;
\item
the first and fourth  vertical arrows from the left are the isomorphisms arising from the nondegenerate pairing (\ref{Eqw1}) for $l=1$ and $0$, respectively.
\end{itemize}
Thus, by applying the five lemma to this diagram, we see that the middle vertical arrow  (\ref{Eqq233}) is an isomorphism.
\end{proof}
%--------------------------------------------------------------------------------
\SSP

%--------------------------------------------------------------------------------
\bco  \label{Coo78}
There exists a canonical nondegenerate skew-symmetric pairing
\begin{align} \label{HHH89}
 \mbH^1 (X, \mcK^\bullet [\nabla^\mr{ad}_\mr{par}]) \times 
\mbH^1 (X, \mcK^\bullet [\nabla^\mr{ad}_\mr{par}]) \migi k.
\end{align}
In particular, we obtain  an isomorphism of $k$-vector spaces
\begin{align} \label{Eqqw9}
 \mbH^1 (X, \mcK^\bullet [\nabla^\mr{ad}_\mr{par}]) \isom  \mbH^1 (X, \mcK^\bullet [\nabla^\mr{ad}_\mr{par}])^\vee.
\end{align}
\eco
%--------------------------------------------------------------------------------
\begin{proof}
The required pairing (\ref{HHH89}) is obtained by composing 
(\ref{Eqq1})
  and the inverse of  (\ref{Eq8010}) in the case where ``$\nabla$" is taken to be $\nabla^\mr{ad}$. 
\end{proof}
%--------------------------------------------------------------------------------
\SSP

Next, 
denote by $T_q \C_{G, [0]^{\times r}}$ (resp., $T_q^\vee \C_{G, [0]^{\times r}}$) the tangent space (resp., the cotangent space) of 
$\C_{G, \rho}$ at $q$.
In particular,  $T_q \C_{G, [0]^{\times r}}$ may be identified with  the set of deformations of $(\mcE, \nabla)$ fixing the radii.
Then, we can prove 
 the following assertion.

\SSP
%------------------------------------------------------------------------------------
\bpr \label{Prop556}
Let us keep the assumption in Proposition \ref{Prop99}.
Then, there exists a canonical isomorphism of $k$-vector spaces
\begin{align} \label{K02}
\gamma_{[0]^{\times r}} : T_q \C_{G, [0]^{\times r}} \isom \mbH^1 (X, \mcK^\bullet [\nabla_{\mr{par}}^{\mr{ad}}]).
\end{align}
Moreover,   the following diagram is commutative:
\begin{align} \label{Eq340}
\vcenter{\xymatrix@C=46pt@R=36pt{
    T_q \C_{G, [0]^{\times r}} \ar[r]^-{\mr{inclusion}} \ar[d]^-{\wr}_-{\gamma_{[0]^{\times r}}} &  T_q \C_{G}  \ar[d]^-{\gamma}_-{\wr}
 \\
    \mbH^1 (X, \mcK^\bullet [\nabla_{\mr{par}}^{\mr{ad}}])   \ar[r]_-{\mr{inclusion}} &  \mbH^1 (X, \mcK^\bullet [\nabla^{\mr{ad}}]).
}}
\end{align}
\epr
%------------------------------------------------------------------------------------
\begin{proof}
After possibly replacing $\msE^\spadesuit$ with  its $\,\qq$-normalization,  we may assume  that $\msE^\spadesuit$ is $\,\qq$-normal.
In particular, the assumption implies that the residue of $\nabla$ at every marked point coincides with $q_{-1}$ under the canonical local trivialization of $\mcE \left(= {^\dagger}\mcE_G \right)$  (cf. ~\cite[\S\,2.4.3]{Wak8}).
Let us take a deformation $(\mcE_\varepsilon, \nabla_\varepsilon)$ of $(\mcE, \nabla)$ over $k_\varepsilon :=k [\varepsilon]/(\varepsilon^2)$.
It corresponds to an element $v$ of $\mbH^1 (X, \mcK^\bullet [\nabla^\mr{ad}])$ via $\gamma$.
This deformation  is represented by a collection 
of data
$(\{ \partial_{\alpha \beta}\}_{(\alpha, \beta) \in I_2}, \{ \delta_\alpha \}_{\alpha \in I})$ (as displayed in  (\ref{Eqq321}))  in the \v{C}ech double complex  $\mr{Tot}^\bullet (\check{C}(\msU, \mcK^\bullet [\nabla^\mr{ad}]))$ associated to $\mcK^\bullet [\nabla^\mr{ad}]$ and an affine open covering $\msU := \{ U_\alpha \}_{\alpha \in I}$ of $X$ indexed by a finite set $I$.
Denote by   $\pi$   the projection $X \times_k k_\varepsilon \migi X$.
Then, for each $\alpha \in I$,
the flat $G$-bundle  $(\mcE_\varepsilon, \nabla_\varepsilon) |_{U_\alpha}$ (i.e., the restriction  of $(\mcE_\varepsilon, \nabla_\varepsilon)$ to $U_\alpha$) is isomorphic to $(\pi^*(\mcE)|_{U_\alpha}, \pi^*(\nabla)|_{U_\alpha} + \delta_\alpha \cdot \varepsilon)$.
If $\sigma_i \in U_{\alpha} (k)$ for some $i \in \{1, \cdots, r \}$,
then the residue of $\nabla_\varepsilon$ at $\sigma_i$ coincides with $q_{-1} + \overline{\delta}_\alpha \cdot \varepsilon \in \mfg (k_\varepsilon)$,
where $\overline{\delta}_\alpha$ denotes
 the image of $\delta_\alpha$ in $\mfg$ via  the composite surjection $\Omega \otimes_{\mcO_X} \mfg_\mcE \migisurj \sigma^*_i (\Omega) \otimes_k \sigma_i^* (\mfg_\mcE) \isom \mfg$, where the second morphism is  given by both the residue map $\sigma^*_i (\Omega)\isom k$ and the isomorphism $\sigma^*_i (\mfg_\mcE)  \isom \mfg$ induced from  the canonical trivialization $\sigma_i^*(\mcE) \isom  G$.

Now, suppose that $v$ lies in  $\mbH^1 (X, \mcK^\bullet [\nabla^\mr{ad}_\mr{par}])$.
After  possibly replacing $(\{ \partial_{\alpha \beta}\}_{(\alpha, \beta)}, \{ \delta_\alpha \}_{\alpha})$ with another,
 we may assume
that $\delta_{\alpha} \in H^0 (U_\alpha, {^c}\Omega \otimes_{\mcO_X} \mfg_\mcE)$ for every $\alpha \in I$.
Hence,  the equality  $\overline{\delta}_\alpha = 0$ holds, and the radius at each marked point  coincides with $\left(\chi (q_{-1} + 0 \cdot \varepsilon) =\right) [0]$.
This means 
 that $(\mcE_\varepsilon, \nabla_\varepsilon)$ is classified by $T_q \C_{G, \rho}$.

Conversely, suppose that  $(\mcE_\varepsilon, \nabla_\varepsilon)$ is of radii $[0]^{\times r}$.
Let us choose $i \in \{ 1, \cdots, r \}$ and choose $\alpha_i \in I$ with $\sigma_i \in U_
{\alpha_i} (k)$.
Recall from ~\cite[Theorem 2]{Kos} (or  ~\cite[Lemme 1.2.3]{Ngo}) that if $\mfg^{\mr{reg}}$ denotes the set of regular elements in $\mfg$, then 
the fiber of the projection $\chi |_{\mfg^{\mr{reg}}} : \mfg^{\mr{reg}} \migi \mfc$ over $[0]$ forms   a homogenous space 
with respect to the adjoint $G$-action.
Since both $q_{-1}$ and $q_{-1} + \overline{\delta}_{\alpha_i} \cdot \varepsilon$ belong to this fiber, there exists an element in $G (k_\varepsilon)$ of the form $e + w_i \cdot \varepsilon$, where $e$ denotes the identity element of $G$ and $w_i$ is an element of $\mfg$, satisfying 
$\mr{Ad}(e + w_i \cdot \varepsilon) (q_{-1}) \left( =q_{-1} + [w_i, q_{-1}]\cdot\varepsilon \right) = q_{-1} + \overline{\delta}_\alpha \cdot \varepsilon$.
Hence, we have $[w_i,  q_{-1}] = \overline{\delta}_\alpha$.
Let us take a section $\lambda_i \in H^0 (U_\alpha, \mfg_\mcE)$  whose image in $\left(H^0(U_\alpha, \mfg_\mcE)/H^0 (U_\alpha, {^c}\mfg_\mcE) = H^0 (U_\alpha, \sigma^*_i(\mfg_\mcE)) = \right) \mfg$ coincides with $w_i$, and 
replace  $\delta_{\alpha_i}$  (resp., $\partial_{\alpha_i \beta}$) with  the section
$\delta_{\alpha_i} - \nabla^\mr{ad}(w_i)$ (resp., $\partial_{\alpha_i \beta} - w_i |_{U_{\alpha_i \beta}}$).
Then, the resulting  collection $(\{ \partial_{\alpha \beta} \}_{(\alpha, \beta)}, \{ \delta_\alpha \}_{\alpha})$ still  represents $v$ and  furthermore lies in  $\mr{Tot}^\bullet (\check{C} (\msU, \mcK^\bullet [\nabla_\mr{par}^\mr{ad}]))$.
This means that $v$ lies in $\mbH^1 (X, \mcK^\bullet [\nabla^\mr{ad}_\mr{par}])$.

As a consequence, we conclude  that  $\gamma$ restricts to an isomorphism of $k$-vector  spaces
$T_q \C_{G, \rho}$ $\isom \mbH^1 (X, \mcK^\bullet [\nabla^\mr{ad}_\mr{par}])$, and
this completes the proofs of the former and latter assertions of this proposition.
\end{proof}
%------------------------------------------------------------------------------------
\SSP

By applying the above proposition together with the isomorphism (\ref{Eqqw9}),
we obtain an isomorphism of $k$-vector spaces
\begin{align} \label{K09}
T_q \C_{G, [0]^{\times r}} \isom T_q^\vee \C_{G, [0]^{\times r}}.
\end{align}

\LSP
%---------------------------[begin subsection]-------------
\subsection{The deformation space of a $G$-oper I} \label{SSa1}

Let $\msE^\spadesuit := (\mcE_B, \nabla)$ be a 
$G$-oper on $\msX$, and write $\mcE := \mcE_B \times^B G$.
%Then, we obtain   the following assertions.
%, some of which can be found in ~\cite[\S\,6.2.4]{Wak8}.

\SSP
%----------------------------------------------------------------------------
\bpr \label{Prop62}
Let $\Box$ denote either the absence or presence of ``$c$". 
Then, there exists 
 a short exact sequence
\begin{align}   
 0 \longmigi H^0 (X, {^\Box}\DV_G)\xrightarrow{'e_{\sharp, \Box}} \mbH^1 (X, \mcK^\bullet [{^\Box}\nabla^\mr{ad}]) \xrightarrow{{'e}_{\flat, \Box}}   H^1(X, \Omega\otimes_{\mcO_X} {^\Box}\DV_G^{\vee})  \longmigi 0.
 \label{HDeq1}
 \end{align}
Moreover,  the following diagram forms an isomorphism of short exact sequences:
\begin{align} \label{Ed3h7}
\vcenter{\xymatrix@C=26pt@R=36pt{
 0 \ar[r] & H^0 (X, {^c}\DV_G)  \ar[r]^-{'e_{\sharp, c}} \ar[d]_-{\wr}^{\mr{Serre \ duality}} & \mbH^1 (X, \mcK^\bullet [{^c}\nabla^\mr{ad}])  \ar[r]^-{{'e}_{\flat, c}} \ar[d]^-{(\ref{Eww3})} &  H^1(X, \Omega\otimes_{\mcO_X} {^c}(\DV_G^{\vee}))  \ar[r] \ar[d]_-{\wr}^{\mr{Serre \ duality}} & 0
 \\
 0 \ar[r] & H^1 (X, \Omega \otimes_{\mcO_X} \DV_G^\vee)^\vee \ar[r]_-{({'e}_{\flat})^\vee} & \mbH^1 (X, \mcK^\bullet [\nabla^\mr{ad}])^\vee \ar[r]_-{('e_\sharp)^\vee} & H^0 (X, \DV_G)^\vee \ar[r]_-{} & 0,
 }}
\end{align}
where the left-hand and right-hand vertical arrows are the isomorphism arising from Serre duality.
\epr
%----------------------------------------------------------------------------
\begin{proof}
After possibly replacing $\msE^\spadesuit$ with its $\,\qq$-normalization,
we may assume that $\msE^\spadesuit$ is $\,\qq$-normal.
The Hodge to de Rham spectral sequence 
\begin{equation}
{'E}^{a,b}_{1} :=H^b (X, \mcK^a[{^\Box}\nabla^\mr{ad}]) \Rightarrow \mbH^{a+b}(X, \mcK^\bullet[{^\Box}\nabla^\mr{ad}]),  \label{HDsp}
\end{equation}
of the complex $\mcK^\bullet [{^\Box}\nabla^\mr{ad}]$  
 yields a short exact sequence
\begin{align}   
 0 \longmigi \mr{Coker}(H^0({^\Box}\nabla^\mr{ad})) \longmigi \mbH^1(X, \mcK^\bullet [{^\Box}\nabla^\mr{ad}]) \longmigi   \mr{Ker}(H^1({^\Box}\nabla^\mr{ad}))   \longmigi 0, \label{HDeq}
 \end{align}
where  $H^j ({^\Box}\nabla^\mr{ad})$ ($j = 0, 1$) denotes the morphism $H^j (X, \mcK^0[{^\Box}\nabla^\mr{ad}]) \migi H^j(X, \mcK^1[{^\Box}\nabla^\mr{ad}])$ 
induced by ${^\Box}\nabla^\mr{ad}$.
(Namely, this is obtained from the dual of the lower horizontal sequence in (\ref{Eq3i90}).)

Here, observe that the following square diagram is commutative:
\begin{align} \label{Ed3g7}
\vcenter{\xymatrix@C=46pt@R=36pt{
   H^0(X,  {^c}\DV_G) \ar[r]^-{H^0 (\varsigma |_{{^c}\DV_G})} \ar[d]^-{\wr}_-{\mr{Serre \    duality}} & H^0(X, \Omega\otimes_{\mcO_X} {^c}\mfg_{\mcE}) \ar[r]^-{\mr{quotient}} \ar[d]^-{\wr}_-{\mr{Serre \    duality}} & \mr{Coker}(H^0 ({^c}\nabla^{\mr{ad}}))   \ar[d]^-{\wr}
   \\
 H^1 (X, \Omega \otimes_{\mcO_X} \DV_G^\vee)^\vee \ar[r]_-{H^1 (\varsigma^\veebar)^\vee}   & H^1 (X, \mfg_\mcE^\vee)^\vee \ar[r]_-{\mr{quotient}} & \mr{Ker}(H^1 (\nabla^\mr{ad}))^\vee,
 }}
\end{align}
where 
\begin{itemize}
\item
 $\varsigma^{\veebar}$ denotes  the  morphism
$\mfg_{\mcE}^\vee \migi \Omega\otimes_{\mcO_X} \DV_G^{\vee}$ 
defined as 
the tensor product of  the dual $\varsigma^\vee : \Omega^\vee \otimes_{\mcO_X} \mfg_\mcE^\vee \migi \DV_G^\vee$ of $\varsigma$ and the identity morphism of $\Omega$;
\item
the left-hand  and middle vertical arrows denote the isomorphisms
 arising from Serre duality;
\item
the right-hand vertical arrow 
denotes the isomorphism induced from (\ref{Eq3i90}) (via $\mfg_\mcE \cong \mfg_\mcE^\vee$).
\end{itemize}
It follows from  ~\cite[Lemma 6.4, (ii)]{Wak8} that the composite of the lower horizontal arrows is  an isomorphism, so the composite of the upper horizontal arrows turns out to be an isomorphism.
Also, by a similar argument together with  ~\cite[Lemma 6.4, (i)]{Wak8}, we obtain the following commutative square diagram  all of whose arrows are isomorphisms:
\begin{align} \label{Ed3gg7}
\vcenter{\xymatrix@C=46pt@R=36pt{
 \mr{Ker} (H^1 ({^c}\nabla^\mr{ad})) \ar[r]^-{\sim} \ar[d]_-{\wr} & H^1 (X, \Omega \otimes_{\mcO_X} {^c}(\DV_G^\vee)) \ar[d]^-{\wr} 
 \\
 \mr{Coker} (H^0 (\nabla^\mr{ad}))^\vee \ar[r]_-{\sim} & H^0 (X, \DV_G)^\vee.
 }}
\end{align}
Under the identifications of $k$-vector spaces given by  the isomorphisms in (\ref{Ed3g7}) and (\ref{Ed3gg7}), the short exact sequence 
(\ref{HDeq}) coincides with 
the desired  sequence (\ref{HDeq1}).
This completes the proof of the former assertion.
Moreover, the latter assertion
follows 
  from the various definitions of morphisms involved.
\end{proof}
%----------------------------------------------------------------------------
\SSP

In the resp'd portion of the following discussion, we suppose that the $G$-oper $\msE^\spadesuit$ is of radii $\rho \in \mfc (k)^{\times r}$ (where $\rho := \emptyset$ if $r = 0$).
Denote by $q$ the $k$-rational point of $\mcO p_G$ (resp., $\mcO p_{G, \rho}$) classifying  $\msE^\spadesuit$;
we use the same notation ``$q$" to denote the image of this point via the immersion  $\mr{Imm}_G : \mcO p_G \migiincl \C_G$ (resp., $\mr{Imm}_{G, \widetilde{\rho}} : \mcO p_{G, \rho} \migiincl \C_{G, \widetilde{\rho}}$).
The affine structure on $\mcO p_G$ (resp., $\mcO p_{G, \rho}$) recalled in \S\,\ref{SS036} gives an identification between 
$H^0 (X, \DV_G)$  (resp., $H^0 (X, {^c}\DV_G)$) and the space of first-order deformations of the $G$-oper $\msE^\spadesuit$ (resp., the $G$-oper $\msE^\spadesuit$ fixing  the radii).
Hence, 
if $T_q \mcO p_G$ (resp., $T_q \mcO p_{G, \rho}$) denotes the tangent space  of $\mcO p_G$ (resp., $\mcO p_{G, \rho}$) at $q$, then there exists a canonical  isomorphism of $k$-vector spaces
\begin{align} \label{Eq300}
\gamma_\sharp : T_q \mcO p_G \isom H^0 (X, \DV_G) \ \left(\text{resp.,} \ 
 \gamma_{\sharp, \rho} : T_q \mcO p_{G, \rho} \isom H^0 (X, {^c}\DV_G) 
 \right).
\end{align}
Moreover, 
$\gamma$ and 
$\gamma_\sharp$
(resp., $\gamma_{\widetilde{\rho}}$ and $\gamma_{\sharp, \rho}$)
 make the following square diagram commute:
\begin{align} \label{Ed37}
\vcenter{\xymatrix@C=33pt@R=36pt{
 T_q \mcO p_G  \ar[r] \ar[d]_-{\gamma_\sharp}^-{\wr} &  T_q \C_G  \ar[d]^-{\gamma}_-{\wr}
 \\
 H^0 (X, \DV_G)  \ar[r]_-{{'}e_\sharp} &  \mbH^1 (X, \mcK^\bullet [\nabla^{\mr{ad}}]) 
}} \ \left(\text{resp.,} \
\vcenter{\xymatrix@C=33pt@R=36pt{
 T_q \mcO p_{G, \rho}  \ar[r] \ar[d]_-{\gamma_{\sharp, \rho}}^-{\wr} &  T_q \C_{G, \widetilde{\rho}}  \ar[d]^-{\gamma_{\widetilde{\rho}}}_-{\wr}
 \\
 H^0 (X, {^c}\DV_G)  \ar[r]_-{{'}e_{\sharp, c}} &  \mbH^1 (X, \mcK^\bullet [{^c}\nabla^{\mr{ad}}])
}}
 \right),
\end{align}
where the upper horizontal arrow denotes the differential at $q$ of  $\mr{Imm}_G$ (resp., $\mr{Imm}_{G, \widetilde{\rho}}$).
Thus, we have obtained the following assertion.

\SSP
%------------------------------------------------------------------------------------------
\bt \label{Th9m}
There exist  canonical short exact sequences
\begin{align} \label{WWW45}
0 \longmigi  T_q \mcO p_G &\longmigi T_q \C_{G} \longmigi T_q^\vee \mcO p_{G, \rho} \longmigi 0, \\
 0 \longmigi  T_q \mcO p_{G, \rho} &\longmigi T_q \C_{G, \widetilde{\rho}} \longmigi T_q^\vee \mcO p_{G} \longmigi 0,  \notag
\end{align}
where 
the second arrows in the upper and lower sequences are  
 the differentials at $q$ of   $\mr{Imm}_G$ and  $\mr{Imm}_{G, \rho}$, respectively.
Moreover,
the second sequence in  (\ref{WWW45}) is compatible with the dual of the first one via the isomorphism (\ref{Er67}).
\et
%------------------------------------------------------------------------------------------
\begin{proof}
The  assertion follows from Proposition \ref{Prop62} together with 
the commutativity of the diagrams displayed  in  (\ref{Ed37}).
%the definitions of    $\gamma$, $\gamma_{\widetilde{\rho}}$, $\gamma_{\sharp}$, and $\gamma_{\sharp, \rho}$.
\end{proof}
%------------------------------------------------------------------------------------------

\LSP
%---------------------------[begin subsection]-------------
\subsection{The deformation space of a  $G$-oper II} \label{SSa5}

Next,  we prove the self-duality of  the parabolic de Rham cohomology group $\mbH^1 (X, \mcK^\bullet [\nabla^\mr{ad}_\mr{par}])$.
To this end, let us first  consider  the following lemma. 

\SSP
%------------------------------------------------------------------------------------
\ble \label{Lem55}
Suppose that 
$\msE^\spadesuit$ is $\,\qq$-normal and  of radii $[0]^{\times r}$.
For each $j \in \mbZ$, we shall set 
\begin{align} \label{Er40}
\Omega_\mr{par}^{0, \Box} (\mfg_\mcE)^j := \Omega_\mr{par}^{0, \Box} (\mfg_\mcE) \cap \mfg_\mcE^j \ \ \text{and} \ \  
\Omega_\mr{par}^{1, \Box} (\mfg_\mcE)^j := \Omega_\mr{par}^{1, \Box} (\mfg_\mcE) \cap (\Omega \otimes_{\mcO_X} \mfg_\mcE^j),
\end{align}
 where $\Box$ denotes either the absence or presence of ``$c$".
 Also, denote by 
 \begin{align} \label{Er11}
 {^\Box}\nabla_\mr{par}^{\mr{ad}(j)} : \Omega_\mr{par}^{0, \Box} (\mfg_\mcE)^j \migi  \Omega_\mr{par}^{0, \Box} (\mfg_\mcE)^{j-1} 
 \end{align}
 the morphism obtained by restricting $\nabla_\mr{par}^\mr{ad}$.
 Then, the following assertions hold:
\begin{itemize}
\item[(i)]
For each integer $j$, the composite
\begin{align} \label{Er10}
H^0 (X, {^c}\DV_G^j) \migi H^0 (X, \Omega^{1, c}_\mr{par}(\mfg_\mcE)^j) \xrightarrow{\mr{quotient}} 
\mr{Coker}(H^0 ({^c}\nabla_\mr{par}^{\mr{ad} (j+1)}))
\end{align}
is an isomorphism,
where the first arrow arises  from the natural inclusion  ${^c}\DV_G^j \migiincl \Omega_\mr{par}^{1, c} (\mfg_\mcE)^j$.
In particular,  we have
\begin{align} \label{Er2}
H^0 (X, {^c}\DV_G) \cong \mr{Coker}(H^0 ({^c}\nabla_\mr{par}^\mr{ad})).
\end{align}
\item[(ii)]
For each integer $j$, the composite
\begin{align} \label{Er17}
\mr{Ker}(H^1 (\nabla^{\mr{ad} (j)}_\mr{par})) \xrightarrow{\mr{inclusion}}  H^1 (X, \Omega_\mr{par}^{0}(\mfg_{\mcE})^j)
\migi H^1 (X, \Omega \otimes_{\mcO_X} \DV_G^{\vee, j})
\end{align}
is an isomorphism, where $ \DV_G^{\vee, j} := (\DV_G/\DV_G^{-j +1})^\vee \left(\subseteq \DV_G^\vee \right)$ and  the second  arrow arises  from the natural surjection 
$\Omega_\mr{par}^{0}(\mfg_{\mcE})^j \migisurj \Omega \otimes_{\mcO_X} \DV_G^{\vee, j}$.
In particular, we have
\begin{align} \label{Er18}
\mr{Ker} (H^1 (\nabla^\mr{ad}_\mr{par})) \cong H^1 (X, \Omega \otimes_{\mcO_X} \DV_G^\vee).
\end{align}
\end{itemize}
\ele
%------------------------------------------------------------------------------------
\begin{proof}
We only consider assertion (ii) since assertion (i) follows from a similar argument.
For each integer $j$,
we shall write $\nabla_\mr{par}^{\mr{ad}(j/j+1)}$ for  the $\mcO_X$-linear  morphism
$\Omega_\mr{par}^0 (\mfg_\mcE^j)/\Omega_\mr{par}^0 (\mfg_\mcE^{j+1}) \migi \Omega_\mr{par}^1 (\mfg_\mcE)^{j-1}/\Omega_\mr{par}^1 (\mfg_\mcE)^{j}$ induced from $\nabla_\mr{par}^\mr{ad}$.

First, let us consider the case  of $j \geq 0$.
Since the residue of $\nabla$ at every marked point coincides with $q_{-1}$ (by the assumption that $\msE^\spadesuit$ is $\,\qq$-normal  and of radii $[0]^{\times r}$),
one may verify that the resp'd portion of ~\cite[(768)]{Wak8} restricts to a {\it split} exact sequence
\begin{align} \label{Er31}
0 \longmigi 
\Omega_\mr{par}^0 (\mfg_\mcE)^j/\Omega_\mr{par}^0 (\mfg_\mcE)^{j+1}
\xrightarrow{\nabla_\mr{par}^{\mr{ad}(j/j+1)}}
\Omega_\mr{par}^1 (\mfg_\mcE)^{j-1}/\Omega_\mr{par}^1 (\mfg_\mcE)^{j}
\longmigi
\DV_{G, j-1}
\longmigi 0.
\end{align}
In particular, the equality $\mr{Ker}(H^1 (\nabla_\mr{par}^{\mr{ad}(j/j+1)})) =0$ holds (by which we obtain $\mr{Ker}(H^1 (\nabla_\mr{par}^{\mr{ad}(j)})) = 0$).
Hence, the assertion for  $j \geq 0$ follows from the fact that
$\DV_G^{\vee, j} \left(= (\DV_G^{}/\DV_G^{-j +1})^\vee \right) = 0$.

Next, we shall  prove the case of $j \leq -1$.
Consider the following  morphism of sequences:
\begin{align} \label{MorSeq}
\vcenter{\xymatrix@C=21pt@R=36pt{
0 \ar[r]& 
 H^1 (\Omega_\mr{par}^0(\mfg_{\mcE})^{j+1})
\ar[r] \ar[d]^-{ H^1(\nabla_\mr{par}^{\mr{ad}(j+1)})}
&H^1 (\Omega_\mr{par}^0 (\mfg_\mcE)^{j})
 \ar[r] \ar[d]^-{H^1(\nabla_\mr{par}^{\mr{ad}(j)})}
& 
H^1 (\Omega_\mr{par}^0 (\mfg_\mcE^j)/\Omega_\mr{par}^0 (\mfg_\mcE^{j+1}))
\ar[r] \ar[d]^-{H^1(\nabla_\mr{par}^{\mr{ad}(j/j+1)})} & 0
\\
0 \ar[r]&
H^1 (\Omega_\mr{par}^1 (\mfg_\mcE)^j)
   \ar[r] &
H^1 (\Omega_\mr{par}^1 (\mfg_\mcE)^{j-1})
\ar[r]
& H^1 (\Omega_\mr{par}^1 (\mfg_\mcE)^{j-1}/\Omega_\mr{par}^1 (\mfg_\mcE)^{j})
 \ar[r]& 0.
}}
\end{align}
Since the residue of $\nabla$ at every marked point coincides with $q_{-1}$,
the line bundles 
$\Omega_\mr{par}^0 (\mfg_\mcE^j)/\Omega_\mr{par}^0 (\mfg_\mcE^{j+1})$
 and 
 $\Omega_\mr{par}^1 (\mfg_\mcE)^{j-1}/\Omega_\mr{par}^1 (\mfg_\mcE)^{j}$
 coincide with $\mfg_\mcE^j/\mfg_\mcE^{j+1}$ and $\Omega \otimes_{\mcO_X} (\mfg_\mcE^{j-1}/\mfg_\mcE^{j})$, respectively.
In particular, these 
are isomorphic to  direct sums of finite  copies of $\Omega^{\otimes j}$ (cf. ~\cite[\S\,2.1.4]{Wak8}), 
so we have 
\begin{equation} 
\label{f=0}
H^0 (X, \Omega_\mr{par}^0 (\mfg_\mcE^j)/\Omega_\mr{par}^0 (\mfg_\mcE^{j+1}))
= H^0 (X, \Omega_\mr{par}^1 (\mfg_\mcE)^{j-1}/\Omega_\mr{par}^1 (\mfg_\mcE)^{j}) = 0.
\end{equation}
By 
the above equalities
together with  the fact  that $\mr{dim}(X) =1$ (which implies  $H^2 (-) =0$), 
both the upper and lower horizontal sequences in  (\ref{MorSeq}) turn out to be  exact.
Also,  
the right-hand vertical  arrow in (\ref{MorSeq}) 
 is surjective because of   the long exact sequence arising from the non-resp'd portion of 
 ~\cite[(768)]{Wak8}.
By descending induction on $j$, 
one  can verify the surjectivity  of  both the left-hand  and middle  vertical  arrows  in (\ref{MorSeq}).
Thus,   
the snake lemma  applied to (\ref{MorSeq}) shows that
 the natural   sequence
\begin{equation} 0 \migi   \mr{Ker}(H^1 (\nabla_\mr{par}^{\mr{ad}(j+1)}))\migi  \mr{Ker}(H^1 (\nabla_\mr{par}^{\mr{ad}(j)})) \migi  \mr{Ker}(H^1 (\nabla_\mr{par}^{\mr{ad}(j/j+1)})) \migi 0  \label{Ker} \hspace{-3mm}
\end{equation}
is  exact.
Let us consider  the  morphism of short exact sequences
\begin{align} \label{MorSeq2}
\vcenter{\xymatrix@C=21pt@R=36pt{
0 \ar[r]& \mr{Ker}(H^1(\nabla_\mr{par}^{\mr{ad}(j+1)})) \ar[d] \ar[r] 
&\mr{Ker}(H^1(\nabla_\mr{par}^{\mr{ad}(j)}))  \ar[r] \ar[d]^-{(\ref{Er17})} &  \mr{Ker}(H^1 (\nabla_\mr{par}^{\mr{ad}(j/j+1)}))
\ar[d] \ar[r] & 0
\\
0 \ar[r]& H^1(\Omega\otimes_{\mcO_X} \DV_G^{\vee, j+1})  \ar[r]& H^1(\Omega\otimes_{\mcO_X} \DV_G^{\vee, j})   \ar[r]
& H^1(\Omega\otimes_{\mcO_X} (\DV_{G, -j})^\vee) \ar[r]& 0,
}}\hspace{-4mm}
\end{align} 
where  the left-hand vertical arrow is  (\ref{Er17}) with $j$ replaced by $j+1$ and the exactness of the lower horizontal sequence follows from  the natural  decomposition $\DV_G^{\vee, j} = \DV_G^{\vee, j+1} \oplus (\DV_{G, -j})^\vee$.
 The right-hand vertical  arrow in 
this  diagram is an isomorphism because it coincides with the inverse of the isomorphism $\mbR^1 f_*(\Omega \otimes_{\mcO_X} (\DV_{G, -j})^\vee) \isom \mr{Ker}(\mbR^1 f_*(\nabla^{\mr{ad}(j/j+1)}))$ arising from the splitting  of 
the non-resp'd portion of 
~\cite[(768)]{Wak8}.
Hence, by descending induction on $j$, we see that  the middle vertical  arrow (\ref{Er17})   is   an isomorphism.
This completes the proof of the assertion.
\end{proof}
%------------------------------------------------------------------------------------
\SSP

By using the above lemma, we obtain the following proposition.

\SSP
%------------------------------------------------------------------------------------
\bpr \label{Er35}
Suppose that 
the 
$G$-oper $\msE^\spadesuit$ is $\,\qq$-normal and  of radii $[0]^{\times r}$ (where $[0]^{\times r} := \emptyset$ if $r = 0$).
Then, the following assertions hold:
\begin{itemize}
\item[(i)]
$\mbH^0 (X, \mcK^\bullet [\nabla_{\mr{par}}^\mr{ad}]) = \mbH^2 (X, \mcK^\bullet [\nabla_{\mr{par}}^\mr{ad}]) = 0$.
\item[(ii)]
There exists a short exact sequence
\begin{align} \label{Eqr}
0 \longmigi H^0 (X, {^c}\DV_G) \xrightarrow{{'}e_{\sharp, \mr{par}}} \mbH^1 (X, \mcK^\bullet [\nabla^\mr{ad}_\mr{par}]) \xrightarrow{{'}e_{\flat, \mr{par}}} H^1 (X, \Omega \otimes_{\mcO_X} \DV_G^\vee) \longmigi 0,
\end{align}
which fits into the following isomorphism of short exact sequences:
\begin{align} \label{Ed3hh7}
\vcenter{\xymatrix@C=31pt@R=36pt{
 0 \ar[r] & H^0 (X, {^c}\DV_G)  \ar[r]^-{'e_{\sharp, \mr{par}}} \ar[d]^-{\mr{Serre \ duality}}_-{\wr} & \mbH^1 (X, \mcK^\bullet [\nabla^\mr{ad}_\mr{par}])  \ar[r]^-{{'e}_{\flat, \mr{par}}} \ar[d]^-{(\ref{Eqqw9})} &  H^1(X, \Omega\otimes_{\mcO_X} \DV_G^{\vee})  \ar[r] \ar[d]^-{\mr{Serre \ duality}}_-{\wr} & 0
 \\
 0 \ar[r] & H^1 (X, \Omega \otimes_{\mcO_X} \DV_G^\vee)^\vee \ar[r]_-{({'e}_{\flat, \mr{par}})^\vee} & \mbH^1 (X, \mcK^\bullet [\nabla^\mr{ad}_\mr{par}])^\vee \ar[r]_-{('e_{\sharp, \mr{par}})^\vee} & H^0 (X, {^c}\DV_G)^\vee \ar[r]_-{} & 0.
 }}
\end{align}
In particular, $\mbH^1 (X, \mcK^\bullet [\nabla_{\mr{par}}^\mr{ad}])$ is a $k$-vector space of dimension  $(2g-2 +r) \cdot \mr{dim}(\mfg) - r \cdot \mr{rk} (\mfg)$.
\end{itemize}
\epr
%------------------------------------------------------------------------------------
\begin{proof}
First, we shall prove assertion (i).
To this end, 
it suffices to verify the equality $\mbH^2 (X, \mcK^\bullet [\nabla_{\mr{par}}^\mr{ad}]) = 0$ because the proof of the remaining one 
is similar.
Since $H^1 (X, \DV_{G, j-1}) = 0$ for every $j$, it follows from the exactness of  (\ref{Er31})
 that the morphism 
\begin{align}
H^1 (\nabla_{\mr{par}}^{\mr{ad} (j/j+1)}) :  H^1 (X, \Omega_\mr{par}^0 (\mfg_\mcE)^{j}/\Omega_\mr{par}^0 (\mfg_\mcE)^{j+1}) \migi H^1 (X, \Omega_\mr{par}^1 (\mfg_\mcE)^{j-1}/\Omega_\mr{par}^1 (\mfg_\mcE)^{j})
\end{align}
induced from  $\nabla_{\mr{par}}^{\mr{ad} (j/j+1)}$ is surjective when $j \geq 0$.
Next,  consider the following morphism of short exact sequences defined for each $j \in \mbZ$:
\begin{align} \label{f}
\vcenter{\xymatrix@C=31pt@R=36pt{
 0 \ar[r] & \Omega_\mr{par}^0 (\mfg_\mcE)^{j+1} \ar[r] \ar[d]^-{\nabla_\mr{par}^{\mr{ad}(j+1)}} &  \Omega_\mr{par}^0 (\mfg_\mcE)^{j} \ar[r] \ar[d]^-{\nabla_\mr{par}^{\mr{ad}(j)}} & \Omega_\mr{par}^0 (\mfg_\mcE)^{j}/\Omega_\mr{par}^0 (\mfg_\mcE)^{j+1} \ar[r] \ar[d]^-{\nabla_{\mr{par}}^{\mr{ad} (j/j+1)}} & 0
  \\
 0 \ar[r] & \Omega_\mr{par}^1 (\mfg_\mcE)^{j} \ar[r] & \Omega_\mr{par}^1 (\mfg_\mcE)^{j-1} \ar[r] & \Omega_\mr{par}^1 (\mfg_\mcE)^{j-1}/\Omega_\mr{par}^1 (\mfg_\mcE)^{j} \ar[r] & 0.
 }}
\end{align}
By descending induction on $j$, 
the surjectivity of $H^1 (\nabla_{\mr{par}}^{\mr{ad} (j/j+1)})$ implies 
 that of $H^1 (\nabla_{\mr{par}}^{\mr{ad}(j)}) : H^1 (X, \Omega_\mr{par}^0 (\mfg_\mcE)^j) \migi H^1 (X, \Omega_\mr{par}^0 (\mfg_\mcE)^{j-1})$.
Then, the  required equality  $\mbH^2 (X, \mcK^\bullet [\nabla_{\mr{par}}^\mr{ad}]) = 0$ follows from
the surjectivity  of  $H^1 (\nabla_{\mr{par}}^{\mr{ad}(j)})$ for  $j= -\mr{rk}(\mfg)$ together with the Hodge to de Rham spectral sequence $E_1^{a, b} = H^b(X, \mcK^a [\nabla_{\mr{par}}^{\mr{ad}}]) \Rightarrow \mbH^{a+b} (X, \mcK^\bullet [\nabla_\mr{par}^\mr{ad}])$.

Assertion (ii)  follows from the commutativity of the diagram (\ref{Eqq301}), the definition of  (\ref{Eqqw9}), Lemma \ref{Lem55}, (i) and (ii), and ~\cite[Proposition 2.23, (ii)]{Wak8}. 
\end{proof}
%------------------------------------------------------------------------------------
\SSP

Let us keep the assumption in the above proposition.
As mentioned at the end of \S\,\ref{SS036},
there exists a canonical affine  structure on 
%Recall that the affine structure on $\mcO p_G$ discussed in restricts to an affine structure on
 $\mcO p_{G, [0]^{\times r}}$ modeled on $H^0 (X, {^c}\DV_G)$.
This affine structure  induces an isomorphism
\begin{align}
\gamma_{\sharp, [0]^{\times r}} : T_q  \mcO p_{G, [0]^{\times r}} \isom H^0 (X, {^c}\DV_G).
\end{align}
It follows from the  various constructions involved that 
  the following square diagram is commutative:
\begin{align} \label{Eq340}
\vcenter{\xymatrix@C=46pt@R=36pt{
 T_q \mcO p_{G, [0]^{\times r}}  \ar[r]
  \ar[d]_-{\gamma_{\sharp, [0]^{\times r}}}^-{\wr} &  T_q \C_{G, [0]^{\times r}}  \ar[d]^-{\gamma_{[0]^{\times r}}}_-{\wr}
 \\
 H^0 (X, {^c}\mcV_G)  \ar[r]_-{{'}e_{\sharp, \mr{par}}} &  \mbH^1 (X, \mcK^\bullet [\nabla_\mr{par}^{\mr{ad}}]), 
}}
\end{align}
where the upper horizontal arrow denotes the differential at $q$ of $\mr{Imm}_{G, [0]^{\times r}}$.

\SSP
%------------------------------------------------------------------------------------
\bt \label{Prc1}
\begin{itemize}
\item[(i)]
There exists a canonical short exact sequence
\begin{align}
0 \longmigi  T_q \mcO p_{G, [0]^{\times r}} \longmigi T_q \C_{G, [0]^{\times r}} \longmigi T_q^\vee \mcO p_{G, [0]^{\times r}} \longmigi 0,
\end{align}
where 
the second arrow is 
the differential of  $\mr{Imm}_{G, [0]^{\times r}}$ at $q$.
Moreover, this sequence 
 is, in an evident sense,  invariant  under taking the dual via 
 (\ref{K09}).
\item[(ii)]
Let us consider $T_q \mcO p_{G, [0]^{\times r}}$ as a $k$-vector subspace of $T_q \C_{G, [0]^{\times r}}$ via the second arrow in  the lower sequence of  (\ref{WWW45}).
Then, $T_q \mcO p_{G, [0]^{\times r}}$ is Lagrangian with respect to the nondegenerate  pairing $T_q \C_{G, [0]^{\times r}} \times T_q \C_{G, [0]^{\times r}} \migi k$ induced by (\ref{K09}).
 \end{itemize}
\et
%------------------------------------------------------------------------------------
\begin{proof}
Assertion (i)  follows from  Proposition \ref{Er35}, (ii), together with the isomorphisms  $\gamma_{[0]^{\times r}}$ and  $\gamma_{\sharp, [0]^{\times r}}$.

Next, by the definition of ${'}e_{\sharp, \mr{par}}$, we see that  $T_q \mcO p_{G, [0]^{\times r}}$ is isotropic.
Hence,  assertion (ii)  follows from the equality $\mr{dim}(T_q \mcO p_{G, [0]^{\times r}}) = \frac{1}{2} \cdot \mr{dim}(T_q \C_{G, [0]^{\times r}})$ resulting from assertion (i). 
\end{proof}
%------------------------------------------------------------------------------------

%%%%%%%%%%%%%%%%%%%%%%%%%%%%%%%%%%%%%%%%%%%%%
%%%%%%%%%%%%%---[begin section]---%%%%%%%%%%%%%%%%%%%%%%%
\vspace{10mm}
\section{Cohomology of the symmetric products of opers}\label{S23g0}\SSP

In this section,
we examine  the (parabolic) de Rham cohomology of symmetric products of a $\mr{PGL}_2$-oper (or an $\mr{SL}_2$-oper).
We will show that, if a given  $G$-oper comes from a $\mr{PGL}_2$-oper via change of structure group, 
then the (parabolic) de Rham cohomology group  of the induced adjoint bundle  decomposes into the direct sum of such cohomology groups (cf. (\ref{Eq83f4})).

\LSP
%---------------------------[begin subsection]-------------
\subsection{The de Rham cohomology of an $\mr{SL}_n$-oper} \label{SS08700}

We shall fix an integer  $n$ with $1 < n < p$, and 
 let $\msF^\heartsuit := (\mcF, \nabla, \{ \mcF^j \}_{i=0}^n)$ be   an $\mr{SL}_n$-oper on $\msX$ (cf. Remark \ref{Rem667}).
 Write $\msF := (\mcF, \nabla)$ and  ${^c}\msF := ({^c}\mcF, {^c}\nabla)$, where ${^c}\nabla$ denotes the log connection on ${^c}\mcF$ obtained by restricting $\nabla$.
 Let $\Box$ denote either the absence or presence of ``$c$".
 Then, 
 the natural morphisms of complexes $\Omega\otimes_{\mcO_X} {^\Box}\mcF^{n-1}[-1] \migi \mcK^\bullet [{^\Box}\nabla]$ and $\mcK^\bullet [{^\Box}\nabla] \migi {^\Box}\mcF/{^\Box}\mcF^1 [0]$ together  induce a  sequence of $k$-vector spaces
\begin{align} \label{Eq116}
0 \longmigi  H^0 (X, \Omega\otimes_{\mcO_X} {^\Box}\mcF^{n-1}) \longmigi H_{\mr{dR}}^1 (X, {^\Box}\msF)
\longmigi H^1 (X, {^\Box}\mcF/{^\Box}\mcF^1) \longmigi 0.
\end{align}
Also, the  natural morphisms  $\Omega \otimes_{\mcO_X} {^c}\mcF^{n-1}[-1] \migi \mcK^\bullet [\nabla_{\mr{per}}]$ and $\mcK^\bullet [\nabla_{\mr{per}}] \migi \mcF/\mcF^1 [0]$  induce a  sequence of $k$-vector spaces
\begin{align} \label{Eq648}
0 \longmigi  H^0 (X, \Omega\otimes_{\mcO_X} {^c}\mcF^{n-1}) \longmigi H_{\mr{dR}, \mr{per}}^1 (X, \msF)
\longmigi H^1 (X, \mcF/\mcF^1) \longmigi 0.
\end{align}

Then,  the following assertion holds.

\SSP
%------------------------------------------------------------------------------
\bpr \label{Prop15}
\begin{itemize}
\item[(i)]
The sequence (\ref{Eq116}) is exact.
\item[(ii)]
Suppose that 
$\msF^\heartsuit$ is of radii $[0]^{\times r} \in \mfc (k)^{\times r}$ (where $[0]^{\times r} := \emptyset$ if $r = 0$).
Then, the sequence (\ref{Eq648}) is exact.
\end{itemize}
\epr
%------------------------------------------------------------------------------
\begin{proof}
We only consider assertion (ii) because the proof of assertion (i) is relatively 
simpler than that of assertion (ii).
For each $j = 0, \cdots, n$, we shall set 
$\Omega_{\mr{par}}^0 (\mcF)^j := \mcF^j$ and 
$\Omega_{\mr{par}}^1  (\mcF)^j := 
\Omega \otimes_{\mcO_X} {^c}\mcF^j + \nabla (\mcF^{j+1})$, where $\mcF^{n+1} := 0$.
Also,  denote by $\nabla^j_{\mr{par}}$ ($j = 1, \cdots, n$) the morphism $\Omega^0_{\mr{par}} (\mcF)^j \migi \Omega^1_{\mr{par}} (\mcF)^{j-1}$ obtained by restricting $\nabla$.
By descending induction on $j \geq 1$, we shall prove the claim that 
the morphism 
\begin{align} \label{Eq400}
H^l (X, \Omega \otimes_{\mcO_X}  {^c}\mcF^{n-1}) \migi  \mbH^{l+1} (X, \mcK^\bullet [\nabla^j_{\mr{par}}])
\end{align}
induced by the  inclusion  $\Omega \otimes_{\mcO_X} {^c}\mcF^{n-1}[-1] \migiincl \mcK^\bullet [\nabla^j_{\mr{par}}]$ is an isomorphism 
 for every $l$.
The base step, i.e., the case of $j= n$ is clear because $\nabla^n_\mr{par}$ coincides with the zero map $0 \migi \Omega \otimes_{\mcO_X} {^c}\mcF^{n-1}$.
To prove the induction step, we assume that the case of $j = j_0$ ($1 < j_0 \leq n$) has been proved.
 Note that the Kodaira-Spencer map $\mr{KS}_{\msF^\heartsuit}^{j_0 -1}$  (cf. (\ref{Eqq212})) decomposes into the composite of morphisms between line bundles
 \begin{align}
\left(\mcF^{j_0-1}/\mcF^{j_0}  = \right)\Omega_{\mr{par}}^0 (\mcF)^{j_0-1} /\Omega_{\mr{par}}^0 (\mcF)^{j_0}
&\xrightarrow{\mr{KS}_{\msF^\heartsuit, \mr{par}}^{j_0 -1}} \Omega_{\mr{par}}^1 (\mcF)^{j_0-2} /\Omega_{\mr{par}}^1 (\mcF)^{j_0-1} \\
& \xrightarrow{\mr{inclusion}}
 \Omega  \otimes_{\mcO_X} (\mcF^{j_0-2}/\mcF^{j_0-1}), \notag
  \end{align}
  where the first arrow, i.e., $\mr{KS}_{\msF^\heartsuit, \mr{par}}^{j_0 -1}$,  is the morphism induced from $\nabla^j_{\mr{par}}$.
Since  $\mr{KS}_{\msF^\heartsuit}^{j_0 -1}$ is an isomorphism by assumption, 
the morphism $\mr{KS}_{\msF^\heartsuit, \mr{par}}^{j_0 -1}$ is verified to be  an isomorphism.
The morphism  $\mbH^{l+1} (X, \mcK^\bullet [\nabla^{j_0 }_{\mr{par}}]) \migi  \mbH^{l+1} (X, \mcK^\bullet [\nabla^{j_0 -1}_{\mr{par}}])$  induced by the inclusion $\mcK^\bullet [\nabla^{j_0 }_{\mr{par}}] \migiincl \mcK^\bullet [\nabla^{j_0 -1}_{\mr{par}}]$ 
is an isomorphism because the following diagram forms a morphism of short exact sequences:
\begin{align} \label{Eq340}
\vcenter{\xymatrix@C=36pt@R=36pt{
0\ar[r] & \Omega_{\mr{par}}^0 (\mcF)^{j_0}\ar[r]^-{\mr{inclusion}} \ar[d]^-{\nabla^{j_0}_{\mr{par}}} &\Omega_{\mr{par}}^0 (\mcF)^{j_0-1} \ar[r]^-{\mr{quotient}} \ar[d]^-{\nabla^{j_0-1}_{\mr{par}}} & \Omega_{\mr{par}}^0 (\mcF)^{j_0-1}/\Omega_{\mr{par}}^0 (\mcF)^{j_0} \ar[r] \ar[d]_{\wr}^-{\mr{KS}_{\msF^\heartsuit, \mr{par}}^{j_0-1}} & 0
\\
0\ar[r] & \Omega_{\mr{par}}^1 (\mcF)^{j_0 -1}\ar[r]_-{\mr{inclusion}} & \Omega_{\mr{par}}^1 (\mcF)^{j_0-2} \ar[r]_-{\mr{quotient}} & \Omega_{\mr{par}}^1 (\mcF)^{j_0-2}/\Omega_{\mr{par}}^1 (\mcF)^{j_0-1}\ar[r] & 0.
}}
\end{align}
Hence, the induction hypothesis implies that
the morphism $H^l (X, \Omega \otimes_{\mcO_X} {^c}\mcF^{n-1}) \migi \mbH^{l+1} (X, \mcK^\bullet [\nabla^{j_0-1}_{\mr{par}}])$ induced by 
 the composite of natural inclusions 
\begin{align}
\Omega \otimes_{\mcO_X} {^c}\mcF^{n-1}[-1] \migiincl \mcK^\bullet [\nabla_{\mr{par}}^{j_0}] \migiincl \mcK^\bullet [\nabla_{\mr{par}}^{j_0-1}]
\end{align}
 is  an isomorphism.
This completes the proof of the claim.
In particular,  for every $l$,  we obtain an isomorphism
\begin{align} \label{Eq5423}
H^l (X, \Omega \otimes_{\mcO_X} {^c}\mcF^{n-1}) \isom 
\mbH^{l+1} (X, \mcK^\bullet [\nabla^1_{\mr{par}}]).
\end{align}

We go back to the proof of assertion (ii).
According to ~\cite[Proposition 4.55]{Wak8},
we may assume that the $\mr{PGL}_n$-oper induced by $\msF^\heartsuit$ via projectivization is normal in the sense of ~\cite[Definition 4.53]{Wak8}.
Let us fix $i \in \{1, \cdots, r \}$, and moreover, choose a local function $t$ defining 
$\sigma_i \in X (k)$.
Then,  around the marked point $\sigma_i$, the log connection $\nabla$ may be described as 
\begin{align} \label{Eq998}
\nabla = d +  \frac{dt}{t} \otimes  \begin{pmatrix} 0 & 0 & 0 & \cdots & 0 & a_n
\\
1 & 0 & 0 & \cdots & 0 & a_{n-1}
\\
0 & 1 & 0 & \cdots & 0 & a_{n-2}
\\
0 & 0 & 1 & \cdots & 0 & a_{n-3}
\\
\vdots & \vdots & \vdots & \ddots & \vdots & \vdots
\\
0 & 0 & 0 & \cdots & 1 & a_1
\end{pmatrix}
\end{align}
for some local functions $a_1, \cdots, a_{n}$ (cf. ~\cite[Remark 4.30]{Wak8}).
Then, the assumption in (ii) together with ~\cite[Theorem 4.49]{Wak8} implies that
$a_j \equiv 0$ mod $(t)$ for every $j = 1, \cdots, n$.
This implies 
 $\Omega^1_{\mr{par}}(\mcF)^{0} =
 \Omega^1_{\mr{par}}(\mcF)$.
Hence,
the diagram
 \begin{align} \label{Eq340}
\vcenter{\xymatrix@C=36pt@R=36pt{
0\ar[r] & \Omega_{\mr{par}}^0 (\mcF)^{1}\ar[r]^-{\mr{inclusion}} \ar[d]^-{\nabla^{1}_{\mr{par}}} & \Omega_{\mr{par}}^0 (\mcF) \ar[r]^-{\mr{quotient}} \ar[d]^-{\nabla_{\mr{par}}} & \mcF^{}/\mcF^{1} \left(=  \Omega_{\mr{par}}^0 (\mcF)/ \Omega_{\mr{par}}^0 (\mcF)^{1} \right) \ar[r] \ar[d] & 0
\\
0\ar[r] &\Omega^1_{\mr{par}}(\mcF)^{0}\ar[r]_-{\sim} &\Omega_{\mr{par}}^1 (\mcF)\ar[r]_-{} & 0\ar[r] & 0
}}
\end{align}
forms a short exact sequence of morphisms, which 
induces the long exact sequence of $k$-vector spaces
\begin{align} \label{Eq399}
H^0 (X, \mcF/\mcF^1)  &\migi  \mbH^1 (X, \mcK^\bullet [\nabla^1_{\mr{par}}]) \migi \mbH^1 (X, \mcK^\bullet [\nabla_{\mr{par}}])   \\
&
\migi H^1 (X, \mcF/\mcF^1) \migi  \mbH^2 (X, \mcK^\bullet [\nabla^1_{\mr{par}}]).  \notag
\end{align}
Since $\mr{deg}(\mcF/\mcF^1) < 0$,  the equality  $H^0 (X, \mcF/\mcF^1) = 0$ holds.
On the other hand, by the claim proved above,
the $k$-vector space $\mbH^2 (X, \mcK^\bullet [\nabla^1_{\mr{par}}]) \left(\cong H^1 (X, \Omega \otimes_{\mcO_X} {^c}\mcF^{n-1}) \cong H^0 (X, (\mcF^{n-1})^\vee)^\vee \right)$ is verified to be zero because of the fact that $\mr{deg}(\mcF^{n-1}) >0$.
Hence,   (\ref{Eq399}) together with 
(\ref{Eq5423}) for $l=0$  implies the exactness of  (\ref{Eq648}), thus completing the proof of this proposition.
\end{proof}
%------------------------------------------------------------------------------
%\SSP

\LSP
%---------------------------[begin subsection]-------------
\subsection{Duality for  symmetric products of opers} \label{SSa2}

Next, let us describe the short exact sequences (and their dualities)  discussed in the previous subsection 
for the $\mr{SL}_n$-oper arising from an $\mr{SL}_2$-oper/a $\mr{PGL}_2$-oper via change of structure group.
If the underlying space is a compact Riemann surface, then 
the  assertion corresponding to the following theorem can be found in ~\cite[Theorem 6]{Gun} (for $n=3$) and  ~\cite[\S\,4.4.4]{Wen}.
The case of
%where  $\msF^\heartsuit_\odot$ is 
 the uniformizing  $\mr{SL}_2$-opers on modular curves
 was
 already proved in ~\cite[Theorem 2.7]{Sch}.
 (Note that, by  the discussion  at the beginning of \S\,2.4 in {\it loc.\,cit.},
 the radius of the uniformizing  $\mr{SL}_2$-oper at every marked point coincides with $[0]$;    hence, the following Theorem \ref{Th67}, (ii), can be regarded as a generalization of that result.)
  Also, see ~\cite[Theorem 3.1, Remark 3.2]{Og2} for a version of elliptic $F$-$T$-crystals.

\SSP
%----------------------------------------------------------------------------------------------------
\bt \label{Th67}
Let  $\msF^\heartsuit_\odot := (\mcF_\odot, \nabla_\odot, \{ \mcF^j_\odot \}_{j=0}^2)$  
be an $\mr{SL}_2$-oper on $\msX$.
 Write $\varTheta$ for the theta characteristic associated to $\msF^\heartsuit$ (cf. \S\,\ref{SS03g61}).
 Also, denote by  $\mr{Sym}^{n-1}\msF_\odot$ 
 the $(n-1)$-st symmetric product
of the underlying flat $\mr{SL}_2$-bundle $\msF_\odot$ of $\msF^\heartsuit_\odot$.
Then, the following assertions hold:
\begin{itemize}
\item[(i)]
Let $\Box$ denote either the absence or presence of ``$c$".
Then, there exists 
a canonical short exact sequence of $k$-vector spaces
\begin{align} \label{Eq8177}
0 \migi H^0 (X, {^\Box}(\varTheta^{\otimes (n+1)}))
\migi
%\xrightarrow{{^\Box}h_{\sharp}^n [\msF_\odot^\heartsuit]}
H^1_{\mr{dR}} (X, {^\Box}(\mr{Sym}^{n-1}\msF_\odot))
%\xrightarrow{{^\Box}h_{\flat}^n [\msF_\odot^\heartsuit]}
 \migi
 H^1 (X, {^\Box}(\varTheta^{\otimes (-n+1)}))
\migi 0.
\end{align}
Moreover, there exists a canonical isomorphism
\begin{align} \label{KKK98}
H^1_{\mr{dR}} (X, {^c}(\mr{Sym}^{n-1}\msF_\odot)) \isom 
H^1_{\mr{dR}} (X, \mr{Sym}^{n-1}\msF_\odot)^\vee
\end{align}
via which the short exact sequences  (\ref{Eq8177}) for $\Box$'s in  both cases are compatible with each other.
\item[(ii)]
Suppose further that
$\msF^\heartsuit_\odot$ is of radii $[0]^{\times r} \in \mfc (k)^{\times r}$ (where $[0]^{\times r} := \emptyset$ if $r = 0$).
Then, there exists a canonical short exact sequence
 \begin{align} \label{Eq810}
0 \migi H^0 (X, \varTheta^{\otimes (n+1)}(-D))
%\xrightarrow{h_{\sharp, \mr{par}}^n [\msF_\odot^\heartsuit]}
\migi
H^1_{\mr{dR}, \mr{par}} (X, \mr{Sym}^{n-1}\msF_\odot)
%\xrightarrow{h_{\flat, \mr{par}}^n [\msF_\odot^\heartsuit]}
\migi
H^1 (X, \varTheta^{\otimes (-n+1)})
\migi 0. 
\end{align}
Moreover, 
 there exists a canonical isomorphism
\begin{align} \label{Ett38}
H^1_{\mr{dR}, \mr{par}} (X, \mr{Sym}^{n-1}\msF_\odot) \isom H^1_{\mr{dR}, \mr{par}} (X, \mr{Sym}^{n-1}\msF_\odot)^\vee
\end{align}
via which the short exact sequence (\ref{Eq810}) is compatible with its dual.
\end{itemize}
\et
%----------------------------------------------------------------------------------------------------
\begin{proof}
We only prove assertion (ii) because the proof of assertion (i) is similar.
For simplicity, we shall write $(\mcF, \nabla) := (\mr{Sym}^{n-1}\mcF_\odot,  \mr{Sym}^{n-1}\nabla_\odot)$.

First, we shall  consider the former assertion.
Under the identifications 
$\mcF^{n-1} = \varTheta^{\otimes (n-1)}$ and $\mcF/ \mcF^1 = \varTheta^{\otimes (-n+1)}$ given by (\ref{Eq901}),
the sequence 
 (\ref{Eq648}) defined for the $\mr{SL}_n$-oper  $\mr{Sym}^{n-1} \msF_\odot^\heartsuit$ becomes    a short exact sequence
\begin{align} \label{Eq888}
0 \migi H^0 (X, \Omega\otimes_{\mcO_X} \varTheta^{\otimes (n-1)}(-D))
\migi 
H_{\mr{dR}, \mr{par}}^1 (X, \mr{Sym}^{n-1}\msF_\odot)
\migi H^1 (X, \varTheta^{\otimes (-n+1)})
\migi 0.
 \end{align}
Moreover, since (\ref{Eq803}) gives an identification $\Omega  = \varTheta^{\otimes 2}$,
we have 
$H^0 (X, \Omega \otimes_{\mcO_X} \varTheta^{\otimes (n-1)}(-D)) \cong H^0 (X, \varTheta^{\otimes (n+1)}(-D))$.
Hence,  
the assertion follows immediately from Proposition \ref{Prop15},  (ii).

Next, we shall prove the latter assertion of (ii).
To do this, we may assume
that the $\mr{SL}_2$-oper $\msF_\odot^\heartsuit$ comes, via the isomorphism ``$\Lambda_\vartheta^{\clubsuit \Rightarrow \diamondsuit}$" obtained in ~\cite[\S\,4.6.6, (553)]{Wak8},  from
a $(2, \vartheta)$-projective connection (cf. ~\cite[Definition 4.37, (ii)]{Wak8}), where $\vartheta$ denotes 
 a $2$-theta characteristic of $X^\mr{log}$ (cf.  ~\cite[Definition 4.31, (i)]{Wak8}) determined by $\varTheta$ in the manner of ~\cite[Example 4.34]{Wak8}.
Then, by considering the matrix form of the log connection $\nabla_\odot$ (cf. ~\cite[Remarks 4.30, 4.39]{Wak8}),
 the dual $\msF_\odot^\vee$ of $\msF_\odot$ is isomorphic to $\msF_\odot$ itself.
We  fix an  isomorphism $\msF_\odot \isom \msF_\odot^\vee$, which 
 induces an isomorphism
$\varpi : (\mcF, \nabla) \isom (\mcF^\vee, \nabla^\vee)$
  (cf. ~\cite[\S\,3]{Can}).

Now, let $U$ and $\iota$ be as in \S\,\ref{SS055g}, and
consider the $\mcO_X$-bilinear pairing 
\begin{align}
\widetilde{\varpi}^{\blacktriangleright} : \iota_* (\iota^* (\mcF)) \times  \iota_* (\iota^* (\Omega \otimes_{\mcO_X} \mcF))  \migi \iota_* (\iota^* (\Omega))
\end{align}
arising from the isomorphism $\varpi$.
In the following discussion, we shall prove (by an argument entirely similar to the proof of Proposition \ref{Prop99}) the claim that {\it $\widetilde{\varpi}^{\blacktriangleright}$ restricts to a nondegenerate $\mcO_X$-bilinear pairing}
\begin{align}
\varpi^{\blacktriangleright} : \Omega_\mr{par}^{0, c}(\mcF) \times \Omega_\mr{par}^1 (\mcF) \migi {^c}\Omega.
\end{align}
Since the  restriction of $\widetilde{\varpi}^{\blacktriangleright}$  to $U$ is nondegenerate,
the problem is  reduced to examining  the pairing $\widetilde{\varpi}^{\blacktriangleright}$  around   the marked points.
Let us choose $i \in \{1, \cdots, r \}$ and choose a local function  $t \in \mcO_X$ defining the closed subscheme  $\mr{Im}(\sigma_i) \subseteq X$.
Then, the formal neighborhood $Q$ of $\sigma_i$ is naturally isomorphic to $\mr{Spec}(k[\![t]\!])$, and we have  ${^c}\Omega |_Q = k[\![t]\!] dt$.
After  fixing a trivialization $\varTheta |_Q = k[\![t]\!]$, 
we can associate, to the local  function $t$, 
a canonical identification $\mcF |_Q = k[\![t]\!]^{\oplus n}$ 
(cf. ~\cite[Remark 4.30]{Wak8});
it restricts to  $\mcF^j |_Q = k[\![t]\!]^{\oplus n -j} \oplus \{ 0 \}^{\oplus j} \left( \subseteq k[\![t]\!]^{\oplus n} \right)$.
This identification 
 gives
 \begin{align}\label{Egio}
 \iota_* (\iota^* (\mcF)) |_{Q} = k (\!(t)\!)^{\oplus n} \
\left(\text{resp.,} \  \iota_* (\iota^* (\Omega \otimes_{\mcO_X} \mcF)) |_{Q} = (k (\!(t)\!) dt)^{\oplus n} \right).
 \end{align}
Moreover, by
the equality  $\mu_i^{\mr{Sym}^{n-1}\nabla_\odot} = q_{-1}$,
(\ref{Egio}) restricts to an identification
\begin{align} \label{Eqq291g}
 \Omega_\mr{par}^{0, c} (\mcF) |_Q = (k[\![t]\!]t)^{\oplus (n-1)} \oplus  k[\![t]\!]  
   \ \left(\text{resp.,} \ \Omega_\mr{par}^{1} (\mcF) |_Q = k[\![t]\!]dt \oplus  \left(k[\![t]\!] \frac{dt}{t}\right)^{\oplus (n-1)}\right).  
\end{align}
 According to this description, 
 each  section $v$ (resp., $u$) of 
$\Omega_\mr{par}^{0, c} (\mcF) |_Q$ (resp., $\Omega_\mr{par}^{1} (\mcF) |_Q$) may be described as the sum 
$v = t \cdot \sum_{i=1}^{n-1} v_i +  v_n$
 (resp., $u = u_1 + \frac{1}{t} \cdot \sum_{i=2}^n u_i$)
  for some $v_1, \cdots, v_n  \in k[\![t]\!]$
   (resp., $u_1, \cdots, u_n  \in k[\![t]\!] dt$). 
Then, up to a multiplicative constant factor, the following equalities hold: 
\begin{align}
\widetilde{\varpi}^\blacktriangleright (v, u) &=
\widetilde{\varpi}^\blacktriangleright \left(t \cdot \sum_{i=1}^{n-1} v_i,  \frac{1}{t} \cdot \sum_{i=2}^n u_i \right) + \widetilde{\varpi}^\blacktriangleright (v_n,  u_1) \\
 &= \widetilde{\varpi}^\blacktriangleright \left(\sum_{i=1}^{n-1} v_i,   \sum_{i=2}^n u_i\right) + \widetilde{\varpi}^\blacktriangleright (v_n,  u_1) \notag \\
 &=  \sum_{i=1}^n\widetilde{\varpi}^\blacktriangleright (v_i, u_{n+1 - i}). \notag
\end{align} 
It follows that  $\widetilde{\varpi}^\blacktriangleright (v, u)$ lies in  $k[\![t]\!] dt \left(= {^c}\Omega |_Q\right)$, and  that 
 the restriction of  $\widetilde{\varpi}^\blacktriangleright$  over $Q$ is nondegenerate.
This completes the proof of the claim.

Just as in the case of (\ref{Eqq1}),  the pairing $\varpi^\blacktriangleright$ together with the natural nondegenerate pairing $\varpi^\triangleright : \Omega_\mr{par}^{1, c}(\mcF) \times \Omega_\mr{par}^{1}(\mcF) \migi {^c}\Omega$ arising from $\varpi$ yields a $k$-bilinear pairing
\begin{align}
%\varpi_H : 
\mbH^1 (X, \mcK^\bullet [{^c}\nabla_\mr{par}])
\times \mbH^1 (X, \mcK^\bullet [\nabla_\mr{par}])
 \migi k.
\end{align}
The induced morphism $\varpi_H : \mbH^1 (X, \mcK^\bullet [{^c}\nabla_\mr{par}]) \migi \mbH^1 (X, \mcK^\bullet [\nabla_\mr{par}])^\vee$ fits into the morphism of exact sequences
\begin{align} \label{Eqqwer}
\vcenter{\xymatrix@C=11pt@R=36pt{
 H^0 (X, \Omega_\mr{par}^{0, c}(\mcF)) \ar[r] \ar[d]^-{\wr} & H^0 (X, \Omega_\mr{par}^{1, c} (\mcF)) \ar[r] \ar[d]^-{\wr} &  \mbH^1 (X, \mcK^\bullet [{^c}\nabla_\mr{par}]) \ar[d]^-{\varpi_{H}}\ar[r]  & H^1 (X, \Omega_\mr{par}^{0, c}(\mcF)) \ar[d]^-{\wr} \ar[r] & H^1 (X, \Omega_\mr{par}^{1, c} (\mcF)) \ar[d]^-{\wr} 
\\
H^1 (X, \Omega_\mr{par}^1 (\mcF))^\vee \ar[r] &H^1 (X, \Omega_\mr{par}^0 (\mcF))^\vee \ar[r]&    \mbH^1 (X, \mcK^\bullet [\nabla_\mr{par}])^\vee \ar[r] & H^0 (X, \Omega_\mr{par}^1 (\mcF))^\vee \ar[r] & H^0 (X, \Omega_\mr{par}^0 (\mcF))^\vee
}}  
\end{align}
(cf. (\ref{Eqq301}) for a similar diagram), where
the horizontal sequences arise from the Hodge to de Rham spectral sequences involved, and the vertical arrows except for the middle one arise from the pairings $\varpi^\blacktriangleright$ and $\varpi^\triangleright$.
By applying the five lemma to this diagram, we see that the middle vertical arrow $\varpi_H$ is an isomorphism.
Under  the equality $\mbH^1 (X, \mcK^\bullet [{^c}\nabla_\mr{par}]) = H^1_{\mr{dR}, \mr{par}} (X, \mr{Sym}^{n-1}\mcF_\odot) \left(= \mbH^1 (X, \mcK^\bullet [\nabla_\mr{par}])  \right)$ resulting from  Proposition \ref{Prp7},
$\varpi_H$ defines an isomorphism 
\begin{align} \label{EHju}
H^1_{\mr{dR}, \mr{par}} (X, \mr{Sym}^{n-1}\msF_\odot) \isom H^1_{\mr{dR}, \mr{par}} (X, \mr{Sym}^{n-1}\msF_\odot)^\vee.
\end{align}
Moreover,  by  the various definitions involved and the local description of $\widetilde{\varpi}^\blacktriangleright$ described above,
we see that  that the short exact sequence (\ref{Eq810})  is compatible with its dual via (\ref{EHju}).
This completes the proof of this proposition.
\end{proof}
%----------------------------------------------------------------------------------------------------
\SSP

Next,  if $n = 2l+1$ for some positive integer $l$, then the following assertion holds.
(We leave to the reader the formulation of relevant duality as asserted in the above theorem.)

\SSP
%----------------------------------------------------------------------------------------------------
\bt \label{Th66}
Let   $\msE^\spadesuit$   be
 a $G^\odot$-oper on $\msX$.
 Recall from the discussion in \S\,\ref{SS03f61} that the $2l$-th symmetric product $\mr{Sym}^{2l}\msE_\odot^\spadesuit$ of  $\msE^\spadesuit$ can be obtained  (cf. (\ref{Eq9d00})).
 Then, the following assertions hold:
 \begin{itemize}
 \item[(i)]
There exists  a canonical short exact sequence of $k$-vector spaces
\begin{align} \label{Eq817}
0 \longmigi H^0 (X, \Omega^{\otimes (l+1)})
%\xrightarrow{h_\sharp^l [\msE_\odot^\spadesuit]}
\longmigi
H^1_{\mr{dR}} (X, \mr{Sym}^{2l}\msE)
%\xrightarrow{h_\flat^l [\msE_\odot^\spadesuit]}
\longmigi
 H^1 (X, \Omega^{\otimes (-l)})
\longmigi 0,
\end{align}
where $\mr{Sym}^{2l}\msE_\cdot$ denotes the underlying flat $\mr{SL}_{2l}$-bundle of $\mr{Sym}^{2l}\msE^\spadesuit_\cdot$.

\item[(ii)]
Suppose further that 
$\msE^\spadesuit_\odot$ is of radii $[0]^{\times r} \in \mfc (k)^{\times r}$ (where $[0]^{\times r} := \emptyset$ if $r = 0$).
Then, 
there exists a canonical short exact sequence
 \begin{align} \label{Eq810F}
0 \longmigi H^0 (X, \Omega^{\otimes (l+1)}(-D))
%\xrightarrow{h_{\sharp, \mr{par}}^l [\msE_\odot^\spadesuit]}
\longmigi
H^1_{\mr{dR}, \mr{par}} (X, \mr{Sym}^{2l}\msE)
%\xrightarrow{h_{\flat, \mr{par}}^l [\msE_\odot^\spadesuit]}
\longmigi 
 H^1 (X, \Omega^{\otimes (-l)})
\longmigi 0. 
\end{align}
\end{itemize}
Finally, when $\msE^\spadesuit$ comes from an $\mr{SL}_2$-oper via projectivization,  these short exact sequences coincide with the corresponding ones in the statement of Theorem \ref{Th67}.
\et
%----------------------------------------------------------------------------------------------------
\begin{proof}
Assertion (i) (resp., (ii)) follows from Proposition \ref{Prop15}, (i) (resp., (ii)) and an argument similar to the proof of Theorem \ref{Th67}, (i) (resp., (ii)). 
\end{proof}
%----------------------------------------------------------------------------------------------------
%\SSP

\LSP
%----------------------------------------------------------------------[begin subsection]-------------
\subsection{The irreducible decomposition of the adjoint bundle} \label{SSa4}

Let $\msE^\spadesuit_\odot := (\mcE_{B^\odot}, \nabla_\odot)$ be a $\,\qq$-normal  $\mr{PGL}_2$-oper.
We shall write $\msE^\spadesuit \left(= (\mcE_B, \nabla) \right) :=\iota_{G*}(\msE^\spadesuit_\odot)$, which is a $G$-oper.
Also,   write $\mcE_\odot := \mcE_{B^\odot} \times^{B^\odot} G^\odot$ and  $\mcE := \mcE_B \times^B G$.
In what follows, we consider a decomposition of the flat adjoint  bundle $(\mfg_\mcE, \nabla^\mr{ad})$ according to the irreducible decomposition of $\mfg$, regarded as  an $\mfs \mfl_2$-module   via $\iota_\mfg$.

Let us fix an integer $l$ with $- \mr{rk}(\mfg) \leq l \leq \mr{rk} (\mfg)$.
Under the natural  identification 
$\Omega^{\otimes l}\otimes_k \mfg_l^{\mr{ad}(q_1)}= \mcV_{G, l} \left(\subseteq \mfg_\mcE \right)$,
 each nonzero element $a \in \mfg_l^{\mr{ad}(q_1)}$ determines an inclusion $\gamma_a : \Omega^{\otimes l} \migiincl \mfg_\mcE$.
Note that, for each $j = 0, \cdots, 2l+1$,  the $k$-vector subspace
\begin{align}
H_a^j := \sum_{s=0}^{2l-j} k \cdot  \mr{ad}(q_{-1})^s(a) \left(= \bigoplus_{s=0}^{2l-j} k \cdot  \mr{ad}(q_{-1})^s(a) \right)
\end{align}
 of $\mfg$ is
 closed under the
 adjoint $B^\odot$-action via $\iota_B$, and  that
  $H_a^0$ 
  defines  an $\mfs \mfl_2$-submodule  of $\mfg$.
Hence, $H_a^0$ determines a rank $(2l+1)$  flat  subbundle $(\mcH_a, \nabla_a)$ of $(\mfg_\mcE, \nabla^\mr{ad})$, and  
the subspaces  $\{ H_a^{j}\}_{j=0}^{2l+1}$ of  $H_a^0$ determine  a decreasing filtration $\{ \mcH_a^j \}_{j=0}^{2l+1}$ on  $\mcH_a$.
 One may verify  that the resulting collection
 \begin{align} \label{Eqq100}
 \msH_a^\heartsuit  := (\mcH_a, \nabla_a, \{ \mcH_a^j \}_{j=0}^{2l+1})
 \end{align}
 form an $\mr{SL}_{2l+1}$-oper on $\msX$.
 %, and that it is isomorphic to $\mr{Sym}^{2l}\msE^\spadesuit_\odot$.

Next, let $\mcD$ (resp., $\mcD_{< m}$ for each nonnegative integer $m$) denote the sheaf  of logarithmic differential operators (resp., logarithmic differential operators of order $< m$) of $X^\mr{log}$ (cf. ~\cite[\S\,4.2.1]{Wak8}).
Then, 
  it follows from the definition of an $\mr{SL}_{2l+1}$-oper that  the composite
 \begin{align}
\widetilde{\gamma}_a  : \mcD_{< 2l+1} \otimes_{\mcO_X} \Omega^{\otimes l} 
\longmigi
  \mcD \otimes_{\mcO_X}  \mcH_a^0 \longmigi \mcH_a  \hspace{20mm} \\
   \left(\text{resp.,} \  \widetilde{\gamma}_{\mr{sym}} : 
  \mcD_{< 2l+1} \otimes_{\mcO_X}\Omega^{\otimes l}
  \longmigi
  \mcD \otimes_{\mcO_X}  \mr{Sym}^{2l}\mcE \longmigi \mr{Sym}^{2l}\mcE_\odot \notag
   \right)
 \end{align}
 is an isomorphism, where 
 the first arrow denotes the injection arising from both the natural  inclusion $\mcD_{< 2l+1} \migiincl \mcD$ and $\gamma_a$ (resp., the injection $\Omega^{\otimes l} \migiincl  \mr{Sym}^{2l}\mcE_\odot$ deduced from  (\ref{Eq913})), and 
 the second arrow arises from the $\mcD$-module structure on $\mcH_a^0$ (resp.,  $\mr{Sym}^{2l}\mcE_\odot$) determined by $\nabla_a$ (resp., $\mr{Sym}^{2l}\nabla_\odot$).
The composite
 \begin{align}
 \mr{Sym}^{2l} \mcE_\odot \xrightarrow{\widetilde{\gamma}_{\mr{sym}}^{-1}} \mcD_{< 2l+1} \otimes_{\mcO_X} \Omega^{\otimes l} \xrightarrow{\widetilde{\gamma}_a} \mcH_a \xrightarrow{\mr{inclusion}} \mfg_\mcE
 \end{align}
 is compatible with the respective connections, i.e., $\mr{Sym}^{2l}\nabla_\odot$ and $\nabla^\mr{ad}$.
That is,  we have obtained  an injective morphism of flat bundles
$\xi_a : \mr{Sym}^{2l} \msE_\odot  \migiincl  (\mfg_\mcE, \nabla^\mr{ad})$.
The morphism of flat bundles
\begin{align}
\xi_l : \mr{Sym}^{2l} \msE_\odot \otimes_k \mfg_l^{\mr{ad}(q_1)} \migiincl  (\mfg_\mcE, \nabla^\mr{ad})
\end{align}
given by  $v \otimes a \mapsto \xi_a (v)$ for  $a \in \mfg_l^{\mr{ad}(q_1)}$ and $v \in \mr{Sym}^{2l}\mcE_\odot$ is well-defined, and the various $\xi_l$'s define   an isomorphism 
\begin{align} \label{Eq160}
\bigoplus_{l=1}^{\mr{rk}(\mfg)} \xi_l : \bigoplus_{l=1}^{\mr{rk}(\mfg)} \mr{Sym}^{2l}\msE_\odot \otimes_k \mfg^{\mr{ad}(q_1)}_l \isom 
(\mfg_\mcE, \nabla^\mr{ad}). 
\end{align}
By applying the $1$-st hypercohomology functor  to this isomorphism,
we obtain an isomorphism of $k$-vector spaces
\begin{align} \label{Eq834}
\bigoplus_{l=1}^{\mr{rk}(\mfg)} H^1_{\mr{dR}} (X, \mr{Sym}^{2l}\msE_\odot) \otimes_k \mfg^{\mr{ad}(q_1)}_l \isom \mbH^1 (X,  \mcK^\bullet [\nabla^\mr{ad}]).
\end{align}
This isomorphism fits into the  isomorphism of short exact sequences
\begin{align} \label{Eq16tt}
\vcenter{\xymatrix@C=11pt@R=36pt{
0 \ar[r] & {\displaystyle \bigoplus_{l=1}^{\mr{rk}(\mfg)}} H^0 (X, \Omega^{\otimes (l+1)})
 \ar[r]
 %^-{\bigoplus_l h^l_\sharp [\msE^\spadesuit_\odot]}
  \ar[d]^-{\wr}
  & 
 {\displaystyle \bigoplus_{l=1}^{\mr{rk}(\mfg)}}
 H^1_{\mr{dR}} (X, \mr{Sym}^{2l}\msE) \otimes_k \mfg_l^{\mr{ad}(q_1)} \ar[r]
 %^-{\bigoplus_l h^l_\flat [\msE^\spadesuit_\odot]}
  \ar[d]_-{\wr}^-{(\ref{Eq834})} & 
 {\displaystyle \bigoplus_{l=1}^{\mr{rk}(\mfg)}}
  H^1 (X, \Omega^{\otimes (-l)}) \ar[r] \ar[d]^-{\wr} & 0 
\\
0 \ar[r] & H^0 (X, \DV_G)\ar[r]_-{{'}e_\sharp} & \mbH^1 (X, \mcK^\bullet [\nabla^\mr{ad}]) \ar[r]_-{{'}e_\flat} & H^1 (X, \Omega \otimes_{\mcO_X} \DV_G^\vee) \ar[r]  & 0 
}}
\end{align}
 (cf. (\ref{Eq817}) for the upper horizontal sequence and  (\ref{HDeq1}) for the lower horizontal sequence), where the left-hand and right-hand vertical arrows are the isomorphism induced from $\DV_G \cong \bigoplus_{l=1}^{\mr{rk}(\mfg)} \Omega \otimes_k \mfg^{\mr{ad}(q_1)}$ (cf. (\ref{QR020})).

Moreover, suppose that $\msE^\spadesuit_\odot$ is of radii $[0]^{\times r}$ (where $[0]^{\times r} := \emptyset$ if $r = 0$).
Then, (\ref{Eq160}) yields an isomorphism
\begin{align} \label{Eq83f4}
\bigoplus_{l=1}^{\mr{rk}(\mfg)} H^1_{\mr{dR}, \mr{par}} (X, \mr{Sym}^{2l}\msE_\odot) \otimes_k \mfg^{\mr{ad}(q_1)}_l \isom \mbH^1 (X,  \mcK^\bullet [\nabla_\mr{par}^\mr{ad}]),
\end{align}
which fits into the isomorphism of short exact sequences
\begin{align} \label{Eq16ttt}
\vcenter{\xymatrix@C=11pt@R=36pt{
0 \ar[r] & {\displaystyle \bigoplus_{l=1}^{\mr{rk}(\mfg)}} H^0 (X, \Omega^{\otimes (l+1)}(-D))
\ar[r]
 \ar[d]^-{\wr}
  & 
 {\displaystyle \bigoplus_{l=1}^{\mr{rk}(\mfg)}}
 H^1_{\mr{dR}, \mr{par}} (X, \mr{Sym}^{2l}\msE_\odot) \otimes_k \mfg_l^{\mr{ad}(q_1)} 
 \ar[r]
 \ar[d]_-{\wr}^-{(\ref{Eq83f4})} & 
 {\displaystyle \bigoplus_{l=1}^{\mr{rk}(\mfg)}}
  H^1 (X, \Omega^{\otimes (-l)}) \ar[r] \ar[d]^-{\wr} & 0 
\\
0 \ar[r] & H^0 (X, {^c}\DV_G)\ar[r]_-{{'}e_{\sharp, \mr{par}}} & \mbH^1 (X, \mcK^\bullet [\nabla_\mr{par}^\mr{ad}]) \ar[r]_-{{'}e_{\flat, \mr{par}} } & H^1 (X, \Omega \otimes_{\mcO_X} \DV_G^\vee) \ar[r]  & 0 
}} 
\end{align}
 (cf. (\ref{Eq810}) for the upper horizontal sequence and  (\ref{Eqr}) for the lower horizontal sequence).
Once we have constructed a canonical self-duality for $ H^1_{\mr{dR}, \mr{par}} (X, \mr{Sym}^{2l}\msE_\odot)$ as asserted  in Theorem \ref{Th67}, (ii),  this duality  will be  verified to be   compatible with 
(\ref{Eqqw9})
 via (\ref{Eq83f4}).

%%%%%%%%%%%%%%%%%%%%%%%%%%%%%%%%%%%%%%%%%%%%%
%%%%%%%%%%%%%---[begin section]---%%%%%%%%%%%%%%%%%%%%%%%
\vspace{10mm}
\section{Dormant opers and their deformations}\SSP

This section discusses opers in positive characteristic, especially, dormant opers (= opers with vanishing $p$-curvature).
The study of the moduli space of dormant opers  was substantially developed  in  ~\cite{Jo14}, ~\cite{JP}, ~\cite{Wak}, ~\cite{Wak2}, and  ~\cite{Wak8}. 

One of the important features of   this moduli space 
 is the generic \'{e}taleness (cf. ~\cite[Theorem G]{Wak8}).
In other words, if the underlying 
 curve is sufficiently general (and $G$ satisfies a certain  additional assumption), then  the moduli spaces of opers and $p$-flat bundles  intersect transversally.
This fact  
  will be essential in the proof of our main theorem, as discussed in the next section.

In the rest of the present paper, 
we  assume that
$\mr{char}(k)= p> 0$ for a prime $p$, and that
 either ``$p >2h_G$" or ``$G = \mr{PGL}_n$ with $1 < n < p$" is fulfilled.

\LSP
%----------------------------------------------------------------------[begin subsection]-------------
\subsection{Dormant opers on pointed stable curves} \label{SS02}

%Denote by $X^{(1)}$ the Frobenius twist of $X$ over $k$ and by $F_{X/k}$ the relative Frobenius morphism $X \migi X^{(1)}$.
%Write $\Omega^{(1)} := \Omega_{X^{(1)\mr{log}}}$.
 
Let $(\mcE, \nabla)$ be   a flat $G$-bundle on $\msX$.
Recall (cf. ~\cite[Definition 3.8]{Wak8}) that the $p$-curvature of $\nabla$ is  the $\mcO_X$-linear morphism ${^p}\psi^\nabla : \mcT^{\otimes p} \migi \widetilde{\mcT}_{\mcE^\mr{log}}$ uniquely determined by the condition that any local section in $\mcT^{\otimes p}$ of the form  $\partial^{\otimes p}$ for some $\partial \in \mcT$ is mapped to $\nabla (\partial)^{[p]}- \nabla (\partial^{[\partial]})$, where $(-)^{[p]}$ denotes the result of applying the  $p$-power  operations (cf. ~\cite[\S\,3.2]{Wak8}).
In the case where $G$ is a matrix group, this morphism is equivalent to the classical definition of $p$-curvature described in terms of vector bundles (cf. e.g., ~\cite[\S\,5]{Kal}).

By a {\bf $p$-flat $G$-bundle}, we mean a flat $G$-bundle $(\mcE, \nabla)$ such that $\nabla$ has  vanishing $p$-curvature.
We shall denote by 
\begin{align} \label{Eq220}
\C_G^{\psi  = 0}  \ \left(\text{resp.,} \  \C_{G, \mu}^{\psi  = 0} \ \text{for each $\mu \in \mfg^{\times r}$} \right)
\end{align}
 the closed substack of $\C_G$ (resp., $\C_{G, \mu}$)  classifying $p$-flat $G$-bundles.

Next, 
a $G$-oper $\msE^\spadesuit := (\mcE_B, \nabla)$
 is called  {\bf dormant}  if $\nabla$ has vanishing $p$-curvature (cf. ~\cite[Definition 3.15]{Wak8}).
Then, we obtain a closed  subscheme 
\begin{align} \label{Eq212}
\mcO p^{^\mr{Zzz...}}_{G} \  \left(\text{resp.}, \ \mcO p^{^\mr{Zzz...}}_{G, \rho}  \right)
\end{align}
of $\mcO p_{G}$ (resp., $ \mcO p_{G, \rho}$ for each $\rho \in \mfc (k)^{\times r}$) classifying
dormant $G$-opers.
The following equality between substacks of $\C_G$ (resp., $\C_{G, \widetilde{\rho}}$) holds:
\begin{align} \label{Eq211}
\mcO p_G^{^\mr{Zzz...}} = \mcO p_G \cap \C_G^{\psi = 0}
\ \left(\text{resp.,} \ 
\mcO p_{G, \rho}^{^\mr{Zzz...}} = \mcO p_{G, \rho} \cap \C_{G, \widetilde{\rho}}^{\psi = 0} \right).
\end{align} 
Since $G$ has been assumed to be   of adjoint type, 
we can apply  ~\cite[Theorem C]{Wak8} to see  that
    $\mcO p^{^\mr{Zzz...}}_G$ (resp., $\mcO p^{^\mr{Zzz...}}_{G, \rho}$)  is a nonempty finite $k$-scheme (resp., a possibly empty finite $k$-scheme).

Now, let $\msE^\spadesuit_\odot := (\mcE_{B^\odot}, \nabla_\odot)$ be a $G^\odot$-oper on $\msX$.
The log connection $\iota_{G*}(\nabla)$ of the $G$-bundle 
$\mcE_G := \mcE_{B^\odot} \times^{B^\odot, \iota_B}G$ (cf. (\ref{associated}))
 satisfies the equality
${^p}\psi^{\iota_{G*}(\nabla)} = (d \iota_G)_\mcE \circ {^p}\psi^{\nabla}$, where $(d \iota_{G})_\mcE$  denotes the natural  injection $\widetilde{\mcT}_{\mcE^\mr{log}_{G^\odot}} \migiincl  \widetilde{\mcT}_{\mcE_{G}^\mr{log}}$ 
 (cf. ~\cite[\S\,3.3.2]{Wak8}).
In particular,  $\msE^\spadesuit_{\odot}$ is dormant if and only if the associated $G$-oper $\iota_{G*}(\msE^\spadesuit_\odot)$ is dormant.
It follows that
(\ref{Eq405})
%  $\iota^{\mcO p}_{G}$
  (resp., 
  (\ref{Eq281}))
%  $\iota^{\mcO p}_{G, \rho_\odot}$)
   restricts to a closed immersion
\begin{align} \label{E2}
%{^p}\iota_G^{\mcO p} : 
\mcO p^{^\mr{Zzz...}}_{G^\odot} \migiincl  \mcO p^{^\mr{Zzz...}}_{G} \ \left(\text{resp.,} \
 %{^p}\iota_{G, \rho_\odot}^{\mcO p} :
  \mcO p^{^\mr{Zzz...}}_{G^\odot, \rho_\odot} \migiincl  \mcO p^{^\mr{Zzz...}}_{G, \iota_\mfc (\rho_\odot)} \right).
\end{align}

\LSP
%---------------------------[begin subsection]-------------
\subsection{Example: a dormant $\mr{PGL}_2$-oper on a Shimura curve} \label{SS0123}

In this subsection, 
we shall illustrate  an example of a dormant $\mr{PGL}_2$-oper  on a Shimura curve resulting from the discussions by H. Reimann (cf. ~\cite{Rei}) and  M. Sheng, J. Zhang, K. Zuo (cf. ~\cite{SZZ}).
See ~\cite{Wak10} for other examples of  dormant $\mr{PGL}_2$-opers,  constructed by using the Gauss maps for  Fermat curves.

Let 
$F$ be a  totally real number field of degree $\geq 2$ in which  $p$ is unramified.
Also,  let $A$
 be a  quaternion  algebra over $F$
 which is split at one infinite place $\tau : F \migiincl \mbR$ and ramified at all remaining infinite places.
  After fixing an embedding $\overline{\mbQ} \migiincl \overline{\mbQ}_p$, one obtains an bijection $\mr{Hom}_\mbQ (F, \overline{\mbQ})$ and $\coprod_{i=1}^r \mr{Hom}_{\mbQ_p} (F_{\mfp_i}, \overline{\mbQ}_p)$, where $\mfp_1 \left(=:\mfp \right), \mfp_2 , \mfp_3, \cdots, \mfp_r$ are all primes in $F$ dividing  $p$ such that $\tau$ lies in $\mr{Hom}_{\mbQ_p} (F_{\mfp}, \overline{\mbQ}_p)$.
Also, we shall choose a pair $(L, K)$ consisting of  a totally imaginary quadratic extension $L$ of $F$ contained in $A$ such that all $\mfp_i$ stay prime in $L$, and 
an imaginary quadratic extension 
$K$  of $F$
  in which 
 all the primes $\mfp_1, \cdots, \mfp_r$ split.
By regarding both  $F^\times$ and $K^\times$ as $\mbQ$-groups,
 one may obtain a  $\mbQ$-group $G$ which makes  the following  sequence exact:
\begin{align}
1 \longmigi F^\times  \xrightarrow{f \mapsto (f, f^{-1})} A^\times \times K^\times \longmigi G \longmigi 1.
\end{align}

Suppose $A$ is split at $\mfp_1$, and write
$\mcO_B := \mcO_A \otimes_{\mcO_F} \mcO_K$, where $\mcO_A$ is an order of $A$ containing $\mcO_L$ with certain additional properties (cf. ~\cite[\S\,2, (2.9)]{Rei}).
Then, 
for every level structure $C = C_p \times C^p \subseteq G (\mbA_f)$ (where  $\mbA_f$ denotes the ring of finite ad\`{e}les of $\mbQ$) with $C_p = G(\mbZ_p)$ and $C^p$ sufficiently small (cf. ~\cite[\S\,2, (2.11)]{Rei}),
there exists a smooth proper $\mcO_{F_\mfp}$-scheme $\mcM_C$ of relative dimension $1$ which is 
 the coarse moduli space of a certain moduli functor of PEL type with the endomorphism algebra $\mcO_B$ (cf.   ~\cite[Definition 2.12, Proposition 2.14, and Corollary 3.14]{Rei}).

Let us take  
 one of the 
 connected components
 of the base-change of $\mcM_C$ to 
 $k$ (via a fixed inclusions $\mcO_{F_\mfp}/p\mcO_{F_\mfp}\migiincl \overline{\mbF}_p \migiincl k$),
  which we denote by $X$.
 After possibly replacing  $C$ with another, we may assume  that  the genus of $X$ is  greater than $1$ and that   there exists a universal abelian scheme  
 $f : Y \migi X$.
The $1$-st  de Rham cohomology sheaf  $\mcH := R^1 f_* (\mcK^\bullet [\mcO_Y \xrightarrow{d}\Omega_{Y /X}])$  on $X$ is equipped with 
the Gauss-Manin connection $\nabla$, as well as  
a $2$-step decreasing filtration $\{ \mcH^j \}_{j=0}^2$
on $\mcH$ 
arising from the Hodge filtration.

If we write $\mfp_i \mcO_K = \mfq_i \overline{\mfq}_i$ ($i=1, \cdots, r$)  for  distinct primes $\mfq_i$, $\overline{\mfq}_i$ in $K$, then
we have $p \mcO_{LK} = \prod_{i=1}^r \mfq_i \overline{\mfq}_i$; it  induces a natural isommorphism  $LK \otimes_\mbQ \mbQ_p \isom \prod_{i=1}^r LK_{\mfq_i} \times LK_{\overline{\mfq}_i}$.
This isomorphism tensored with $\mbQ^{\mr{un}}_p$ (= the maximal unramified extension of $\mbQ_p$) restricts to a decomposition
\begin{align} \label{Eq201rt}
\prod_{\phi \in \Phi} w (\phi) \times \overline{w} (\phi) : \mcO_{LK} \otimes_\mbZ W (\overline{\mbF}_p) \isom \prod_{\phi \in \Phi} W (\overline{\mbF}_p) \times W (\overline{\mbF}_p),
\end{align}
where $\Phi := \mr{Hom}_{\mbQ} (L, \overline{\mbQ})$ and 
for every $\phi \in \Phi$ we denote by $w (\phi)$ (resp., $\overline{w} (\phi)$) the unique extension of $\phi$ in  
$\mr{Hom}_\mbQ (LK, \overline{\mbQ})$
whose restriction 
to $K$
 lies $\mr{Hom}_{\mbQ_p}(K_{\mfq_i}, \overline{\mbQ}_p)$
 (resp., $\mr{Hom}_{\mbQ_p}(K_{\overline{\mfq}_i}, \overline{\mbQ}_p)$)
  for some $i$.
By the injection $\left(\mcO_{LK} \subseteq \right) \mcO_B \migiincl \mr{End}_X (Y)$,
$\mcH$ is equipped with a structure of  $\mcO_{LK}\otimes_\mbZ W (\overline{\mbF}_p)$-module.
Hence
the  decomposition (\ref{Eq201rt}) gives rise to a decomposition of  $(\mcH, \nabla)$ 
 into a direct sum of rank $2$ flat subbundles
\begin{align} \label{Eq107}
(\mcH, \nabla) = \bigoplus_{\phi \in \Phi} (\mcH_{\phi}, \nabla_\phi) \oplus (\mcH_{\overline{\phi}}, \nabla_{\overline{\phi}}),
\end{align}
(cf. ~\cite[Proposition 4.1]{SZZ});
it restricts to a decomposition $\mcH^1 = \bigoplus_{\phi \in \Phi} \mcH^1_\phi \oplus \mcH_{\overline{\phi}}^1$, where $\mcH_{\phi}^{1} := \mcH^{1} \cap \mcH_{ \phi}$ and $\mcH_{\overline{\phi}}^{1} := \mcH^{1} \cap \mcH_{\overline{\phi}}$.
Also, 
the structure of Dieudonn\'{e} crystal on $\mcH$ corresponding to the abelian scheme 
$Y/X$
 determines, 
 via reduction modulo $p$,
an $\mcO_{X}$-linear  morphism
$\xi : F_{X}^* \mcH \migi \mcH$ (where $F_X$ denotes the absolute Frobenius endomorphism of $X$).
Let us take $\widetilde{\tau} \in \Phi$ with $\widetilde{\tau} |_{F} = \tau$.
If $\sigma$ denotes the automorphism induced by the Frobenius automorphism of $\overline{\mbF}_p$ over $\mbF_p$, then 
$\xi$ restricts to an isomorphism
\begin{align}
\xi_{\widetilde{\tau}} : F_{X}^*\mcH_{\sigma^{-1}\widetilde{\tau}} \isom \mcH_{\widetilde{\tau}}
\end{align}
(cf. ~\cite[Proposition 3.1, (iii)]{SZZ}).
Since this isomorphism is compatible with the connection,
we see that $\nabla_{\widetilde{\tau}}$  has vanishing $p$-curvature.
By the latter assertion of ~\cite[Proposition 4.1]{SZZ},
 both $\mcH_{\widetilde{\tau}}^{1}$ and $\mcH_{\widetilde{\tau}}/\mcH_{\widetilde{\tau}}^{1}$  are line bundles.
Moreover, if $p \geq 2g$, then 
the structure of Higgs field defined as the Kodaira-Spencer map 
$\mcH_{\widetilde{\tau}}^{1} \migi \Omega_{X} \otimes_{\mcO_X} \left(\mcH_{\widetilde{\tau}}/\mcH_{\widetilde{\tau}}^{1}\right)$ is an isomorphism (cf. ~\cite[Proposition 4.4]{SZZ}).
It follows that the collection
\begin{align}
(\mcH_{\widetilde{\tau}}, \nabla_{\widetilde{\tau}}, \{ \mcH^j_{\widetilde{\tau}}  \}_{j=0}^2)
\end{align}
specifies a dormant $\mr{GL}_2$-oper.
In particular, by taking its projectivization, one can obtain  a dormant $\mr{PGL}_2$-oper  on  $X$.

\LSP
%---------------------------[begin subsection]-------------
\subsection{The deformation space of a dormant $G$-oper} \label{SS030}

Let $(\mcE, \nabla)$ be a flat $G$-bundle on $\msX$ with vanishing $p$-curvature.
In the resp'd portion of the following discussion, we suppose that
$(\mcE, \nabla)$ is equipped with a marking, and that $\nabla$ has residues $\mu \in \mfg^{\times r}$ (where $\mu := \emptyset$ if $r = 0$).

Consider the conjugate spectral sequence
\begin{align}\label{Csp}
''E_{2}^{a,b} := H^a(X, \mcH^b(\mcK^\bullet [\nabla^\mr{ad}])) \Rightarrow \mbH^{a+b}(X, \mcK^\bullet [\nabla^\mr{ad}])  \hspace{4mm} \\
\left(\text{resp.,} \ ''E_{2}^{a,b} := H^a(X, \mcH^b(\mcK^\bullet [{^c}\nabla^\mr{ad}])) \Rightarrow \mbH^{a+b}(X, \mcK^\bullet [{^c}\nabla^\mr{ad}])  \right)
\end{align}
associated to  $\mcK^\bullet [\nabla^\mr{ad}]$ (resp., $\mcK^\bullet [{^c}\nabla^\mr{ad}]$), 
where $\mcH^b(\mcK^\bullet [-])$ denotes the $b$-th cohomology sheaf of the complex $\mcK^\bullet [-]$.
This spectral sequence 
 induces a  short exact sequence
\begin{align} 
\label{Ceq}
   0 \migi H^1(X, \mr{Ker}(\nabla^\mr{ad})) \xrightarrow{{''e_\sharp}}
 \mbH^1(X, \mcK^\bullet[\nabla^\mr{ad}]) \xrightarrow{{''e}_\flat}
  H^0 (X, \mr{Coker}(\nabla^\mr{ad})) \migi 0  \hspace{10mm}\\
  \left(\text{resp.,} \  0 \migi H^1(X, \mr{Ker}({^c}\nabla^\mr{ad})) \xrightarrow{{''e_{\sharp, c}}}
 \mbH^1(X, \mcK^\bullet[{^c}\nabla^\mr{ad}]) \xrightarrow{{''e}_{\flat, c}}
  H^0 (X, \mr{Coker}({^c}\nabla^\mr{ad})) \migi 0 \right). \notag
\end{align}

  Recall that the deformation of $p$-curvature may be described as the Cartier map 
   of the flat bundle $(\mfg_\mcE, \nabla^\mr{ad})$  (cf.   ~\cite[Proposition 6.11]{Wak8}), and that 
   this map (resp., this map restricted to ${^c}\mfg_\mcE$)
    factors through the quotient $\Omega \otimes_{\mcO_X} \mfg_\mcE \migisurj \mr{Coker}(\nabla^\mr{ad})$ (resp., the quotient $\Omega \otimes_{\mcO_X} {^c}\mfg_\mcE \migisurj \mr{Coker}({^c}\nabla^\mr{ad})$) (cf. ~\cite[Proposition 2.2.4]{Og}).
  Hence, 
the space of first-order deformations of the flat $G$-bundle $(\mcE, \nabla)$  preserving the condition of vanishing $p$-curvature (resp., the condition of vanishing $p$-curvature and fixing residues) may be identified with $H^1 (X, \mr{Ker}(\nabla^\mr{ad}))$ (resp., $H^1 (X, \mr{Ker}({^c}\nabla^\mr{ad}))$).
That  is, 
the isomorphism $\gamma$ (resp., $\gamma_\mu$)
 restricts, via ${''}e_\sharp$ (resp., ${''}e_{\sharp, c}$),  to an 
 isomorphism  of $k$-vector spaces 
\begin{align} \label{Eq303}
{^p} \gamma : T_q \C_G^{\psi = 0} \isom 
H^1 (X, \mr{Ker}(\nabla^\mr{ad}))
 \ \left(\text{resp.,} \  {^p} \gamma_\mu : T_q \C_{G, \mu}^{\psi = 0} \isom H^1 (X, \mr{Ker}({^c}\nabla^\mr{ad}))
  \right).
\end{align}

Next, suppose further that there exists a structure of $B$-reduction $\mcE_B$ on $\mcE$ for which $\msE^\spadesuit := (\mcE_B, \nabla)$ specifies a  dormant $\,\qq$-normal  $G$-oper (resp., a dormant $\,\qq$-normal $G$-oper of radii $\rho := \chi (\mu)$, where $\rho := \emptyset$ if $r = 0$).
Denote by $q$ the $k$-rational point of $\mcO p_G^{^\mr{Zzz...}}$ (resp., $\mcO p_{G, \rho}^{^\mr{Zzz...}}$) classifying $\msE^\spadesuit$;
we use the same symbol ``$q$" to denote the point of 
 $\C_G^{\psi = 0}$ (resp., $\C_{G, \widetilde{\rho}}^{\psi = 0}$) classifying $(\mcE, \nabla)$.

We already know the following assertion.

\SSP
%-----------------------------------------------------------------------------
\bpr \label{Eqw3}
The following equalities hold:
\begin{align} \label{WWW8}
\mr{dim}(T_q \C_{G}^{\psi = 0}) &=
2 \cdot \mr{dim}(T_q \mcO p_G) - r \cdot \mr{rk}(\mfg)\\
& = 2 \cdot \mr{dim}(T_q \C_{G}^{\psi = 0}) + r \cdot \mr{rk}(\mfg) \notag \\
& =   (2g-2 +r)\cdot \mr{dim}(\mfg). \notag
\end{align}
In particular, we have 
\begin{align} \label{WWW9}
\mr{dim}(T_q \mcO p_G) +  \mr{dim}(T_q \C_{G}^{\psi = 0})  =  \mr{dim}(T_q \C_{G}^{\psi = 0}).
\end{align}
\epr
%-----------------------------------------------------------------------------
\begin{proof}
The assertion follows from  (\ref{Eq3033}), (\ref{Eq300}), (\ref{Eq303}), and  ~\cite[Proposition 2.23, 6.5,  and 6.18]{Wak8}.
\end{proof}
%-----------------------------------------------------------------------------
\SSP

We shall set
\begin{align} \label{Eq302}
e_\flat^{(0)} : H^0 (X, \DV_G) \xrightarrow{{'}e_\sharp} \mbH^1 (X, \mcK^\bullet [\nabla^\mr{ad}])  \xrightarrow{{''}e_\flat } H^0 (X, \mr{Coker}(\nabla^\mr{ad}))\hspace{10mm}\\
\left(\text{resp.,} \
{^c}e_\flat^{(0)} : H^0 (X, {^c}\DV_G) \xrightarrow{{'}e_{\sharp, c}} \mbH^1 (X, \mcK^\bullet [{^c}\nabla^\mr{ad}])  \xrightarrow{{''}e_{\flat, c} } H^0 (X, \mr{Coker}({^c}\nabla^\mr{ad})). 
 \right) \notag
\end{align}
The discussion preceding  the above proposition  (together with the equality (\ref{Eq211})) also shows that
the tangent space $T_q \mcO p_G^{^\mr{Zzz...}}$ (resp., $T_q \mcO p_{G, \rho}^{^\mr{Zzz...}}$)
of $\mcO p_G^{^\mr{Zzz...}}$ (resp., $\mcO p_{G, \rho}^{^\mr{Zzz...}}$) at $q$
may be identified with
 $\mr{Ker}(e_\flat^{(0)}) = \mr{Im}({'}e_\sharp) \cap \mr{Im}({''}e_\sharp)$ (resp., $\mr{Ker}({^c}e_\flat^{(0)}) = \mr{Im}({'}e_{\sharp, c}) \cap \mr{Im}({''}e_{\sharp, c})$).
 That is to say,  we obtain  an isomorphism of $k$-vector spaces
\begin{align} \label{Eq301}
{^p}\gamma_\sharp := \gamma_\sharp \cap {^p}\gamma  : T_q \mcO p_G^{^\mr{Zzz...}} \!\isom \mr{Ker}(e_\flat^{(0)})
\
\left(\text{resp.,} \  
{^p}\gamma_{\sharp, \rho} := \gamma_{\sharp, \widetilde{\rho}} \cap {^p}\gamma_\rho  : T_q \mcO p_{G, \rho}^{^\mr{Zzz...}}\! \isom \mr{Ker}({^c}e_\flat^{(0)})
\right)  \hspace{-3mm}
\end{align}
(cf. ~\cite[Corollary 6.20]{Wak8}).

\LSP
%---------------------------[begin subsection]-------------
\subsection{Ordinariness for dormant $G$-opers} \label{SS0101}

Next, we shall denote by 
\begin{align} \label{Eq100}
\mcO p^{^\mr{Zzz...}}_{G, \text{\'{e}t}}
\end{align}
the \'{e}tale locus  of $\mcO p^{^\mr{Zzz...}}_{G}$ over $k$.
It follows from ~\cite[Theorem G]{Wak8} that
$\mcO p^{^\mr{Zzz...}}_{G, \text{\'{e}t}} = \mcO p^{^\mr{Zzz...}}_{G}$ if $\msX$ is a general pointed stable curve and $G$ is either $\mr{PGL}_n$ with $2n < p$, $SO_{2l+1}$ with $4l +2 < p$, or $PCSp_{2m}$ with $4m< p$.
For example, if $G$ is one of these groups, then  this equality holds for 
every totally  degenerate  curve, in the sense of ~\cite[Definition 7.15]{Wak8}.

\SSP
%---------------------------------------------------------------------------------
\bde \label{Def3}
We shall say that a dormant $G$-oper $\msE^\spadesuit$ on $\msX$ is {\bf ordinary}
if it is classified by a point of $\mcO p^{^\mr{Zzz...}}_{G, \text{\'{e}t}}$, or equivalently, 
the morphism  $e_\flat^{(0)}$ associated to $\msE^\spadesuit$ is injective (cf. (\ref{Eq301})).
\ede
%---------------------------------------------------------------------------------
\SSP

%---------------------------------------------------------------------------------
\begin{rema}[Previous studies]
In some of the  papers written  by the author, the ordinariness for opers and  related notions were  discussed in several  settings.
For example, 
we refer the reader to  ~\cite[Definitions 2.1.2, 3.4.1]{Wak9} for the ordinariness of  $\mr{PGL}_n$-opers.
On the other hand,
the  case of dormant  indigenous $(G, H)$-bundles (where $G$ denotes a connected smooth algebraic group and $H$ denotes a suitable closed subgroup of $G$) on a smooth algebraic variety can be found in 
~\cite[Definition 6.7.1]{Wak11}.
\end{rema}
%---------------------------------------------------------------------------------
%\SSP

\SSP
%-------------------------------------------------------------
\bpr \label{Thm4}
Let $\msE^\spadesuit := (\mcE_B, \nabla)$ be a  dormant $\,\qq$-normal  $G$-oper on $\msX$.
Then, the  following three conditions are equivalent to each other:
\begin{itemize}
\item[(a)]
$\msE^\spadesuit$ is ordinary in the sense of Definition \ref{Def3}.
\item[(b)]
The morphism  $e_\flat^{(0)}$ associated to $\msE^\spadesuit$ is an isomorphism.
\item[(c)]
 The composite
\begin{align} \label{Eq876}
e_\sharp^{(0)} : H^1 (X, \mr{Ker}(\nabla^\mr{ad})) \xrightarrow{{''}e_\sharp} \mbH^1 (X, \mcK^\bullet [\nabla^\mr{ad}]) \xrightarrow{{'}e_\flat } H^1 (X, \Omega  \otimes_{\mcO_X} \DV_G^\vee)
\end{align}
is an isomorphism.
\end{itemize}
\epr
%-------------------------------------------------------------
\begin{proof}
The equivalence (b) $\Leftrightarrow$ (c) follows from the exactness of the sequences (\ref{HDeq1}) and  (\ref{Ceq}).
Also, there is nothing to prove for the implication (b) $\Rightarrow$ (a).
Finally, the implication  (a) $\Rightarrow$ (b) follows from the equality
 \begin{align}
 \mr{dim}H^0 (X, \DV_G) = \mr{dim}H^0 (X, \mr{Coker}(\nabla^\mr{ad})) \left( = \frac{2g-2 + r}{2} \cdot \mr{dim}( \mfg ) + \frac{r}{2} \cdot \mr{rk}(\mfg) \right),
 \end{align}
 that was already proved in ~\cite[Propositions 2.23 and 6.18]{Wak8} (see also Proposition \ref{Eqw3} described above).
\end{proof}
%-------------------------------------------------------------
\SSP

%\SSP
%-------------------------------------------------------------
\bco \label{COg467}
Let $\msE^\spadesuit$ be an ordinary dormant $G$-oper on $\msX$.
Then, there exists a canonical  decomposition of $\mbH^1 (X, \mcK^\bullet [\nabla^\mr{ad}])$ into a direct sum
\begin{align} \label{Eq110}
\mbH^1 (X, \mcK^\bullet [\nabla^\mr{ad}]) 
&=  H^0 (X, \DV_G) \oplus H^1 (X, \Omega \otimes_{\mcO_X} \DV_G^\vee)  
\left(\cong H^0 (X, \DV_G) \oplus H^0 (X, {^c}\DV_G)^\vee \right).  
\end{align}
\eco
%-------------------------------------------------------------
\begin{proof}
After possibly replacing $\msE^\spadesuit$ with its $\,\qq$-normalization (cf. ~\cite[Proposition 2.19]{Wak8}), we may assume that $\msE^\spadesuit$ is $\,\qq$-normal.
Since $e_\flat^{(0)}$ is bijective by the equivalence (a) $\Leftrightarrow$ (b) asserted in  Proposition  \ref{Thm4},
we obtain the composite surjection
\begin{align} \label{Eq308}
\mbH^1 (X, \mcK^\bullet [\nabla^\mr{ad}]) \xrightarrow{{''}e_\flat} H^1 (X, \mr{Coker}(\nabla^\mr{ad})) \xrightarrow{(e_\flat^{(0)})^{-1}} H^0 (X, \DV_G).
\end{align}
It specifies  a split surjection of the short exact sequence (\ref{HDeq1}) in the case where $\Box$ denotes the absence of ``$c$".
In particular, this split surjection determines the  desired decomposition (\ref{Eq110}).
\end{proof}
%-------------------------------------------------------------
\SSP

Let  $\msE^\spadesuit$ and $q$ be as in \S\,\ref{SS030}.
 Suppose further  that $\msE^\spadesuit$ is ordinary, which means that
the two substacks $\mcO p_G$, $\C_G^{\psi =0}$ of $\C_G$ intersect transversally at $q$.
Since $e_\flat^{(0)}$ is an isomorphism (cf. Proposition \ref{Thm4}), the morphism
\begin{align} \label{Eq304}
({'}e_\sharp, {''}e_\sharp) : H^0 (X, \DV_G) \oplus H^1 (X, \mr{Ker}(\nabla^\mr{ad})) \migi \mbH^1 (X, \mcK^\bullet [\nabla^\mr{ad}])
\end{align}
turns out to be  an isomorphism.
Under the identifications 
between various $k$-vector spaces  given by $\gamma$ (cf. (\ref{Eq3033})), $\gamma_\sharp$ (cf. (\ref{Eq300})),  and ${^p}\gamma$ (cf. (\ref{Eq303})),
the isomorphism (\ref{Eq304}) reads 
 a decomposition 
\begin{align} \label{Eq306}
T_q \mcO p_G \oplus T_q \C_G^{\psi = 0} \isom T_q \C_G
\end{align}
of the tangent space  $T_q \C_G$;  it can be  obtained  by differentiating 
the immersions  $\mcO p_G \migiincl \C_G$ and  $\C_G^{\psi = 0} \migiincl \C_G$ at $q$.

%%%%%%%%%%%%%%%%%%%%%%%%%%%%%%%%%%%%%%%%%%%%%
%%%%%%%%%%%%%---[begin section]---%%%%%%%%%%%%%%%%%%%%%%%
\vspace{10mm}
\section{Eichler-Shimura-type isomorphisms  for dormant opers}\SSP

This  final section is devoted to  proving  Theorem \ref{ThA}, which establishes  an Eichler-Shimura-type decomposition  for a  dormant $G$-oper
 on a general pointed stable curve.
This  is constructed  by relating it to the decomposition of a tangent space already obtained in  (\ref{Eq306}).

\LSP
%----------------------------------------------------------------------[begin subsection]-------------
\subsection{The irreducible decomposition of the deformation space} \label{SS0134}

Let $\msE^\spadesuit := (\mcE_B, \nabla)$ be a dormant $\,\qq$-normal  $G$-oper on $\msX$.
Suppose that $\msE^\spadesuit$ is of radii $\rho := (\rho_i)_{i=1}^r\in \mfc (k)^{\times r}$, where $\rho := \emptyset$ if $r = 0$.

Note that the dual  $(\mfg_\mcE^\vee, \nabla^{\mr{ad}\vee})$ of the flat bundle $(\mfg_\mcE, \nabla^\mr{ad})$ is isomorphic to $(\mfg_\mcE, \nabla^\mr{ad})$ itself because of the nondegeneracy of the Killing form on $\mfg$.
Hence, by applying  ~\cite[Corollary 6.16]{Wak8} to the case where ``$(\mcV, \nabla)$" is taken to be $(\mfg_\mcE, \nabla^\mr{ad})$, we obtain 
  an isomorphism of short exact sequences:
\begin{align} \label{Eq240}
\vcenter{\xymatrix@C=21pt@R=36pt{
0 \ar[r] &  H^1 (X, \mr{Ker}({^c}\nabla^\mr{ad})) \ar[r]^-{{''}e_{\sharp, c}} \ar[d]_-{\wr}^-{} &\mbH^1 (X, \mcK^\bullet [{^c}\nabla^\mr{ad}]) \ar[r]^-{{''}e_{\flat, c}} \ar[d]_-{\wr}^-{(\ref{Eww3})} & H^0 (X, \mr{Coker}({^c}\nabla^\mr{ad})) \ar[r] \ar[d]_-{\wr}^-{} & 0
\\
0 \ar[r]& H^0 (X, \mr{Coker}(\nabla^\mr{ad}))^{\vee} \ar[r]_-{({''}e_\flat)^\vee } &\mbH^1 (X, \mcK^\bullet [\nabla^\mr{ad}])^\vee  \ar[r]_-{({''}e_\sharp)^\vee} &H^1(X, \mr{Ker}(\nabla^\mr{ad}))^\vee   \ar[r]  & 0.
}}
\end{align}

\SSP
%-----------------------------------------------------------------------------------
\ble \label{Tgge}
The isomorphism $e_{\flat}^{(0)}$ associated to $\msE^\spadesuit$ (cf. Proposition \ref{Thm4}) restricts to an isomorphism  of $k$-vector spaces
\begin{align}
%{^c}e_{\flat}^{(0)} : 
H^0 (X, {^c}\DV_G) \isom H^0 (X, \mr{Coker}({^c}\nabla^\mr{ad})).
\end{align}
\ele
%-----------------------------------------------------------------------------------
\begin{proof}
Note that both the log connections $\nabla^\mr{ad}$ and ${^c}\nabla^\mr{ad}$ have  vanishing $p$-curvature.
Hence,
the morphism  $\delta : \mr{Coker}({^c}\nabla^\mr{ad}) \migi \mr{Coker}(\nabla^\mr{ad})$
 induced by the inclusion of complexes $\mcK^\bullet [{^c}\nabla^\mr{ad}] \migiincl \mcK^\bullet [\nabla^\mr{ad}]$ is injective because of the commutativity of  the following square diagram:
\begin{align} \label{Eq24011}
\vcenter{\xymatrix@C=46pt@R=36pt{
\mr{Coker}({^c}\nabla^\mr{ad}) \ar[r]^-{\delta} \ar[d]_-{\wr} &
\mr{Coker}(\nabla^\mr{ad}) \ar[d]^-{\wr}
\\
\Omega^{(1)} \otimes_{\mcO_{X^{(1)}}} \mr{Ker}({^c}\nabla^\mr{ad})
\ar[r]_-{\mr{inclusion}} &
\Omega^{(1)} \otimes_{\mcO_{X^{(1)}}} \mr{Ker}(\nabla^\mr{ad}),
}}
\end{align}
where
\begin{itemize}
\item
$X^{(1)}$ denotes the Frobenius twist of $X$ over $k$, and $\Omega^{(1)}$ denotes the  pull-back of $\Omega$ along the base-change morphism $X^{(1)} \isom X$ induced by the Frobenius automorphism of $k$;
\item
we regard both $\mr{Ker}({^c}\nabla^\mr{ad})$ and  $\mr{Ker}(\nabla^\mr{ad})$ as $\mcO_{X^{(1)}}$-modules via the underlying homeomorphism of the relative Frobenius morphism $X \migi X^{(1)}$;
\item
the right-hand and left-hand vertical arrows denote the isomorphisms induced by the Cartier operators of 
 $(\mfg_\mcE, \nabla^\mr{ad})$ and $({^c}\mfg_\mcE, {^c}\nabla^\mr{ad})$, respectively
(cf. the discussion following ~\cite[Proposition 1.2.4]{Og}).
\end{itemize}
In particular, 
the morphism 
\begin{align}
H^0(\delta) : H^0 (X, \mr{Coker}({^c}\nabla^\mr{ad})) \migi  H^0 (X, \mr{Coker}(\nabla^\mr{ad}))
\end{align}
 induced by $\delta$ is injective.
It follows that the composite of natural morphisms ${^c}\DV_G [-1] \migi \mcK^\bullet [{^c}\nabla^\mr{ad}] 
\migi \mr{Coker}({^c}\nabla^\mr{ad})[-1]$
 gives a restriction 
\begin{align} \label{Fgr4}
%{^c}e_{\flat}^{(0)} : 
H^0 (X, {^c}\DV_G) \migiincl  H^0 (X, \mr{Coker}({^c}\nabla^\mr{ad}))
\end{align}
of $e_{\flat}^{(0)}$ via the injection $H^0 (\delta)$ and the natural inclusion $H^0 (X, {^c}\DV_G) \migiincl H^0 (X, \DV_G)$.
Since the dual  $(\mfg_\mcE^\vee, \nabla^{\mr{ad}\vee})$ of the flat bundle $(\mfg_\mcE, \nabla^\mr{ad})$ is isomorphic to $(\mfg_\mcE, \nabla^\mr{ad})$ itself because of the nondegeneracy of the Killing form on $\mfg$,
it follows from ~\cite[Corollary 6.16]{Wak8}
that $H^0 (X, \mr{Coker}({^c}\nabla^\mr{ad}))$ is isomorphic to $H^1 (X, \mr{Ker}(\nabla^\mr{ad}))^\vee$.
Hence, 
this fact together with
~\cite[Propositions 2.23 and 6.18]{Wak8} implies the equality 
\begin{align} \label{Eq325}
\mr{dim}( H^0 (X, {^c}\DV_G)) =  \mr{dim} (H^0 (X, \mr{Coker}({^c}\nabla^\mr{ad}))) \left(= \frac{2g-2 +r}{2} \cdot \mr{dim}(\mfg) - \frac{r}{2} \cdot \mr{rk}(\mfg) \right).
\end{align}
Thus,  (\ref{Fgr4}) is verified to be an isomorphism, as desired. 
\end{proof}
%-----------------------------------------------------------------------------------
\SSP

Next, suppose that $\msE^\spadesuit$ is ordinary.
The dual of the composite  
\begin{align}
H^0 (X, \mr{Coker}(\nabla^\mr{ad})) \xrightarrow{(e_\flat^{(0)})^{-1}} H^0(X, \DV_G) \xrightarrow{{'}e_\sharp} \mbH^1 (X, \mcK^\bullet [\nabla^\mr{ad}])
\end{align}
 determines a split surjection  of the lower horizontal sequence in (\ref{Eq240}).
Hence, it induces  a splitting of the upper horizontal one.
Under the identifications between various $k$-vector spaces given by $\gamma_{\widetilde{\rho}}$ (cf. (\ref{Eq350})), ${^p}\gamma_{\widetilde{\rho}}$ (cf. (\ref{Eq303})), and ${^c}e_\flat^{(0)}\circ \gamma_{\sharp, \rho}$ (cf. (\ref{Eq300}), (\ref{Eq302})), the resulting decomposition of $\mbH^1 (X, \mcK^\bullet [{^c}\nabla^\mr{ad}])$ determines a direct sum decomposition
\begin{align} \label{Ek8}
T_q \C_{G, \widetilde{\rho}}^{\psi = 0} \oplus T_q \mcO p_{G, \rho}    \isom   T_q \C_{G, \widetilde{\rho}}
\end{align}
of $T_q \C_{G, \widetilde{\rho}}$; it 
 coincides with the morphism induced by differentiating  the inclusions $\mcO p_{G, \rho} \migiincl  \C_{G, \widetilde{\rho}}$  and  $\C_{G, \widetilde{\rho}}^{\psi = 0} \migiincl  \C_{G, \widetilde{\rho}}$.

Then, the following assertion holds.

\SSP
%-----------------------------------------------------------------------------------
\bt \label{Th77}
Let $\msE^\spadesuit$ be an ordinary dormant $\,\qq$-normal  $G$-oper on $\msX$.
Let us consider  the following two composite isomorphisms:
\begin{align} \label{Ek1}
T_q \C_{G, \widetilde{\rho}}^{\psi = 0}
& \xrightarrow{{^p}\gamma_{\widetilde{\rho}}} H^1 (X, \mr{Ker}({^c}\nabla^\mr{ad}))   \xrightarrow{\sim} H^0(X, \mr{Coker}(\nabla^\mr{ad}))^\vee  \\
& \xrightarrow{(e_\flat^{(0)})^\vee} H^0 (X, \DV_G)^\vee  
\xrightarrow{\gamma_\sharp^\vee} T^\vee_q \mcO p_G, \notag
\end{align}
\begin{align} \label{Ek2}
T_q \mcO p_{G, \rho} &\xrightarrow{\gamma_{\sharp, \rho}}H^0 (X, {^c}\DV_G) \xrightarrow{{^c}e_\flat^{(0)}}H^0 (X, \mr{Coker}({^c}\nabla^\mr{ad}))  \\ & \xrightarrow{\sim} H^1 (X, \mr{Ker}(\nabla^\mr{ad}))^\vee \xrightarrow{({^p}\gamma)^\vee} T_q^\vee \C_{G}^{\psi = 0} \notag
\end{align}
(cf. (\ref{Eq240}) for the definitions of the second arrow in (\ref{Ek1}) and the third arrow in (\ref{Ek2})).
Then,  these isomorphisms make the following square diagram commute:
 \begin{align} \label{Eq223}
\vcenter{\xymatrix@C=46pt@R=36pt{
 T_q \C_{G, \widetilde{\rho}}^{\psi = 0} \oplus T_q \mcO p_{G, \rho} 
 \ar[r]^-{(\ref{Ek8})}_-{\sim} \ar[d]_-{(\ref{Ek1}) \oplus (\ref{Ek2})}^-{\wr} & 
 T_q \C_{G, \widetilde{\rho}}
  \ar[d]^-{(\ref{Er67})}_-{\wr}
 \\
 T_q^\vee \mcO p_{G} \oplus   T_q^\vee \C_{G}^{\psi = 0}  
  \ar[r]_-{(\ref{Eq306})^\vee}^-{\sim} & T_q^\vee \C_{G}.
}}
\end{align}
In particular, if $r = 0$, then
the subspace $T_q \C_G^{\psi = 0}$ of $T_q \C_G$ is  Lagrangian with respect to the bilinear form on $T_q \C_G$ determined by (\ref{Er67}).
\et
%-----------------------------------------------------------------------------------
\begin{proof}
The assertion follows  from  Theorem \ref{Prc1}, (ii), and  the definitions of various morphisms involved.
\end{proof}
%-----------------------------------------------------------------------------------

\SSP

%------------------------------------------------------------------------------
\begin{rema}[Case of the parabolic de Rham cohomology] \label{Rem88}
When $r > 0$,
we cannot  construct a direct sum decomposition of $T_q \C_{G, [0]^{\times r}}$  using a dormant $G$-oper  in the same manner as (\ref{Eq306}) or (\ref{Ek8}).
In fact, as  discussed  in ~\cite[Remark 3.37]{Wak8},
any dormant $G$-oper cannot have radii $[0]^{\times r}$.
\end{rema}
%------------------------------------------------------------------------------

\LSP
%----------------------------------------------------------------------[begin subsection]-------------
\subsection{The canonical decomposition derived  from ordinariness} \label{SS0893}

Next, let $\msE^\spadesuit_\odot := (\mcE_{B^\odot}, \nabla_\odot)$ be a dormant $G^\odot$-oper on $\msX$,
and write $\mcE_\odot := \mcE_{B^\odot} \times^{B^\odot} G^\odot$ and $\msE_\odot := (\mcE_\odot, \nabla_\odot)$.

%----------------------------------[begin lemma]------------------
\SSP
\bpr\label{T01}
The following two conditions are equivalent to each other:
\begin{itemize}
\item[(a)]
The dormant $G$-oper $\iota_{G *}(\msE^\spadesuit_\odot)$ (cf. (\ref{associated})) is ordinary;
\item[(b)]
 For each integer  $l$ with $1\leq l \leq \mr{rk}(\mfg)$, 
the composite  morphism
\begin{align} \label{Eq838}
\kappa_{l} : H^1 (X, \mr{Ker} (\mr{Sym}^{2l}\nabla_\odot)) 
\xrightarrow{\theta_l }
 H^1_{\mr{dR}} (X, \mr{Sym}^{2l}\msE_\odot) \longmigi
  %\xrightarrow{h^l_\flat [\msE^\spadesuit_\odot]}
   H^1 (X, \Omega^{\otimes (-l)})
\end{align}
(cf. (\ref{Eq817}) for the definition of the second arrow) is an isomorphism, where the first arrow $\theta_l$
 %${''}e_\sharp [\mr{Sym}^{2l}\nabla_\odot]$ 
 arises from the natural inclusion $\mr{Ker} (\mr{Sym}^{2l}\nabla_\odot) [0] \migi \mcK^\bullet [\mr{Sym}^{2l}\nabla_\odot]$.
\end{itemize}
 \epr
%------------------------------[begin proof]-------------------
\begin{proof}
After possibly replacing $\msE^\spadesuit_\odot$ with its $\,\qq$-normalization,
we may assume that $\msE^\spadesuit_\odot$ is $\,\qq$-normal.
We shall  write $\iota_{G*}(\msE^\spadesuit_\odot) := (\mcE_B, \nabla)$.
The isomorphism (\ref{Eq160}) restricts to an isomorphism
\begin{align}
\bigoplus_{l=1}^{\mr{rk}(\mfg)}\mr{Ker}(\mr{Sym}^{2l}\nabla_\odot)  \otimes_k \mfg_l^{\mr{ad}(q_1)} \isom \mr{Ker}(\nabla^\mr{ad}).
\end{align}
It  induces 
a commutative square diagram
\begin{align} \label{Eq166}
\vcenter{\xymatrix@C=26pt@R=36pt{
\displaystyle{\bigoplus_{l=1}^{\mr{rk}(\mfg)}}H^1 (X, \mr{Ker}(\mr{Sym}^{2l}\nabla_\odot))  \otimes_k \mfg_l^{\mr{ad}(q_1)}
\ar[r]^-{\sim} \ar[d]_-{\bigoplus_l \theta_l \otimes \mr{id}} & H^1 (X,  \mr{Ker}(\nabla^\mr{ad})) \ar[d]^-{{''}e_\sharp }
\\
\displaystyle{\bigoplus_{l=1}^{\mr{rk}(\mfg)}} H^1_{\mr{dR}}(X, \mr{Sym}^{2l}\msE_\odot) \otimes_k\mfg_l^{\mr{ad}(q_1)}\ar[r]_-{(\ref{Eq834})} & \mbH^1_{} (X, \mcK^\bullet [\nabla^\mr{ad}]).
}}
\end{align}
By composing (\ref{Eq166}) and the right-hand square diagram in  (\ref{Eq16tt}),
we obtain a commutative square
\begin{align} \label{Eq16ji}
\vcenter{\xymatrix@C=26pt@R=36pt{
\displaystyle{\bigoplus_{l=1}^{\mr{rk}(\mfg)}}H^1 (X, \mr{Ker}(\mr{Sym}^{2l}(\nabla)))  \otimes_k \mfg_l^{\mr{ad}(q_1)}
\ar[r]^-{\sim} \ar[d]_-{\bigoplus_{l}\kappa_{l} \otimes \mr{id}} & H^1 (X,  \mr{Ker}(\nabla^\mr{ad})) \ar[d]^-{e_\sharp^{(0)}}
\\
\displaystyle{\bigoplus_{l=1}^{\mr{rk}(\mfg)}}H^1 (X, \Omega^{\otimes (-l)})\otimes_k\mfg_l^{\mr{ad}(q_1)}\ar[r]_-{\sim} & H^1 (X, \Omega\otimes \DV_G^\vee)
}}
\end{align}
(cf. (\ref{Eq876}) for the definition of $e_\sharp^{(0)}$).
The commutativity of this diagram implies that $\kappa_{l}$ is an isomorphism for every $l$ if and only if $e_\sharp^{(0)}$
 is an isomorphism.
 Hence, 
 the assertion  follows  from the equivalence (a)  $\Leftrightarrow$ (c) asserted in Proposition \ref{Thm4}.
 \end{proof}
\SSP
%-------------------------------------------------------------

%\SSP
%-------------------------------------------------------------
\bco \label{Co34}
Let $\msE_\odot^\spadesuit$ be a dormant $G^\odot$-oper on $\msX$.
Also, let $G'$ be another algebraic group  
satisfying the same conditions imposed on $G$ (cf. \S\,\ref{Ssde}).
Suppose that $\mr{rk}(\mfg) \geq  \mr{rk}(\mfg')$ and that the $G$-oper $\iota_{G*}(\msE^\spadesuit)$ is ordinary.
Then, the  dormant $G'$-oper $\iota_{G'*}(\msE^\spadesuit)$ is ordinary.
\eco
%-------------------------------------------------------------
\begin{proof}
The assertion can be  immediately proved  by  applying Proposition \ref{T01} to both
$\iota_{G*}(\msE^\spadesuit)$ and  $\iota_{G'*}(\msE^\spadesuit)$.
\end{proof}
%-------------------------------------------------------------
\SSP

%-------------------------------------------------------------
\bco \label{Co34f}
Let $\msE_\odot^\spadesuit$ be a dormant $G^\odot$-oper on $\msX$.
Suppose that $\msX$ is general in the moduli stack $\overline{\mcM}_{g, r}$ of $r$-pointed stable curves of genus $g$ over $k$.
(To be precise, $\msX$ specifies a point of $\overline{\mcM}_{g, r}$ that lies outside some fixed closed substack.)
\begin{itemize}
\item[(i)]
Suppose that $\mr{rk}(\mfg) \leq \frac{p-3}{2}$.
Then,
the dormant $G$-oper $\iota_{G *}(\msE^\spadesuit_\odot)$ associated to  any dormant $G^\odot$-oper $\msE^\spadesuit_\odot$ on $\msX$  is ordinary.
\item[(ii)]
For each integer $l$ with   $1 \leq l \leq \frac{p-3}{2}$,
the morphism $\kappa_{l}$ associated to $\msE^\spadesuit_\odot$ (introduced in Proposition \ref{T01}) is an isomorphism.
\end{itemize}
\eco
%-------------------------------------------------------------
\begin{proof}
First, let us  prove  assertion (i).
For simplicity, we shall set  $G^\circledast := \mr{PGL}_\frac{p-1}{2}$, which has rank $\frac{p-3}{2}$.
By applying  Corollary \ref{Co34} to both $G$ and $G^\circledast$,
we can reduce 
 the problem to the case of $G = G^\circledast$.
Denote by $\mfO \mfp^{^\mr{Zzz...}}_{G^\odot}$ (resp., $\mfO \mfp^{^\mr{Zzz...}}_{G^\circledast}$) the Deligne-Mumford stack  classifying pairs $(\msX, \msE^\spadesuit)$ consisting of $\msX \in \mr{Ob}(\overline{\mcM}_{g, r})$ and a dormant $G^\odot$-oper  (resp., a dormant $G^\circledast$-oper) $\msE^\spadesuit$ on $\msX$.
According to ~\cite[Chap.\,II, Theorem 2.8]{Mzk2}, the stack   $\mfO \mfp^{^\mr{Zzz...}}_{G^\odot}$ is   irreducible, and the natural projection $\mfO \mfp^{^\mr{Zzz...}}_{G^\odot} \migi \overline{\mcM}_{g, r}$ is finite and faithfully flat.
Hence, one can find  a unique  irreducible component  $\mcN$ of  $\mfO \mfp_{G^\circledast}^{^\mr{Zzz...}}$ containing the image of the $\overline{\mcM}_{g, r}$-morphism 
 \begin{align} \label{Mor44}
 %\iota_{G^\circledast}^{\mfO \mfp} : 
 \mfO \mfp_{G^\odot}^{^\mr{Zzz...}} \migi \mfO \mfp_{G^\circledast}^{^\mr{Zzz...}}
 \end{align}
    given by assigning $(\msX, \msE^\spadesuit) \mapsto (\msX, \iota_{G*}(\msE^\spadesuit))$.
If $\mfO \mfp_{G^\circledast, \text{\'{e}t}}^{^\mr{Zzz...}}$ denotes the \'{e}tale locus of $\mfO \mfp_{G^\circledast}^{^\mr{Zzz...}}$ over $\overline{\mcM}_{g, r}$, then
the intersection  $\mfO \mfp_{G^\circledast, \text{\'{e}t}}^{^\mr{Zzz...}} \cap \mcN$ defines a dense open substack of $\mcN$ (cf.  ~\cite[Theorem G]{Wak8}).
Since  $\mfO \mfp_{G^\circledast}^{^\mr{Zzz...}}$ is finite over $\overline{\mcM}_{g, r}$ (cf. ~\cite[Theorem 3.34]{Wak8}), the restriction 
%\begin{align}
$\mfO \mfp_{G^\odot}^{^\mr{Zzz...}} \migi \mcN$ of 
(\ref{Mor44})
%$\iota_{G^\circledast}^{\mfO \mfp}$
 is dominant.
This implies
 that the inverse image  
 $\mcQ$ of  $\mfO \mfp_{G^\circledast, \text{\'{e}t}}^{^\mr{Zzz...}} \cap \mcN$ via (\ref{Mor44})
 %$\mcQ :=  (\iota_{G^\circledast}^{\mfO \mfp})^{-1} (\mfO \mfp_{G^\circledast, \text{\'{e}t}}^{^\mr{Zzz...}} \cap \mcN)$
  is a dense open substack of $\mfO \mfp_{G^\odot}^{^\mr{Zzz...}}$.
The complement in $\overline{\mcM}_{g, r}$ of the image of $\mfO \mfp_{G^\odot}^{^\mr{Zzz...}}  \setminus \mcQ$ forms a dense open substack of $\overline{\mcM}_{g, r}$ and  by definition classifies pointed stable curves $\msX$
such that $\iota_{G^\circledast *}(\msE^\spadesuit)$ is ordinary for any dormant $G^\odot$-oper $\msE^\spadesuit$ on $\msX$.
This completes the proof of assertion (i).

Moreover, assertion (ii) follows from assertion (i) together with Proposition \ref{T01}.
\end{proof}
%-------------------------------------------------------------
\SSP

We conclude this paper  with the following assertions.

\SSP
%-------------------------------------------------------------------------------------------------------
\bt \label{Co12}
Let $\msE^\spadesuit_\odot$ be a dormant $G^\odot$-oper on $\msX$, and write
$\msE_\odot$ for the underlying flat $G^\odot$-bundle of $\msE^\spadesuit_\odot$.
Also, let $l$ be an integer with $1 \leq l \leq \frac{p-3}{2}$.
Suppose that  $\msX$ is general in $\overline{\mcM}_{g, r}$.
Then, 
 there exists a canonical 
 splitting of (\ref{Eq817}).
 In particular, we obtain a canonical 
 decomposition
\begin{align}
H^1_{\mr{dR}} (X, \mr{Sym}^{2l}\msE_\odot) &= H^0 (X, \Omega^{\otimes (l+1)}) \oplus H^1 (X, \Omega^{\otimes (-l)}) \\
& \hspace{-1mm}\left(\cong  H^0 (X, \Omega^{\otimes (l+1)}) \oplus H^0 (X, \Omega^{\otimes (l+1)}(-D))^\vee\right) \notag
\end{align}
of the $1$-st de Rham cohomology group $H^1_{\mr{dR}} (X, \mr{Sym}^{2l}\msE_\odot)$ of $\mr{Sym}^{2l}\msE_\odot$.
\et
%-------------------------------------------------------------------------------------------------------
\begin{proof}
It follows from Corollary \ref{Co34f}, (ii),  that 
$\kappa_{l}$ is an isomorphism for every  $l = 1, \cdots, \frac{p-3}{2}$.
Then,  
the composite injection 
\begin{align} 
H^1 (X, \Omega^{\otimes (-l)}) \xrightarrow{\kappa_{l}^{-1}} H^1 (X, \mr{Ker}(\mr{Sym}^{2l}\nabla_\odot)) \xrightarrow{\theta_l}H_\mr{dR}^1 (X, \mr{Sym}^{2l}\msE_\odot).
\end{align}
specifies 
the desired splitting. 
\end{proof}
%-------------------------------------------------------------------------------------------------------
\SSP

%\SSP
%---------------------------------------------------------------------[begin theorem]----------------------
\bco[cf. Theorem \ref{ThA}]\label{ThAgg}  
Let $\msF^\heartsuit_\odot$ be   a  dormant $\mr{SL}_2$-oper on $\msX$, and  $l$  an integer with $1 \leq l \leq \frac{p-3}{2}$.
Denote by $\varTheta$ the theta characteristic of $\msX$ associated to $\msF^\heartsuit_\odot$, and 
by  $\msF_\odot$ the underlying flat $\mr{SL}_2$-bundle of $\msF^\heartsuit_\odot$.
Suppose that $\msX$ is general in $\overline{\mcM}_{g, r}$.
Then, there exists a canonical splitting  of (\ref{Eq810}) for $n= 2l +1$.
In particular,  we obtain  
a canonical decomposition 
\begin{align} \label{Eq131}
H^1_{\mr{dR}} (X, \mr{Sym}^{2l}\msF) & =   H^0  (X, \varTheta^{\otimes 2(l+1)}) \oplus H^1 (X, \varTheta^{\otimes (-2l)}) \\
& \hspace{-1mm} \left(\cong H^0  (X, \varTheta^{\otimes 2(l+1)}) \oplus H^0 (X, \varTheta^{\otimes 2(l+1)}(-D))^\vee  \right) \notag
\end{align}
of  the $1$-st de Rham cohomology $H^1_{\mr{dR}} (X, \mr{Sym}^{2l}\msF_\odot)$ of $\mr{Sym}^{2l}\msF_\odot$.
 \eco
%------------------------------------------------------------------------------[end theorem]-------------
\begin{proof}
The assertion follows from Theorem \ref{Co12} and the comment at the end of  the statement in Theorem \ref{Th66}.
\end{proof}
%--------------------------------------------------------------------------

\SSP
%---------------------------------[begin remark]------------------
\begin{rema} \label{R037}
In ~\cite[\S\,2.4]{Wak7}, the author constructed a $p$-adic version of the Eichler-Shimura isomorphism for the $2$-nd symmetric product of  an ordinary nilpotent indigenous bundle, i.e., a certain $G^\odot$-oper with nilpotent $p$-curvature (cf. ~\cite[Chap.\,II,  Definition 3.1]{Mzk1}).
This result  is an essential ingredient for proving  the main theorem in ~\cite{Wak7}, which  compares  canonical symplectic structures on the related moduli spaces.
Despite the similarity in concept, this result seems to have no (at least direct) relation with the decomposition asserted in the above theorem because the $G^\odot$-opers treated there never have vanishing $p$-curvature.
However, we also expect that the previous  $p$-adic version extends to  the higher symmetric products,  generalizing Faltings' construction of the $p$-adic Eichler-Shimura isomorphisms for modular curves (cf. ~\cite[Theorem 6]{Fal}).
 \end{rema}
%-----------------------------------[end remark]-------------------

\LSP
%%%%%%%%%%%%%%%%%%%%%%%%%%%%%%%---[begin section]---%%%%%%%%%%%%%%
\subsection*{Acknowledgements} 
We would like to thank modular curves for their helpful comments on 
the Eichler-Shimura isomorphism.
Our work was partially supported by Grant-in-Aid for Scientific Research (KAKENHI No. 21K13770).

%%%%%%%%%%%%%%%%%%%%%%%%%%%%%%%%%%%%%%%%%%%%%
\end{document}